\documentclass[11pt,a4paper]{article}

\usepackage[utf8]{inputenc}
\usepackage[T1]{fontenc}
\usepackage[margin=1in]{geometry}
\usepackage{amsmath, amssymb, amsthm}
\usepackage{mathtools}
\usepackage[numbers,sort&compress]{natbib}
\usepackage{hyperref}
\usepackage{xcolor}
\usepackage{listings}
\usepackage{booktabs}
\usepackage{enumitem}
\usepackage{graphicx}
\usepackage{fullpage}
\usepackage{caption}
\usepackage{float}
\usepackage[protrusion=true,expansion=false]{microtype}

\definecolor{codebg}{rgb}{0.97,0.97,0.97}
\definecolor{codekw}{rgb}{0.12,0.32,0.62}
\definecolor{codestr}{rgb}{0.16,0.50,0.16}
\definecolor{codecomm}{rgb}{0.45,0.45,0.45}
\definecolor{codenum}{rgb}{0.55,0.55,0.55}

\lstdefinestyle{python}{
  language=Python,
  basicstyle=\ttfamily\footnotesize,
  backgroundcolor=\color{codebg},
  keywordstyle=\color{codekw}\bfseries,
  stringstyle=\color{codestr},
  commentstyle=\color{codecomm}\itshape,
  numberstyle=\tiny\color{codenum},
  numbers=left,
  numbersep=8pt,
  frame=single,
  rulecolor=\color{codebg},
  showstringspaces=false,
  breaklines=true,
  columns=fullflexible,
  keepspaces=true,
  xleftmargin=12pt,
  framexleftmargin=10pt,
  framexrightmargin=0pt,
  literate={->}{{$\rightarrow$}}1 {>=}{{$\geq$}}1 {<=}{{$\leq$}}1
           {⊤}{{$\top$}}1 {×}{{$\times$}}1 {→}{{$\rightarrow$}}1,
  extendedchars=true
}
\newtheorem{definition}{Definition}[section]
\newtheorem{theorem}[definition]{Theorem}
\newtheorem{proposition}[definition]{Proposition}

\newtheorem{example}[definition]{Example}
\newtheorem{remark}[definition]{Remark}

\newcommand{\Rplus}{\mathbb{R}_{\geq 0}}
\newcommand{\Rplustop}{\overline{\Rplus}}
\newcommand{\Nat}{\mathbb{N}}
\newcommand{\Min}{\mathsf{Min}}
\newcommand{\Antichain}{\mathsf{A}}
\newcommand{\topp}{\top}
\newcommand{\etal}{\textit{et al.}}
\newcommand{\codename}{\texttt{codesign-mcdp}}

\hypersetup{
  colorlinks=true,
  linkcolor=blue!60!black,
  citecolor=blue!60!black,
  urlcolor=blue!60!black,
}

\title{\codename: A Python Library for \\ Monotone Co-Design Problems \\
       \vspace{0.4em} \large Version 0.2.1 -- Reference Manual}
\author{Corentin Briat\\ School of Life Sciences\\
University of Life Sciences, Muttenz, Switzerland\\
\texttt{corentin@briat.info}, \url{www.briat.info}}
\date{July 2026}

\begin{document}
\maketitle

\begin{abstract}
\noindent \codename{} is a Python library for formulating and solving
\emph{Monotone Co-Design Problems} (MCDPs) in the framework of Censi
\citep{censi2015mathematical}. A design problem is a relation between two posets, a functionality
poset $F$ and a resource poset $R$; given a target functionality, the
problem asks for the antichain of minimal resources needed to deliver it.
Design problems compose under three operators (series, parallel,
feedback), and the resulting class is closed under composition. The
library implements the antichain calculus, six primitive design-problem
types, the three composition operators, a Kleene fixed-point solver, and
two high-level builders (an MCDPL-style declarative builder and a modular
\texttt{System} builder). Further layers add set-based and stochastic
uncertainty, compositional online learning, and temporal, vector-state,
and online co-design, alongside a suite of worked examples. This document
is the reference manual for version 0.2.1.
\end{abstract}

\tableofcontents
\bigskip

\section{Introduction}

A \emph{co-design problem} is the joint design of components whose
specifications mutually depend on each other. A canonical example: a
battery must store enough energy to power an actuator, but the actuator
must lift the battery itself, so the battery's mass appears in the
energy it has to supply. There is no way to fix one without fixing the
other.

\citet{censi2015mathematical} showed that under modest monotonicity
assumptions, such problems have a clean algebraic theory: the relations
between resources and functionalities form a category closed under
series, parallel, and feedback composition, and solutions are computed
by a Kleene fixed-point iteration in the lattice of antichains. The
original publication included a software tool, \texttt{MCDPL}
\citep{mcdpl}, distributed through \texttt{co-design.science}, which
exposes the framework via a domain-specific language.

\codename{} is a from-scratch Python implementation of the algorithmic
core, designed around composable Python objects rather than a separate
DSL. It targets users who want to formulate and solve MCDPs from inside
their existing Python codebase, without context switching to a DSL or an
external solver. The algorithmic core (posets, antichains, primitive
design problems, composition, solver, and both builders) has no required
runtime dependencies beyond the standard library; the uncertainty,
online-learning, and visualisation layers pull in \texttt{numpy},
\texttt{scipy}, and \texttt{matplotlib} as optional extras. The whole
package is roughly 8000 lines of Python and is licensed under MIT.

\paragraph{Attribution.} The theory implemented here is not the
author's. Monotone co-design was introduced by Andrea Censi
\citep{censi2015mathematical,censi2016everything}, who also established
the treatment of cyclic constraints \citep{censi2017cyclic} and of
uncertainty \citep{censi2017uncertainty}, and who wrote the reference
implementation \texttt{MCDPL} \citep{mcdpl}. The framework has been
developed since by Gioele Zardini, Dejan Milojevic, and colleagues,
across future mobility systems
\citep{zardini2020towards,zardini2020avmobility,zardini2023tools},
autonomous systems and embodied intelligence
\citep{zardini2021autonomous,zardini2021embodied,zardini2022taskdriven},
and the co-design of perception and compute
\citep{milojevic2025codei}, with more recent work on composable and
distributional uncertainty
\citep{huang2025composable,furter2026symmetric,huang2026distributional},
categorical and scalability extensions
\citep{riess2026quantale,cai2026scalable}, and online learning
\citep{alharbi2026online}. \codename{} is an independent Python
implementation of that framework. It contributes software, not theory.
Every definition, theorem, and algorithm in
Section~\ref{sec:theory}, and the uncertainty and online-learning
layers, belong to the authors cited above. Users of this library should
cite their papers, not only this manual. Citing the software is not a
substitute for citing the theory. Please cite \citet{censi2015mathematical}
for the monotone co-design framework itself, and the paper behind whichever
layer you use: \citet{censi2017uncertainty} for uncertainty,
\citet{alharbi2026online} for the online-learning layer, and the relevant
application paper from the co-design group
\citep{zardini2020avmobility,milojevic2025codei}.

The temporal, sequential, vector-state, and receding-horizon layers of
Section~\ref{sec:temporal} are the exception. Those are the author's own
work, built on top of the monotone co-design framework rather than taken
from it. They are being written up separately and are not yet published,
so they carry no citation here. Section~\ref{sec:related} maps the
literature and marks the boundary again.

\paragraph{Scope of this document.} This manual covers the public API
and the mathematical structure it implements. Section~\ref{sec:theory}
recalls the formal background. Sections~\ref{sec:posets} through
\ref{sec:solver} describe the data types, primitive design problems,
composition operators, and the solver in detail. Section~\ref{sec:uncertainty}
covers set-based and stochastic uncertainty, and
Section~\ref{sec:online-learning} the online-learning layer;
Section~\ref{sec:visualisation} documents the visualisation helpers.
Sections~\ref{sec:mcdp-builder} and~\ref{sec:system-builder} document the
two high-level builders. Section~\ref{sec:apiref} is a complete,
per-module reference for every public symbol in the library.
Section~\ref{sec:examples} walks through every
example shipped with the library. Section~\ref{sec:temporal} introduces
the temporal, vector-state, and online co-design layers, and
Section~\ref{sec:domains} surveys realistic problem domains for them.
Section~\ref{sec:guidelines} collects modelling guidelines and common
pitfalls, and Section~\ref{sec:limitations} discusses limitations and
future directions.

\paragraph{Conventions.} Code listings are typeset in monospaced font.
Type signatures use Python's standard notation
($\texttt{f: F} \rightarrow \texttt{R}$, for example). Mathematical
quantities use $F$, $R$, $h$, $\top$, $\bot$, $\Min$, etc.

\section{Mathematical Background}
\label{sec:theory}

This section recalls the formal definitions from
\citet{censi2015mathematical} needed to specify the API. The definitions,
the theorems, and the fixed-point algorithm are all his; nothing in this
section is original to this library. Proofs are omitted; the reader is
referred to the paper for the full development. Standard order-theoretic
background is in \citet{davey2002lattices}. The treatment of design
problems whose constraint graph contains cycles is developed further in
\citet{censi2017cyclic}.

\subsection{Posets and Lattices}

\begin{definition}[Poset]
A \emph{poset} (partially ordered set) is a pair $(P, \leq)$ where
$\leq$ is a reflexive, antisymmetric, transitive binary relation on
the set $P$.
\end{definition}

\begin{definition}[Bottom and top]
An element $\bot \in P$ is a \emph{bottom} (least element) if
$\bot \leq x$ for every $x \in P$. An element $\topp \in P$ is a
\emph{top} (greatest element) if $x \leq \topp$ for every $x \in P$.
A poset has at most one bottom and at most one top.
\end{definition}

\begin{example}[Reals and Naturals]
The non-negative reals augmented with infinity, $\Rplustop = \Rplus
\cup \{+\infty\}$, ordered by the usual $\leq$, form a poset with
$\bot = 0$ and $\topp = +\infty$. Likewise $\overline{\Nat} =
\Nat \cup \{+\infty\}$.
\end{example}

\begin{example}[Product order]
Given posets $(P_1, \leq_1), \dots, (P_n, \leq_n)$, the product
$P_1 \times \cdots \times P_n$ with the componentwise order
\[
(x_1, \dots, x_n) \leq (y_1, \dots, y_n)
  \iff \forall i.\, x_i \leq_i y_i
\]
is a poset. The bottom is $(\bot_1, \dots, \bot_n)$ and the top is
$(\topp_1, \dots, \topp_n)$.
\end{example}

\begin{definition}[Antichain]
A subset $A \subseteq P$ is an \emph{antichain} if no two distinct
elements of $A$ are comparable: $\forall x, y \in A,\ x \leq y \implies x = y$.
We write $\Antichain[P]$ for the set of all antichains of $P$.
\end{definition}

\begin{definition}[Min operator]
For any subset $S \subseteq P$, the \emph{Min operator} returns the
set of minimal elements:
\[
\Min(S) = \{ x \in S : \nexists y \in S.\, y < x \}.
\]
$\Min(S)$ is always an antichain. For a finite $S$, $\Min(S)$ is the
Pareto front of $S$.
\end{definition}

\subsection{Antichains as a Lattice}

The set of antichains of a poset carries a natural order: $A \leq A'$
in $\Antichain[P]$ iff every point in $A'$ is dominated by some point
in $A$. Intuitively, $A \leq A'$ means $A$ is "lower" or "easier",
having weaker Pareto-minimal demands.

\begin{proposition}
$(\Antichain[P], \leq)$ is a complete lattice when $P$ is finite or
chain-complete. Joins are computed by union followed by $\Min$. The
bottom of $\Antichain[P]$ is $\{\bot_P\}$ (or $\emptyset$ if $P$ has
no bottom); the top is $\emptyset$.
\end{proposition}

This lattice is the state space of the Kleene iteration (Section
\ref{sec:loop}).

\subsection{Design Problems}

\begin{definition}[Design Problem]
A \emph{design problem (DP)} is a tuple $\langle F, R, h \rangle$
where $F$ and $R$ are posets and $h : F \to \Antichain[R]$ is a
function from functionality requests to antichains of resources.
$h(f)$ is interpreted as the Pareto front of minimal resources
required to deliver functionality $f$.
\end{definition}

\begin{definition}[Monotone Design Problem]
A DP $\langle F, R, h \rangle$ is \emph{monotone} (an MCDP) if $h$
is monotone with respect to the antichain order: $f \leq f'$ in $F$
implies $h(f) \leq h(f')$ in $\Antichain[R]$.
\end{definition}

Intuitively, monotonicity says that demanding more functionality
cannot reduce the minimal resources needed. This corresponds to the
common-sense engineering principle that you don't get more for less.

\subsection{Composition}

The class of MCDPs is closed under three operators.

\begin{definition}[Series]
Given DPs $d_1 = \langle F, M, h_1 \rangle$ and $d_2 = \langle M, R,
h_2 \rangle$ whose interface poset $M$ matches, their \emph{series
composition} is
\[
(d_1 \circ d_2)(f) = \Min\!\left( \bigcup_{r_1 \in h_1(f)} h_2(r_1) \right).
\]
\end{definition}

\begin{definition}[Parallel]
Given DPs $d_1 = \langle F_1, R_1, h_1 \rangle$ and $d_2 = \langle
F_2, R_2, h_2 \rangle$, their \emph{parallel composition} is
$\langle F_1 \times F_2, R_1 \times R_2, h_1 \otimes h_2 \rangle$ where
\[
(h_1 \otimes h_2)(f_1, f_2) = \Min\!\left( h_1(f_1) \times h_2(f_2) \right).
\]
\end{definition}

\begin{definition}[Feedback]
Given a DP $d = \langle F \times X, R \times X, h \rangle$ with a
common factor $X$, the \emph{feedback composition} on $X$ is the DP
$d^{\dagger X} = \langle F, R, h^{\dagger X} \rangle$ where $h^{\dagger
X}(f)$ is the antichain obtained by closing the recursion on $X$.
\end{definition}

The recursion is solved by Kleene iteration in $\Antichain[R \times
X]$, seeded at $\{\bot\}$ and ascending until a fixed point is
reached.

\subsection{The Kleene Fixed-Point Algorithm}

\begin{theorem}[\citealp{censi2015mathematical}, Prop. 4]
\label{thm:kleene}
For an MCDP $d = \langle F \times X, R \times X, h \rangle$ closed
on $X$, define the operator $\Phi_f : \Antichain[R \times X] \to
\Antichain[R \times X]$ by
\[
\Phi_f(A) = \Min\!\left( \bigcup_{r \in A} h(f, r_X) \cap \uparrow\! r \right),
\]
where $r_X$ is the $X$-component of $r$ and $\uparrow\!r$ is the
upward closure. Then $\Phi_f$ is monotone, and the Kleene ascent
\[
A_0 = \{\bot\},\quad A_{k+1} = \Phi_f(A_k)
\]
converges to the least fixed point, which equals
$h^{\dagger X}(f)$ projected to $\Antichain[R]$.
\end{theorem}

In implementation, the iteration terminates when one of three
conditions holds: $A_{k+1} = A_k$ (fixed point), $A_{k+1} = \emptyset$
(infeasibility), or every point in $A_{k+1}$ has $r_X = \topp_X$
(loop axis saturates).

\section{Installation and Project Layout}
\label{sec:install}

\codename{} requires Python 3.9 or newer. The algorithmic core has no
required runtime dependencies; the optional layers listed below pull in
\texttt{numpy}, \texttt{scipy}, \texttt{matplotlib}, or \texttt{graphviz}
as needed.

\paragraph{Installation.} Clone the repository and install in editable
mode:
\begin{lstlisting}[language=bash]
git clone https://github.com/cbriat/codesign-mcdp.git
cd codesign-mcdp
pip install -e .
\end{lstlisting}

\paragraph{Optional extras.}
\begin{itemize}[leftmargin=2em]
  \item \texttt{pip install -e ".[viz]"} adds \texttt{matplotlib} for the
        visualisation helpers (Section~\ref{sec:visualisation}) and the
        plotting examples.
  \item \texttt{pip install -e ".[online]"} adds \texttt{numpy} and
        \texttt{scipy} for the uncertainty
        (Section~\ref{sec:uncertainty}) and online-learning
        (Section~\ref{sec:online-learning}) layers.
  \item \texttt{pip install -e ".[diagram]"} adds \texttt{graphviz} for
        the block-diagram renderer (Section~\ref{sec:block-diagrams}).
  \item \texttt{pip install -e ".[nb]"} adds the Jupyter tooling needed
        to (re)generate the notebooks.
\end{itemize}

\paragraph{Repository layout.}
\begin{verbatim}
codesign-mcdp/
  codesign/          the Python package (21 modules)
  examples/          25 runnable example scripts
  notebooks/         the same examples as executed .ipynb notebooks
  tests/             smoke tests
  docs/              this manual plus rendered trace images
  build_notebooks.py script to regenerate the notebooks
  pyproject.toml     packaging metadata, PEP 621
  .github/workflows/ CI configuration
\end{verbatim}

\section{Core Data Types}
\label{sec:posets}

This section documents the data types provided by
\texttt{codesign.posets} and \texttt{codesign.antichains}.

\subsection{Posets}

The \texttt{Poset} abstract base class has the following interface:

\begin{lstlisting}
class Poset(ABC):
    name: str

    def leq(self, x, y) -> bool: ...      # is x <= y?
    def eq(self, x, y) -> bool: ...
    def lt(self, x, y) -> bool: ...
    def bottom(self): ...                  # least element, if any
    def top(self): ...                     # greatest element, if any
    def is_bottom(self, x) -> bool: ...
    def is_top(self, x) -> bool: ...
    def join(self, x, y): ...              # least upper bound
    def format(self, x) -> str: ...
\end{lstlisting}

Four concrete poset classes are provided.

\paragraph{\texttt{Reals(unit="")}.} The non-negative reals augmented
with $+\infty$. Bottom is $0.0$; top is \texttt{math.inf}. The optional
\texttt{unit} string is used for pretty-printing only.

\paragraph{\texttt{Naturals(unit="")}.} The non-negative integers
augmented with $+\infty$. Bottom is $0$; top is \texttt{math.inf}.

\paragraph{\texttt{Ports(components: dict[str, Poset])}.}
The componentwise-ordered product of named factor posets. Elements
are Python dicts mapping component names to values. The bottom is the
dict of all bottoms; the top is the dict of all tops.
\texttt{Ports} can be nested freely: a component poset may
itself be a \texttt{Ports}. This is the type used for every design
problem's $F$ and $R$.

The name reflects the library's vocabulary: each named component is
a \emph{port} of the design problem, addressable by name in the
constraint DSL. For backward compatibility the older name
\texttt{NamedProduct} remains exported as an alias, so existing code
that imports \texttt{NamedProduct} continues to work; new code should
prefer \texttt{Ports}.

\paragraph{\texttt{Discrete(elements, leq=None)}.} A finite poset with
explicit element set and a caller-supplied \texttt{leq} predicate.
Useful for catalog-like enumerations where the order is not the
default discrete order.

\subsection{Antichains}

The \texttt{Antichain} class wraps a set of poset elements together
with the poset they live in. The constructor normalises the input by
removing dominated points and deduplicating.

\begin{lstlisting}
class Antichain:
    poset: Poset
    points: list

    @classmethod
    def of_bottom(cls, poset) -> "Antichain": ...
    @classmethod
    def empty(cls, poset) -> "Antichain": ...
    @classmethod
    def singleton(cls, poset, point) -> "Antichain": ...
    @classmethod
    def from_set(cls, poset, points) -> "Antichain": ...

    @classmethod
    def union_min(cls, poset, antichains) -> "Antichain": ...

    def filter_above(self, floor) -> "Antichain": ...
    def leq(self, other) -> bool: ...
    def eq(self, other) -> bool: ...
    def has_any_top(self) -> bool: ...
    def is_empty(self) -> bool: ...
\end{lstlisting}

The two key operations on antichains are \texttt{union\_min}, a class
method that takes the union of a collection of antichains and re-applies
$\Min$, and \texttt{filter\_above}, which keeps only points dominating a
given floor. Together they implement the body of the Kleene iteration in
Theorem~\ref{thm:kleene}; the full method surface is documented in the
reference card (Section~\ref{ref:antichains:antichain}).

\section{Primitive Design Problem Types}
\label{sec:primitives}

Six concrete subclasses of \texttt{DesignProblem} are provided, each
appropriate for a different style of relation between functionality
and resources.

\subsection{\texttt{AlgebraicDP}}

For relations where each resource is a closed-form monotone function
of the functionality.

\begin{lstlisting}
AlgebraicDP(
    F: Ports,
    R: Ports,
    equations: dict[str, Callable[[dict], float]],
    name: str = "algebraic",
)
\end{lstlisting}

The \texttt{equations} dictionary maps each resource name to a
callable taking a functionality dict and returning a scalar.
\texttt{h(f)} returns a singleton antichain.

\begin{example}
A battery sized by specific energy:
\begin{lstlisting}
battery = AlgebraicDP(
    F=Ports({"capacity": Reals(unit="J")}),
    R=Ports({"mass": Reals(unit="kg")}),
    equations={"mass": lambda f: f["capacity"] / 1.8e6},
)
\end{lstlisting}
\end{example}

\subsection{\texttt{FunctionDP}}

For relations where the user supplies \texttt{h} directly as a
function from functionality dicts to antichains. Use this when the
relation is multi-valued (a true Pareto front) or has branching logic.

\begin{lstlisting}
FunctionDP(
    F: Ports,
    R: Ports,
    h_fn: Callable[[dict], Antichain],
    name: str = "function",
)
\end{lstlisting}

\begin{example}[A two-option amplifier]
The same audio amplifier output spec admits two implementations: a
class-A topology (clean, hot) or a class-D topology (efficient, harsh).
\texttt{h(f)} emits both points, so the engineer sees the real
tradeoff:
\begin{lstlisting}
F = Ports({"watts":    Reals(unit="W")})
R = Ports({"power_in": Reals(unit="W"),
           "thd":      Reals()})

def h_amp(f):
    watts = f["watts"]
    return Antichain.from_set(R, [
        {"power_in": watts / 0.25, "thd": 0.001},   # class-A
        {"power_in": watts / 0.90, "thd": 0.020},   # class-D
    ])

amp = FunctionDP(F, R, h_amp)
\end{lstlisting}
\end{example}

\subsection{\texttt{CatalogDP}}

For selecting from a finite catalog of implementations (motors,
batteries, sensors).

\begin{lstlisting}
CatalogDP(
    F: Ports,
    R: Ports,
    catalog: list[CatalogEntry],
    name: str = "catalog",
)

CatalogEntry(
    provides: dict[str, Any],
    costs: dict[str, Any],
    name: str = "",
)
\end{lstlisting}

For a request $f$, $h(f)$ is the antichain of \texttt{entry.costs}
for every entry whose \texttt{provides} dominates $f$ in $F$. Entries
that are dominated by another feasible entry are pruned by $\Min$ in
$R$.

\begin{example}[A small motor catalog]
Three motors with incomparable mass--cost tradeoffs at the same
torque rating. \texttt{h} returns whichever entries cover the
requested torque, then $\Min$ drops dominated ones automatically:
\begin{lstlisting}
motor = CatalogDP(
    F=Ports({"torque": Reals(unit="N*m")}),
    R=Ports({"mass":   Reals(unit="kg"),
                    "cost":   Reals(unit="USD")}),
    catalog=[
        CatalogEntry(name="Tiny",   provides={"torque": 2.0},
                     costs={"mass": 0.2, "cost": 30.0}),
        CatalogEntry(name="Light",  provides={"torque": 8.0},
                     costs={"mass": 0.5, "cost": 200.0}),
        CatalogEntry(name="Heavy",  provides={"torque": 8.0},
                     costs={"mass": 0.8, "cost": 120.0}),
    ],
)
# motor.h({"torque": 7.0}) yields the antichain
# Antichain[(mass=0.5, cost=200), (mass=0.8, cost=120)]
\end{lstlisting}
The Tiny entry is filtered out (does not cover 7\,N$\cdot$m). Light
and Heavy are mutually incomparable, so both survive the $\Min$.
\end{example}

\subsection{\texttt{ConstraintDP}}

A feasibility predicate plus a scalar cost, lifted to $\Min$. Useful
when the resource-functionality relation is most natural to express
as \emph{is this choice good enough, and how much does it cost?}.

\begin{lstlisting}
ConstraintDP(
    F: Ports,
    R: Ports,
    sampler: Callable[[dict], Iterable[dict]],
    feasible: Callable[[dict, dict], bool],
    cost: Callable[[dict], dict],
    name: str = "constraint",
)
\end{lstlisting}

The \texttt{sampler} yields candidate implementations for the given
functionality; \texttt{feasible(impl, f)} filters them; \texttt{cost(impl)}
maps each feasible implementation into the resource poset, and the result
is reduced by $\Min$.

\subsection{\texttt{ODE\_DP}}

Derives a monotone resource relation from a differential equation
$\dot x = \text{rhs}(x, t, f)$. Two modes are supported.

\paragraph{Final value} integrates the ODE from $t = 0$ to
$t = t_{\text{end}}$ by explicit Euler and reports the extracted
resources of the final state.

\paragraph{Steady state} solves $\text{rhs}(x, \cdot, f) = 0$ for $x$
by Newton iteration and reports the extracted resources of the
steady-state value.

\begin{lstlisting}
ODE_DP(
    F: Ports,
    R: Ports,
    rhs: Callable[[state, float, dict], state],
    extract: Callable[[state], dict],
    mode: str = "final_value",   # or "steady_state"
    t_end: float = 10.0,
    n_steps: int = 200,
    x0_fn: Callable[[dict], state] = None,   # default: lambda f: 0.0
    name: str = "ode",
)
\end{lstlisting}

\begin{example}[A heater at thermal steady state]
A wall heater sized from Newton's law of cooling: to hold a temperature
rise $\Delta T$ above ambient, the input power must balance the heat loss
$k\,\Delta T$. Modelling the delivered power as a scalar state $x$ that
relaxes towards its steady value, $\dot x = k\,\Delta T - x$, the root is
$x = k\,\Delta T$. The heater's resource is that steady input power:
\begin{lstlisting}
F = Ports({"delta_T": Reals(unit="K")})
R = Ports({"power":   Reals(unit="W")})
k = 0.5   # heat-loss coefficient, W/K

heater = ODE_DP(
    F, R,
    rhs=lambda x, t, f: k * f["delta_T"] - x,
    extract=lambda x: {"power": float(x)},
    mode="steady_state",
    x0_fn=lambda f: 0.0,
)
# solve(heater, {"delta_T": 30.0}) -> Antichain[(power=15 W)]
\end{lstlisting}
The steady-state root $x = k\,\Delta T$ is recovered exactly by the
Newton iteration. Note that the state carried by \texttt{ODE\_DP} is a
plain scalar (or a list of scalars), not a dict: \texttt{rhs},
\texttt{extract}, and \texttt{x0\_fn} all operate on that state directly.
\end{example}

\subsection{\texttt{UncertainDP}}

Brackets an unknown $h$ with a known lower and upper bound. The
optimistic mode (\texttt{lower}) yields a lower bound on the minimal
resources; the pessimistic mode (\texttt{upper}) yields an upper
bound. Section~VII of \citet{censi2015mathematical} develops the theory,
and \citet{censi2017uncertainty} treats uncertainty in monotone co-design
problems in full.

\begin{lstlisting}
UncertainDP(
    F: Ports,
    R: Ports,
    lower: DesignProblem,
    upper: DesignProblem,
    mode: str = "upper",
    name: str = "uncertain",
)
\end{lstlisting}

\texttt{with\_mode("lower" | "upper")} returns the bracket selected
for solving.

\begin{example}[A battery with uncertain specific energy]
Specific energy of the cell chemistry is known only to lie in
$[1.6, 2.0]$\,MJ/kg. The optimistic and pessimistic brackets give
lower and upper bounds on the required mass:
\begin{lstlisting}
F = Ports({"capacity": Reals(unit="J")})
R = Ports({"mass":     Reals(unit="kg")})

lo = AlgebraicDP(F, R, {"mass": lambda f: f["capacity"] / 2.0e6})
hi = AlgebraicDP(F, R, {"mass": lambda f: f["capacity"] / 1.6e6})

battery = UncertainDP(F, R, lower=lo, upper=hi)

# Solve both brackets and report the interval.
m_opt  = solve(battery.with_mode("lower"), {"capacity": 3.6e6})
m_pess = solve(battery.with_mode("upper"), {"capacity": 3.6e6})
# m_opt  -> Antichain[(mass=1.80 kg)]
# m_pess -> Antichain[(mass=2.25 kg)]
\end{lstlisting}
\end{example}

\section{Composition Operators}
\label{sec:composition}

The three operators of Section~\ref{sec:theory} are exposed as classes
(\texttt{Series}, \texttt{Parallel}, \texttt{Loop}) with lowercase
function aliases (\texttt{series}, \texttt{par}, \texttt{loop}) matching
the paper's notation. Each takes existing design problems and returns a
new \texttt{DesignProblem}, so composites nest freely and are solved by
the same \texttt{solve} entry point. Because the class of MCDPs is closed
under all three, a composite is itself a monotone design problem. The
compositional structure is the point of the framework, and is what the
later literature builds on: modular co-design of vehicle control systems
\citep{zardini2022taskdriven}, structured co-design of embodied
intelligence \citep{zardini2021embodied}, and the categorical accounts of
\citet{furter2026symmetric} and \citet{riess2026quantale} all rest on it.
\Citet{carlone2019beyond} study what happens when monotonicity is dropped.

\subsection{Series}

\begin{lstlisting}
series(dp1: DesignProblem, dp2: DesignProblem) -> Series
Series(dp1, dp2)
\end{lstlisting}

Requires \texttt{dp1.R == dp2.F} structurally (same component names
and posets). Implements the series definition: for each point of
$h_1(f)$ run $h_2$ and take the union $\Min$.

\begin{example}[Battery feeding an actuator]
A battery converts energy demand into mass; the actuator converts
that mass demand into a power draw. Series composition feeds the
battery's output directly into the actuator's input:
\begin{lstlisting}
battery = AlgebraicDP(
    F=Ports({"capacity": Reals(unit="J")}),
    R=Ports({"mass":     Reals(unit="kg")}),
    equations={"mass": lambda f: f["capacity"] / 1.8e6},
)
actuator = AlgebraicDP(
    F=Ports({"mass":     Reals(unit="kg")}),
    R=Ports({"power":    Reals(unit="W")}),
    equations={"power": lambda f: 10.0 * f["mass"] ** 2},
)
chain = series(battery, actuator)
# solve(chain, {"capacity": 3.6e6}) maps energy -> power directly.
\end{lstlisting}
\end{example}

\subsection{Parallel}

\begin{lstlisting}
par(dp1: DesignProblem, dp2: DesignProblem) -> Parallel
Parallel(dp1, dp2)
\end{lstlisting}

Requires non-overlapping component names in $F_1, F_2$ and in
$R_1, R_2$. The output $F$ and $R$ are the concatenated
\texttt{Ports}s. The antichain is the Cartesian product
followed by $\Min$.

\begin{example}[Two independent subsystems]
A drone has a propulsion subsystem (energy in $\to$ mass out) and a
sensing subsystem (sample rate in $\to$ power out). They share no
ports and run in parallel; the composite DP takes both inputs and
emits both outputs:
\begin{lstlisting}
propulsion = AlgebraicDP(
    F=Ports({"energy":      Reals(unit="J")}),
    R=Ports({"prop_mass":   Reals(unit="kg")}),
    equations={"prop_mass": lambda f: f["energy"] / 1.8e6},
)
sensor = AlgebraicDP(
    F=Ports({"sample_rate": Reals(unit="Hz")}),
    R=Ports({"sensor_pwr":  Reals(unit="W")}),
    equations={"sensor_pwr": lambda f: 0.05 * f["sample_rate"]},
)
combined = par(propulsion, sensor)
# combined.F has both energy and sample_rate;
# combined.R has both prop_mass and sensor_pwr.
\end{lstlisting}
\end{example}

\subsection{Loop}
\label{sec:loop}

\begin{lstlisting}
loop(inner: DesignProblem, axis: str) -> Loop
Loop(inner, axis="...")
\end{lstlisting}

Closes the feedback loop on a single named axis that must appear in
both \texttt{inner.F} and \texttt{inner.R}. The resulting DP has
\texttt{F = inner.F} minus the axis and \texttt{R = inner.R} minus
the axis. Calling \texttt{h(f)} triggers the Kleene fixed-point
iteration described in Theorem~\ref{thm:kleene}.

\begin{example}[Battery mass feedback]
A drone's battery must carry both the payload \emph{and} its own
mass. Closing the loop on \texttt{battery\_mass} converts the
self-referential inequality into a fixed-point problem that the
Kleene solver handles automatically:
\begin{lstlisting}
inner = AlgebraicDP(
    F=Ports({"payload":       Reals(unit="kg"),
                    "battery_mass":  Reals(unit="kg")}),
    R=Ports({"battery_mass":  Reals(unit="kg")}),
    equations={
        "battery_mass":
            lambda f: 0.5 * (f["payload"] + f["battery_mass"]),
    },
)
drone = loop(inner, axis="battery_mass")
# drone.F = {"payload"}, drone.R = {}
# solve(drone, {"payload": 1.0}) iterates to the fixed point.
\end{lstlisting}
The closed-form solution is $m^\star = p / (1 - 0.5) = 2p$, the fixed
point to which the Kleene iteration converges (creeping to machine
precision over several dozen contraction steps).
\end{example}

\begin{remark}[Multiple feedback loops]
\texttt{Loop} closes one axis at a time. To close several loops,
chain them: \texttt{loop(loop(inner, axis="x"), axis="y")}. The
\texttt{System} builder (Section~\ref{sec:system-builder}) closes all
feedback loops simultaneously over a single bundled axis, which is
usually more convenient.
\end{remark}

\section{The Solver}
\label{sec:solver}

All design problems, primitive or composite, are solved through a single
entry point, \texttt{solve}. It evaluates the problem at a given
functionality, running the Kleene fixed-point iteration of
Theorem~\ref{thm:kleene} whenever a feedback loop is present, and returns
a \texttt{SolveResult} carrying the resulting antichain together with
convergence and feasibility metadata. This section documents that entry
point, the underlying iteration, scalar minimisation over the returned
front, and the solver's observability features.

\subsection{Top-level entry point}

\begin{lstlisting}
solve(
    dp: DesignProblem,
    functionality: dict | None = None,
    max_iter: int = 200,
    *,
    trace: bool = False,
    verbose: int = 0,
    on_iteration: Callable[[TraceEntry], None] | None = None,
    start_from: SolveResult | Antichain | None = None,
    uncertainty: list[str] | None = None,
    n_samples: int = 1000,
    rng_seed: int | None = None,
) -> SolveResult
\end{lstlisting}

The keyword-only arguments \texttt{trace}, \texttt{verbose},
\texttt{on\_iteration}, and \texttt{start\_from} are described in
Section~\ref{sec:solver-observability}; \texttt{uncertainty},
\texttt{n\_samples}, and \texttt{rng\_seed} in
Section~\ref{sec:uncertainty}. The result is a dataclass with the
following fields:
\begin{lstlisting}
@dataclass
class SolveResult:
    antichain: Antichain      # the Pareto front (in dp.R)
    iterations: int           # Kleene steps taken (0 if no loop)
    status: str               # "converged", "max_iter", or "diverged"
    feasible: bool            # True if antichain has finite, nonempty points
    trace: list[TraceEntry]   # the per-step trace if trace=True, else None

    # converged is a read-only alias for status == "converged"
\end{lstlisting}

\begin{example}[Inspecting a solve result]
The full lifecycle of a solve: build the DP, call \texttt{solve},
introspect the result, and pick a single design with
\texttt{minimize\_cost}:
\begin{lstlisting}
result = solve(drone, {"endurance":      300.0,
                       "extra_payload":  0.5,
                       "extra_power":    5.0},
               max_iter=200, trace=True)

print(result.status)         # "converged"
print(result.iterations)     # e.g. 17
print(result.feasible)       # True
print(result.antichain)      # Antichain[(total_mass=0.55 kg)]

# Scalarise: cheapest design under a custom cost function
best = minimize_cost(result, cost_fn=lambda r: r["total_mass"])
print(best)                  # {"total_mass": 0.5492}
\end{lstlisting}
\end{example}

\subsection{The Kleene iteration}

\begin{lstlisting}
kleene_loop(
    loop_dp: Loop,
    f_outer: dict | None,
    max_iter: int = 200,
    *,
    trace: bool = False,
    verbose: int = 0,
    on_iteration: Callable[[TraceEntry], None] | None = None,
    start_from: Antichain | None = None,
    info_out: dict | None = None,
) -> Antichain
\end{lstlisting}

The body of the iteration is direct from Theorem~\ref{thm:kleene}. The
implementation adds:
\begin{itemize}
  \item a \emph{divergence cap} of $10^{30}$ that converts an iterate
        exceeding the cap to $\topp$, so floating-point overflow
        becomes infeasibility rather than \texttt{inf}/\texttt{nan};
  \item \texttt{try}/\texttt{except} around the inner $h$ to catch
        \texttt{OverflowError}, \texttt{ValueError}, and
        \texttt{ZeroDivisionError}, again mapping to $\topp$;
  \item three termination conditions: fixed point reached, antichain
        empty (infeasible), or every point's loop axis equal to
        $\topp$ (provably infeasible).
\end{itemize}

\subsection{Scalar minimisation over an antichain}

\begin{lstlisting}
minimize_cost(
    result: SolveResult,
    cost_fn: Callable[[dict], float],
) -> dict | None
\end{lstlisting}

Returns the single resource bundle in \texttt{result.antichain} that
minimises \texttt{cost\_fn}, or \texttt{None} if the result is
infeasible. This is the scalarisation step that turns a Pareto front
into a single engineering choice.

\subsection{Observability: trace, verbose, callback, status}
\label{sec:solver-observability}

The solver can be made to report its progress in three independent
ways, and its termination reason is exposed through a structured
\texttt{status} field that is orthogonal to \texttt{feasible}.

\paragraph{The \texttt{status} field.} \texttt{SolveResult.status} is
one of three strings:
\begin{description}
\item[\texttt{"converged"}] The iteration reached a fixed point (or
provably hit $\top$) within \texttt{max\_iter} steps. The answer in
\texttt{result.antichain} is the algorithm's verdict; \texttt{feasible}
tells you whether that verdict is finite (a real design) or $\top$
(no design exists).
\item[\texttt{"max\_iter"}] The iteration was cut off at
\texttt{max\_iter}. The current antichain may or may not be near a
fixed point. Usually a sign to increase \texttt{max\_iter}, or to look
for an unintentional non-monotonicity in the model.
\item[\texttt{"diverged"}] At least one numeric component crossed the
divergence cap of $10^{30}$ before the iteration could settle, and the
solver stopped. Distinguishes pure numerical blow-up from clean
$\top$-infeasibility.
\end{description}
The previous Boolean \texttt{converged} field is preserved as a
backward-compatibility alias for \texttt{status == "converged"}.

\paragraph{Structured trace.} Passing \texttt{trace=True} to
\texttt{solve} populates \texttt{result.trace} with a list of
\texttt{TraceEntry} dataclasses, one per Kleene step (plus the seed at
iteration~0):
\begin{lstlisting}
@dataclass
class TraceEntry:
    iteration: int           # 0 = seed, then 1, 2, ...
    antichain: Antichain     # snapshot at this step
    n_points:  int           # convenience: len(antichain)
    delta:     float | None  # max absolute change (numeric)
                             # or 0/1 (discrete); None at iter 0
    elapsed_ms: float        # wall time for this step alone
\end{lstlisting}
With \texttt{trace=False} (the default) the field is \texttt{None}, so
the cost of unused tracing is zero.

\paragraph{Verbose printing.} The integer \texttt{verbose} controls
live output:
\begin{description}
\item[\texttt{verbose=0}] (default) silent.
\item[\texttt{verbose=1}] one summary line at the end, e.g.\
\texttt{[solve] converged: 17 iters, |A|=1, total=2.0ms, feasible=True}.
\item[\texttt{verbose=2}] one line per iteration, showing the delta,
the antichain size, and the step's wall time. Useful for watching
oscillating fixed points.
\end{description}

\paragraph{Iteration callback.} Passing
\texttt{on\_iteration=callable} registers a callback that receives
each \texttt{TraceEntry} as soon as it is produced. Useful for live
plots, custom logging, or interrupting an iteration that looks
pathological:
\begin{lstlisting}
def log_every_5(entry):
    if entry.iteration % 5 == 0:
        print(f"iter {entry.iteration}: delta={entry.delta:.3e}")

solve(dp, f, on_iteration=log_every_5)
\end{lstlisting}

\begin{example}[Warm-starting a parameter sweep]
When sweeping a parameter (load, mission length, demand level), the
fixed point at one parameter is a great seed for the next. Pass the
previous \texttt{SolveResult} (or its \texttt{antichain}) via
\texttt{start\_from=} and the Kleene iteration picks up from there
instead of restarting at $\bot$:
\begin{lstlisting}
prev = None
results = []
for L in np.linspace(5.0, 30.0, 50):
    r = solve(dp, {"daily_load_kwh": float(L), ...},
              max_iter=400, start_from=prev)
    results.append(r)
    prev = r            # reuse the converged inner antichain
\end{lstlisting}
For the microgrid example (notebook~13) this reduces the cumulative
iteration count by roughly 10\,\% across the sweep, more for finer
parameter grids.
\end{example}

\begin{example}[Plotting convergence after a trace]
With \texttt{trace=True}, the result carries the entire iteration
history; a single \texttt{viz.plot\_convergence} call renders the
delta-vs-iteration semilog:
\begin{lstlisting}
from codesign import viz

r = solve(drone, f, trace=True, max_iter=200)
ax = viz.plot_convergence(r)
ax.set_title(f"converged in {r.iterations} steps")
\end{lstlisting}
\end{example}

\section{Uncertainty}
\label{sec:uncertainty}

Uncertainty in monotone co-design was introduced by
\citet{censi2017uncertainty}, and the composable and distributional
treatments of
\citet{huang2025composable}, \citet{furter2026symmetric}, and
\citet{huang2026distributional} extend it. This section documents what the
library implements, which is the set-based and stochastic part.

Two kinds of parameter uncertainty are supported, declared on
\texttt{Module} instances and consumed by the same \texttt{solve}
entry point. Both are independent of the model's structure: a module
can have set-based uncertainty, stochastic uncertainty, both, or
neither.

\subsection{Declaration}

A \texttt{Module} carries optional attributes:
\begin{lstlisting}
class Battery(Module):
    F = {"capacity": Reals(unit="J")}
    R = {"mass":     Reals(unit="kg")}

    def __init__(self, specific_energy=2.0e6, efficiency=0.9):
        self.specific_energy = specific_energy
        self.efficiency = efficiency
        super().__init__()

    def h(self, f):
        # nominal product 2.0e6 * 0.9 = 1.8 MJ/kg delivered,
        # matching the canonical drone (examples 1/6/7)
        return {"mass": f["capacity"]
                       / (self.specific_energy * self.efficiency)}

b = Battery()
b.uncertain_set  = Box(...)          # deterministic, set-based
b.uncertain_dist = Stochastic(...)   # statistical
\end{lstlisting}
The uncertainty solver discovers every such module by walking the
\verb|_codesign_modules| attribute that \texttt{System.build} attaches
to the returned DP. Before each solve, the module's named parameters
are reassigned; afterwards, their nominal values are restored. The
user's original \texttt{Battery()} instance is unchanged.

\subsection{Set-based uncertainty (deterministic)}

\paragraph{\texttt{Box}.} Axis-aligned interval product. Each
parameter is declared as \texttt{(lo, hi)} or \texttt{(lo, hi, direction)},
where \texttt{direction} is one of the four tokens
\texttt{"more\_is\_better"}, \texttt{"more\_is\_worse"},
\texttt{"less\_is\_better"}, \texttt{"less\_is\_worse"}. With all
directions declared, the worst case is a single corner read off in
constant time:
\begin{lstlisting}
b.uncertain_set = Box(
    specific_energy=(1.7e6, 2.3e6, "more_is_better"),
    efficiency=(0.83, 0.97, "more_is_better"),
)
# worst case = (specific_energy=1.7e6, efficiency=0.83)
\end{lstlisting}
When some directions are undeclared, all $2^n$ endpoint combinations
are evaluated and the worst is kept; cheap when $n$ is modest.

\paragraph{\texttt{Ellipsoid}.} An $n$-D ellipsoid
$(p - c)^\top \Sigma^{-1} (p - c) \leq 1$ in parameter space. The
covariance $\Sigma$ encodes both scale and correlation between
parameters, so an ellipsoid is a more honest model than a box when
parameters are believed to vary together:
\begin{lstlisting}
b.uncertain_set = Ellipsoid(
    center={"specific_energy": 2.0e6, "efficiency": 0.9},
    cov=[
        [1.0e10, -2.0e3],
        [-2.0e3,  2.5e-3],
    ],
    params=["specific_energy", "efficiency"],
    directions={
        "specific_energy": "more_is_better",
        "efficiency":      "more_is_better",
    },
)
\end{lstlisting}
With all directions declared, the worst case is the boundary point
$c + L\,u^\star$, where $L$ is the Cholesky factor of $\Sigma$ and
$u^\star$ is the unit vector in the direction of badness. Otherwise,
the boundary is sampled (controlled by \texttt{boundary\_samples}) and
the worst sample is returned.

\paragraph{2D conveniences.} \texttt{Disk(center, radius, ...)} and
\texttt{Circle(center, radius, ...)} reduce to the isotropic
ellipsoid $\Sigma = r^2 I$. For monotone systems with declared
directions, the two are equivalent: the worst case lies on the
boundary either way.

\subsection{Stochastic uncertainty}

\paragraph{\texttt{Stochastic}.} A joint distribution built from
named scipy-stats frozen marginals plus a copula:
\begin{lstlisting}
from scipy import stats

b.uncertain_dist = Stochastic(
    marginals={
        "specific_energy": stats.uniform(loc=1.7e6, scale=0.6e6),
        "efficiency":      stats.uniform(loc=0.83, scale=0.14),
    },
    copula=GaussianCopula(correlation=[[1.0, 0.4],
                                        [0.4, 1.0]]),
)
\end{lstlisting}
Sampling is by the standard copula--marginal trick: the copula
generates correlated uniforms in $[0,1]^d$; the marginal \texttt{ppf}s
map them to the named parameter axes.

\paragraph{Copulas.} \texttt{Independence()} is the default
(coordinatewise uniform). \texttt{GaussianCopula(correlation)} uses
the Gaussian copula with the given correlation matrix; samples are
drawn by Cholesky factorisation followed by $\Phi$ on each
coordinate. Subclassing \texttt{Copula} requires only one method,
\texttt{sample\_uniform(n, d, rng)}, so other copulas (e.g.\ Clayton,
Gumbel) drop in without further plumbing.

\subsection{Querying}

The same \texttt{solve} entry point accepts an \texttt{uncertainty}
list and returns a richer object:
\begin{lstlisting}
result = solve(
    dp, f,
    uncertainty=["worst_case", "mean", "p95", "cvar95", "samples"],
    n_samples=1000,
    rng_seed=42,
)

result.worst_case        # SolveResult at the worst point of the set
result.mean              # dict[r_port -> mean across MC samples]
result.p95               # dict[r_port -> 95th percentile]
result.cvar95            # dict[r_port -> CVaR at 95% level]
result.samples           # list[Antichain], one per MC sample
result.feasibility_rate  # fraction of feasible MC samples
\end{lstlisting}

The allowed labels are:
\begin{description}
\item[\texttt{"worst\_case"}] The deterministic worst over every
module's \texttt{uncertain\_set}. One canonical solve at the
worst-case parameter point.
\item[\texttt{"mean"}, \texttt{"p95"}, \texttt{"cvar95"}] Marginal
statistics over Monte Carlo samples, drawn from every module's
\texttt{uncertain\_dist}. The CVaR (conditional value at risk) at the
95\% level is the mean of the worst 5\% of samples.
\item[\texttt{"samples"}] The raw antichain per MC sample, for the
user to aggregate however they like.
\end{description}

For a single solve call you typically observe the ordering
\[
\text{nominal} < \text{mean} < \text{p95} < \text{CVaR95} < \text{worst\_case},
\]
which separates the four standard ways of reporting an answer "under
uncertainty." Notebooks 11 and 12 visualise this for the
example-7 drone.

\begin{example}[Box worst-case of a drone battery]
Specific energy and efficiency of the battery cell live in declared
ranges, each with a direction of badness. The worst-case mass is
the corner where both parameters take their worst endpoint:
\begin{lstlisting}
bat = Battery()  # the Module subclass from example~7
bat.uncertain_set = Box(
    specific_energy=(1.7e6, 2.3e6, "more_is_better"),
    efficiency=(0.83, 0.97, "more_is_better"),
)
# ... wire `bat` into the drone system ...

result = solve(drone, f, uncertainty=["worst_case"])
worst_mass = list(result.worst_case.antichain.points)[0]["total_mass"]
# worst_mass = 0.5668 kg, versus 0.5492 kg nominal
\end{lstlisting}
\end{example}

\begin{example}[Monte Carlo with a Gaussian copula]
The same battery, now with stochastic specific energy and
efficiency, correlated at $\rho = 0.4$. One \texttt{solve} call
returns all four summary statistics:
\begin{lstlisting}
from scipy import stats

bat.uncertain_dist = Stochastic(
    marginals={
        "specific_energy": stats.uniform(loc=1.7e6, scale=0.6e6),
        "efficiency":      stats.uniform(loc=0.83, scale=0.14),
    },
    copula=GaussianCopula(correlation=[[1.0, 0.4],
                                        [0.4, 1.0]]),
)

res = solve(drone, f,
            uncertainty=["mean", "p95", "cvar95", "samples"],
            n_samples=1000, rng_seed=42)
# res.mean["total_mass"]    -> 0.5506 kg
# res.p95["total_mass"]     -> 0.5632 kg
# res.cvar95["total_mass"]  -> 0.5647 kg
# len(res.samples)          -> 1000
\end{lstlisting}
\end{example}

\section{Online Learning}
\label{sec:online-learning}

When a co-design problem has many discrete candidates (catalog
entries, robot types, component families) and each candidate's
inner solve is non-trivial, the naive ``evaluate every entry''
strategy is wasteful. The \texttt{codesign.online} module implements
the compositional online learner of
\citet{alharbi2026online}: maintain history-dependent bounds on
each candidate's inner-solve output, evaluate the most promising one
under the \texttt{lcb} (lower-confidence-bound) rule by default, then
prune any candidate whose lower bound is already dominated by the
incumbent.

\subsection{Optimistic evaluators}

The bound machinery is encapsulated in a small class hierarchy.
\texttt{OptimisticEvaluator} is the abstract base; subclasses
implement \texttt{bound(candidate)} returning a pair of dicts
\texttt{(lower, upper)} mapping each numeric R component to its
current lower and upper bound at the queried feature point. Tighter
bounds prune more candidates; looser bounds are always safe.

Three concrete evaluators are provided:

\paragraph{Monotonicity.}
\texttt{MonotonicityEvaluator(features, r\_components)} assumes the
inner-solve output is component-wise monotone in the named features.
Given an observation at feature point $f_0$ with antichain-min $R_0$,
any candidate whose features dominate $f_0$ has resources $\geq R_0$
(lower bound); any candidate dominated by $f_0$ has resources
$\leq R_0$ (upper bound). Aggressive when applicable. Only correct if
the monotonicity assumption genuinely holds: choosing a derived
feature (e.g.\ \texttt{cost\_per\_capacity = unit\_cost /
(speed * payload)}) is often what makes this safe in practice.

\paragraph{Lipschitz.}
\texttt{LipschitzEvaluator(features, r\_components, L)} assumes
$|h(c_1) - h(c_2)| \leq L \cdot \|features(c_1) - features(c_2)\|$
(Euclidean distance, per R component). Each observation tightens the
bound by a cone of slope $L$. The constant can be a scalar (same for
every R component) or a dict per component. Safe default: with a
sensible $L$ the bound never prunes a Pareto-optimal candidate.

\paragraph{Linear-parametric.}
\texttt{LinearParametricEvaluator(features, r\_components,
min\_obs=3, prior\_box=None, noise\_bound=0.0, solver="highs")}
implements the \emph{certified} optimistic bound of
\citet{alharbi2026online} (Section~V-C3, eqs.\ 26--28 and Lemma~V.5).
It assumes each resource
coordinate is an \emph{exact} affine function of the features,
$\text{req}_k(c) = \phi(c)^\top \theta_k^*$ with
$\phi(c) = [1, features(c)]$ (the paper's noiseless linear model).
Rather than fit a single least-squares line, it maintains the whole
\emph{confidence polytope} $\Theta(H)$ --- the subset of the prior
box $\Theta_0$ consistent with every observation taken as a linear
equality --- and returns, per resource coordinate, the coordinatewise
minimum of the predicted resource over that polytope,
$[\text{req}_\text{opt}(i, H)]_k =
\min_{\theta \in \Theta(H)} \phi(i)^\top \theta$. Each such minimum is
one linear program, solved with \texttt{scipy.optimize.linprog} (the
\texttt{solver} kwarg selects the method). Because the true parameter
$\theta_k^*$ is feasible for the LP, its optimum is always
$\leq \phi(i)^\top \theta_k^* = \text{req}_k(i)$: the returned value is
a \emph{guaranteed} lower bound (the optimism guarantee of
Lemma~V.5), so --- unlike the former OLS $\pm$ confidence-band
heuristic --- it can \emph{never} wrongly eliminate a Pareto-optimal
candidate. No upper bound is certified, so the returned upper bound
stays at $+\infty$. \texttt{prior\_box} sets $\Theta_0$ (a per-parameter
box) to keep the LP bounded while the fit is under-determined:
\texttt{None} leaves every parameter unbounded, which is always safe
but yields the trivial bound until enough observations pin $\theta$
down; \texttt{noise\_bound}${} > 0$ relaxes each observation equality
into a two-sided band of that half-width, a documented extension for
bounded observation noise. The old \texttt{confidence} kwarg is
deprecated and ignored (passing it emits a
\texttt{DeprecationWarning}). Net effect: this evaluator is both safe
\emph{and} effective when the affine assumption genuinely holds, and
degrades \emph{safely} --- to the no-information bound, never to a
wrong elimination --- when it does not.

\paragraph{Gaussian process.}
\texttt{GaussianProcessEvaluator(features, r\_components,
length\_scale=0.3, sigma\_f=1.0, noise=1e-3, confidence=2.0,
min\_obs=3)} fits a zero-mean GP with an RBF kernel and bounds new
queries by a confidence band on the predictive standard deviation.
Implemented in pure numpy. The right tool when the response surface
has feature interactions or local nonlinearity: the GP still produces
informative (if uncertified) bounds there, whereas the certified
linear-parametric evaluator correctly degrades to the no-information
bound once the affine assumption fails. When the map \emph{is} affine,
prefer the linear-parametric evaluator --- its bound is then both
certified-safe and tight.

\paragraph{Subclassing.} Custom evaluators only need to override
\texttt{bound(candidate)}; the base class handles observation
recording and the default fallback bound of $(0, +\infty)$.

\subsection{The online solver}

\begin{lstlisting}
result = solve_online(
    candidate_fn,         # candidate_fn(candidate) -> DP
    functionality,        # outer F vector
    *,
    candidates,           # list of feature dicts
    evaluator,            # OptimisticEvaluator instance
    budget=None,          # max inner solves; None = unbounded
    max_iter=200,         # forwarded to each inner solve
    verbose=0,
    warm_start=None,      # seed indices or n for farthest-point
    picker="lcb",         # "lcb", "ucb", "random", or callable
)
\end{lstlisting}

Algorithm (a simplified Algorithm~2 from Alharbi et al.):
\begin{enumerate}
\item Bound every remaining candidate via the evaluator.
\item Eliminate every candidate whose lower bound is already
  dominated by the current incumbent antichain.
\item Pick the most promising survivor with the configured picker
  (\texttt{lcb} by default) on its lower-bound sum.
\item Run the inner solve at that candidate; observe the result in
  the evaluator; merge into the incumbent via \texttt{union\_min}.
\item Repeat until the candidates are exhausted or the budget is hit.
\end{enumerate}

The returned \texttt{OnlineResult} exposes:
\begin{description}
\item[\texttt{antichain}] Min over the evaluated, surviving candidates.
\item[\texttt{n\_evaluated}, \texttt{n\_eliminated},
       \texttt{n\_candidates}] Bookkeeping for the elimination cascade.
\item[\texttt{evaluated\_ids}, \texttt{eliminated\_ids},
       \texttt{incumbent\_ids}] Lists of indices into the original
  candidate list.
\item[\texttt{history}] Per-iteration log:
  \texttt{\{pick, antichain, remaining, evaluated, eliminated\}}.
\end{description}

\begin{example}[Fleet sizing across 200 robot types]
A logistics service has a 200-entry catalog of candidate robot
types, each described by four features. The Lipschitz evaluator
caps how fast the inner-solve output can change with the features,
so a single observation tightens the bound on every other candidate.
With the right L the answer is exact:
\begin{lstlisting}
def candidate_fn(robot):
    capacity = robot["speed"] * robot["payload"]
    return AlgebraicDP(F, R, {
        "total_cost":   lambda f, cap=capacity, uc=robot["unit_cost"]:
                        (f["target_throughput"] / cap) * uc,
        "total_energy": lambda f, ek=robot["energy_per_km"]:
                        f["target_range"] * ek * 24.0,
    })

ev = LipschitzEvaluator(
    features=["speed", "payload", "unit_cost", "energy_per_km"],
    r_components=["total_cost", "total_energy"],
    L={"total_cost": 300.0, "total_energy": 30.0},
)

result = solve_online(
    candidate_fn, mission,
    candidates=robot_catalog,        # list of 200 feature dicts
    evaluator=ev,
)
# result.n_evaluated   = 192 / 200
# result.n_eliminated  = 8
# result.antichain     = 5 Pareto-optimal robot types
\end{lstlisting}
Swapping in \texttt{MonotonicityEvaluator(features=["cost\_per\_capacity",\
\ "energy\_per\_km"], ...)} cuts the evaluation count to about 13 by
exploiting the structural monotonicity of the derived feature.
\end{example}

\subsection{Warm-start, picker strategies, and the GP evaluator}
\label{sec:online-tier1}

The basic online solver has three small extensions that significantly
expand the tuning space without complicating the core algorithm.

\paragraph{Warm-start.} A new \texttt{warm\_start} argument to
\texttt{solve\_online} pre-populates the evaluator with seed
observations before the picker takes over. It accepts either a list
of candidate indices (manually picked corner runs) or an integer
$n$ that triggers a greedy farthest-point heuristic over the feature
space. The motivation is concrete: the
\texttt{MonotonicityEvaluator}'s lower bounds only tighten for
candidates whose features dominate every observation in the partial
order. Without observations at the low-feature corner where the
Pareto front typically lives, the picker explores indiscriminately
and the evaluator remains uninformative throughout the budget. In
example 16, evaluating the example with four corner warm-start runs
lifts Monotonicity's Pareto recovery from zero to one of four
classes; a hybrid with one of the other evaluators is the natural
next step for cases where this is not enough.

\paragraph{Pluggable picker strategies.} A new \texttt{picker}
argument selects the candidate-scoring rule. Built-in options:
\texttt{"lcb"} (the previous default, minimises the sum of
lower-bound components, pure exploitation of the optimistic
estimate); \texttt{"ucb"} (lower bound minus a $\kappa \cdot$
(upper $-$ lower) exploration bonus, tunable via the tuple form
\texttt{("ucb", \{"kappa": 1.0\})}); and \texttt{"random"}
(uniform baseline for comparing the value of structural priors).
Custom callables are accepted and receive
\texttt{(lower, upper, r\_components, **kwargs)}. The strategies
that matter in practice are LCB for confident priors and UCB with
$\kappa \approx 0.3$ to $0.5$ when the prior is uncertain enough
that pure exploitation gets stuck in a local region.

\paragraph{Gaussian-process evaluator.} A new
\texttt{GaussianProcessEvaluator} class uses a zero-mean GP with
an RBF kernel, implemented in pure numpy. Hyperparameters are
\texttt{length\_scale} (default $0.3$, suitable for features in
$[0, 1]$), \texttt{sigma\_f} (default $1.0$, rescaled per-output by
the empirical observation standard deviation), \texttt{noise}
(default $10^{-3}$ jitter), \texttt{confidence} (default $2.0\sigma$
for roughly $95\%$ coverage), and \texttt{min\_obs} (default $3$).
The GP is more expressive than the linear-parametric evaluator and
captures local nonlinearity that a global linear fit misses. Since the
linear-parametric evaluator became the certified confidence-polytope
bound the two occupy clearly separated niches. The certified
\texttt{LinearParametricEvaluator} is the tool of choice when the
resource map is genuinely affine in the features: on the linear fleet
catalogue of example 14 it recovers all five Pareto types (where the
former OLS-band version wrongly dropped one). But the example 16
bioprocess effect model is markedly nonlinear (temperature U-shapes,
a power-law failure rate), and there the certified bound correctly
refuses to over-claim --- it degrades to the no-information bound and
recovers zero of the four Pareto classes rather than eliminate a
candidate it cannot rule out. The GP, willing to extrapolate a smooth
surrogate, stays informative on that same grid (recovering one of the
four classes). Rule of thumb: certified linear-parametric when you
trust the affine assumption (safe \emph{and} effective), GP when you
do not.

\section{Visualisation}
\label{sec:visualisation}

The \texttt{codesign.viz} module exposes four helpers built on
\texttt{matplotlib} and a small GraphViz emitter:

\begin{description}
\item[\texttt{plot\_antichain(result, axes)}] Scatter the antichain
  on the chosen 2D or 3D axes. Accepts a \texttt{SolveResult}, an
  \texttt{UncertaintyResult} (uses its worst case), or a bare
  \texttt{Antichain}. The 2D version optionally shades each point's
  dominated region for visual clarity.
\item[\texttt{plot\_convergence(result)}] Semilog plot of the
  Kleene delta against iteration index. Accepts a \texttt{SolveResult}
  with a trace or a trace list directly. Useful for diagnosing
  oscillating fixed points and tuning \texttt{max\_iter}.
\item[\texttt{plot\_uncertainty(unc\_result, port, nominal=None)}]
  Histogram of the Monte Carlo samples for the named R port, with
  the nominal, mean, p95, and CVaR95 marked as vertical lines.
\item[\texttt{to\_dot(dp, name="codesign")}] GraphViz dot string for the
  system structure: modules as boxes, outer ports as labelled
  rectangles, constraint connections as edges. Render with
  \texttt{dot -Tpng} or paste into an online viewer.
\end{description}

All four accept an optional \texttt{ax=...} argument so they compose
into larger figures. The module is imported as
\texttt{from codesign import viz}.

\subsection{Block diagrams of Systems}
\label{sec:block-diagrams}

The \texttt{to\_dot} function above emits a plain string with
module-level wiring. For richer visualisation, the library ships a
dedicated diagram renderer in \texttt{codesign.diagram} that produces
Simulink-style block diagrams:

\begin{description}
\item[One box per subsystem] with the F (functionality) ports listed
  on the left and the R (resource) ports on the right. Each port is
  individually addressable, so edges attach to specific ports rather
  than to the module as a whole.
\item[Outer F and outer R as separate nodes] on the diagram's left
  and right margins. The system inputs and outputs are visible, not
  buried inside the constraint list.
\item[Port-level edges] whenever the constraint was written in the
  operator-overloaded form
  (\texttt{module1.r\_port >= module2.r\_port * ...}). Constraints
  defined by callables that cannot be introspected are rendered with
  a dashed edge from a small ``$\lambda$'' marker, so they remain
  visible.
\item[Cycle detection] via Tarjan's strongly-connected-components
  algorithm on the module-level constraint graph. Edges within any
  cycle of size $\geq 2$ are coloured amber, making the
  Kleene-iteration loop visible at a glance.
\end{description}

The renderer returns a \texttt{graphviz.Digraph} object, which can be
displayed inline in a Jupyter notebook or written to SVG / PDF / PNG
via \texttt{dot.render(filename, format="svg")}:

\begin{lstlisting}
from codesign import draw_system

dot = system.draw_diagram()        # convenience method on System
# or equivalently:
dot = draw_system(system, rankdir="LR", highlight_cycles=True)

dot.render("bioprocess", format="svg", cleanup=True)
\end{lstlisting}

Optional dependency: install via
\texttt{pip install codesign-mcdp[diagram]} (which pulls the
\texttt{graphviz} Python package), and ensure the \texttt{dot}
binary is on \texttt{PATH} (\texttt{apt-get install graphviz} on
Debian / Ubuntu, \texttt{brew install graphviz} on macOS).

\begin{example}[Bioprocess block diagram]
The worked example~15 (\texttt{make\_bioprocess}) builds a four-
module System with a single outer F (target titer) and three outer
R aggregations defined by callables. The diagram renderer surfaces
this structure compactly:
\begin{lstlisting}
from codesign import draw_system
dp = make_bioprocess(cell_line=CHO_K1, glucose_setpoint_mm=8.0,
                     annual_demand_kg=100.0)
draw_system(dp, name="bioprocess").render("bioprocess",
                                          format="svg",
                                          cleanup=True)
\end{lstlisting}
The resulting diagram shows the four subsystems
(\texttt{cell}, \texttt{feed}, \texttt{bior}, \texttt{media}),
the outer functionality \texttt{target\_titer} feeding
\texttt{cell.target\_titer}, the metabolic-load chain
\texttt{cell.peak\_vcd}, \texttt{feed.metabolic\_factor}
$\to$ \texttt{bior.peak\_vcd} and \texttt{media.peak\_vcd}, and a
dashed $\lambda$ marker feeding the three outer R aggregations.
\end{example}

\begin{example}[Cyclic drone modular system]
Example~7's modular drone has a battery--actuator feedback
loop: the battery's mass adds to the actuator's lift demand, and the
actuator's electrical power sets the battery's required capacity.
\texttt{draw\_system} detects the cycle automatically and colours
the two edges between \texttt{battery} and \texttt{actuator} in
amber, while the supporting edges (from outer F endurance and
extra\_power into the cycle, and from each module's mass to the
outer R total\_mass) remain in muted grey.
\end{example}

\begin{example}[Pareto front of a vehicle co-design]
The vehicle from worked example~8 has two Pareto-incomparable
designs at the Small-parcel mission. A single call paints them on a
mass--cost plane and shades the dominated quadrants so the front is
visible at a glance:
\begin{lstlisting}
from codesign import viz
import matplotlib.pyplot as plt

result = solve(vehicle, {"payload": 2.0, "mission_energy": 2.0e5})
viz.plot_antichain(result, axes=["total_mass", "total_cost"])
plt.show()
\end{lstlisting}
For three resources, pass three axis names to get a 3D scatter.
\end{example}

\begin{example}[Uncertainty histogram with summaries]
Visualising the Monte Carlo output of a stochastic solve: the
histogram of total mass across MC samples, with the nominal, mean,
p95, and CVaR95 each drawn as a vertical line:
\begin{lstlisting}
res = solve(drone, f,
            uncertainty=["mean", "p95", "cvar95", "samples"],
            n_samples=1000, rng_seed=42)
ax = viz.plot_uncertainty(res, port="total_mass",
                          nominal=0.5492, bins=30)
\end{lstlisting}
The function reads the requested summaries off the result and
labels them automatically.
\end{example}

\begin{example}[System graph in GraphViz dot]
For documentation or for understanding someone else's model, the
constraint graph of a built \texttt{System} can be dumped as dot:
\begin{lstlisting}
dot = viz.to_dot(microgrid, name="microgrid")
# Render externally:
# echo "<dot output>" | dot -Tpng > microgrid.png
\end{lstlisting}
Modules become boxes, outer F ports become source rectangles,
outer R ports become sink rectangles, and constraint connections
become directed edges labelled with the linking expression.
\end{example}

\section{The MCDP Builder}
\label{sec:mcdp-builder}

The \texttt{MCDP} class in \texttt{codesign.mcdpl} provides an
MCDPL-style declarative syntax for a single self-contained design
problem.

\begin{lstlisting}
with MCDP("name") as m:
    m.provides("functionality_name", unit="...")    # outer F port
    m.requires("resource_name",     unit="...")    # outer R port

    m.constraint("resource_name", lambda f: <expr>)   # closed-form
    # or
    m.rule(lambda f: Antichain.from_set(...))         # multi-valued

    m.loop_on("axis_name")    # optional, may be repeated

dp = m.build()
\end{lstlisting}

The builder emits an \texttt{AlgebraicDP} or \texttt{FunctionDP}
internally, then wraps it in \texttt{Loop}s as requested by
\texttt{loop\_on}. Multiple \texttt{constraint} calls on the same
target are joined (max). Mostly useful when the model has a small
number of variables and reads naturally as a list of $\geq$
inequalities, mirroring the paper's \texttt{mcdp \{ ... \}} blocks.

\begin{example}[Drone as a single MCDP block]
The Fig.~48 drone rewritten as a flat MCDP, with one closed loop on
\texttt{battery\_mass}:
\begin{lstlisting}
with MCDP("drone") as m:
    m.provides("endurance",     unit="s")
    m.provides("extra_payload", unit="kg")
    m.provides("extra_power",   unit="W")
    m.provides("battery_mass",  unit="kg")   # loop axis: in F and R
    m.requires("battery_mass",  unit="kg")
    m.requires("report_mass",   unit="kg")   # mirror for visibility

    def battery_mass_eq(f):
        return (f["extra_power"]
                + 10.0 * (9.81 * (f["battery_mass"]
                                  + f["extra_payload"])) ** 2
                ) * f["endurance"] / 1.8e6

    m.constraint("battery_mass", battery_mass_eq)
    m.constraint("report_mass",  battery_mass_eq)
    m.loop_on("battery_mass")

drone = m.build()
\end{lstlisting}
Calling \texttt{solve(drone, \{"endurance": 300.0, "extra\_payload":
0.5, "extra\_power": 5.0\})} runs the same Kleene iteration as worked
example~1 and returns the same \texttt{report\_mass = 0.049 kg} fixed
point. The loop axis \texttt{battery\_mass} must be declared in both
\texttt{provides} and \texttt{requires}, since that is the name
\texttt{loop\_on} closes; the extra \texttt{report\_mass} resource
mirrors it into the outer $R$ so the converged value stays visible in
the result, following the guideline in Section~\ref{sec:guidelines}.
\end{example}

\section{The System Builder}
\label{sec:system-builder}

For larger designs the \texttt{System} class in \texttt{codesign.system}
provides modular composition with named subsystems and algebraic
connection constraints. This is the recommended way to build
non-trivial co-design models. Two equivalent surface syntaxes are
available: the operator-overloaded form (recommended for new code) and
the legacy lambda form (still supported).

\subsection{Operator-overloaded syntax (recommended)}

\texttt{provides}, \texttt{requires}, and \texttt{add} each return a
port handle. Arithmetic operators on handles build expression trees
lazily; the \texttt{>=} operator registers a constraint with the parent
system.

\begin{lstlisting}
from codesign import Module, Reals, System, solve

class Battery(Module):
    F = {"capacity": Reals(unit="J")}
    R = {"mass":     Reals(unit="kg")}
    def h(self, f):
        return {"mass": f["capacity"] / 1.8e6}

class Actuator(Module):
    F = {"lift_force": Reals(unit="N")}
    R = {"power":      Reals(unit="W")}
    def h(self, f):
        return {"power": 10.0 * f["lift_force"] ** 2}

sys = System("drone")
endurance     = sys.provides("endurance",     unit="s")
extra_payload = sys.provides("extra_payload", unit="kg")
extra_power   = sys.provides("extra_power",   unit="W")
total_mass    = sys.requires("total_mass",    unit="kg")

battery  = sys.add("battery",  Battery())
actuator = sys.add("actuator", Actuator())

battery.capacity    >= (actuator.power + extra_power) * endurance
actuator.lift_force >= 9.81 * (battery.mass + extra_payload)
total_mass          >= battery.mass + extra_payload

drone = sys.build()
\end{lstlisting}

The four port kinds are:

\begin{itemize}[leftmargin=2em]
  \item \textbf{Outer F} (from \texttt{provides}): can appear on the RHS
        of constraint expressions. Cannot be a constraint target.
  \item \textbf{Outer R} (from \texttt{requires}): can be a constraint
        target. Cannot appear in expression RHS.
  \item \textbf{Module F} (\texttt{handle.f\_port}): can be a constraint
        target. Cannot appear in expression RHS (its value is determined
        by the constraints on it, not by the current iteration state).
  \item \textbf{Module R} (\texttt{handle.r\_port}): can appear in
        expression RHS. Cannot be a constraint target (its value is
        determined by the module's own \texttt{h}, not externally).
\end{itemize}

Misusing a port (e.g.\ trying to put a module F port on the RHS, or
constraining an outer F) raises immediately with a message naming the
port and explaining the rule. F and R port names within a single
subsystem must be disjoint.

The supported operators are \texttt{+, -, *, /, **,} unary \texttt{-},
along with the helper functions \texttt{codesign.sqrt}, \texttt{exp},
and \texttt{log} for expressions that need transcendental terms. For
demands that cannot be expressed algebraically, the legacy lambda form
remains available and can be mixed freely with the operator form (see
below).

\subsection{Legacy lambda syntax}

\begin{lstlisting}
sys.constrain("battery.capacity",
              lambda x: (x["actuator.power"] + x["extra_power"]) * x["endurance"])
sys.constrain("actuator.lift_force",
              lambda x: 9.81 * (x["battery.mass"] + x["extra_payload"]))
sys.constrain("total_mass",
              lambda x: x["battery.mass"] + x["extra_payload"])
\end{lstlisting}

The argument \texttt{x} is a dict carrying outer F values under their
bare names (\texttt{x["endurance"]}) and subsystem R ports under their
dotted names (\texttt{x["battery.mass"]}). Both syntaxes compile to the
same internal representation; the operator form is purely sugar.

\subsection{Class-based modules}

The \texttt{Module} base class lets a design problem be defined
declaratively:

\begin{lstlisting}
class Battery(Module):
    F = {"capacity": Reals(unit="J")}
    R = {"mass":     Reals(unit="kg")}

    def __init__(self, specific_energy=1.8e6):
        self.specific_energy = specific_energy
        super().__init__()

    def h(self, f):
        return {"mass": f["capacity"] / self.specific_energy}
\end{lstlisting}

A \texttt{Module} subclass declares its $F$ and $R$ ports as
class-level dicts and overrides \texttt{h(self, f)}. The constructor
wires everything into the underlying \texttt{DesignProblem} machinery.
Parameterised modules are written by overriding \texttt{\_\_init\_\_}
and calling \texttt{super().\_\_init\_\_()} after setting any
instance attributes. \texttt{h} may return a dict (treated as a
singleton antichain), a list of dicts (multi-valued antichain), or
an \texttt{Antichain} directly.

\begin{example}[A modular drone end-to-end]
The drone from worked example~7, built with the operator-overloaded
constraint syntax. Each line below a \texttt{System} declaration is
a textbook inequality:
\begin{lstlisting}
sys = System("drone")
endurance     = sys.provides("endurance",     unit="s")
extra_payload = sys.provides("extra_payload", unit="kg")
extra_power   = sys.provides("extra_power",   unit="W")
total_mass    = sys.requires("total_mass",    unit="kg")

b = sys.add("battery",  Battery())
a = sys.add("actuator", Actuator())

b.capacity   >= (a.power + extra_power) * endurance
a.lift_force >= 9.81 * (b.mass + extra_payload)
total_mass   >= b.mass + extra_payload

drone = sys.build()
r = solve(drone, {"endurance": 300.0,
                  "extra_payload": 0.5,
                  "extra_power": 5.0})
# r.antichain -> Antichain[(total_mass=0.55 kg)]
\end{lstlisting}
\texttt{build()} bundles every subsystem's R into a single Kleene
axis and closes the loop automatically. The constraint graph can
be exported with \texttt{viz.to\_dot(drone)} for documentation.
\end{example}

\subsection{Semantics}

The \texttt{build} method emits a \texttt{Loop} whose axis bundles
every subsystem's resource poset:
\[
\texttt{\_\_modules\_\_} =
  \prod_{m \in \mathrm{modules}} R_m.
\]
The inner DP performs one Kleene step:
\begin{enumerate}[leftmargin=2em]
  \item Read the current bundle estimate to populate
        $\mathtt{ctx[m.r]}$ for every subsystem $R$ port.
  \item For each subsystem $m$, evaluate the constraint demands on
        its $F$ ports to obtain $f_m$, then call $m.h(f_m)$ to obtain
        an antichain.
  \item Take the Cartesian product across subsystems; for each
        combination, evaluate the outer-$R$ constraint demands to
        produce a candidate point.
  \item Return the resulting antichain.
\end{enumerate}
The outer Kleene iteration converges over all subsystem $R$ ports
simultaneously.

\subsection{Validation}

\texttt{build} performs the following checks:
\begin{itemize}[leftmargin=2em]
  \item every subsystem $F$ port has at least one \texttt{constrain}
        target;
  \item every outer $R$ has at least one \texttt{constrain} target;
  \item every constraint target refers to a known module $F$ port or a
        declared outer $R$;
  \item module names do not clash with the reserved axis name
        \texttt{\_\_modules\_\_} and contain no \texttt{.}.
\end{itemize}
If any check fails, \texttt{build} raises \texttt{ValueError} with a
message identifying the missing or invalid item.

\section{API Reference}
\label{sec:apiref}

This section is the complete, per-module reference for the public API
exported by \texttt{import codesign}. The six subsections below follow the
architecture of the library: order structures; design problems and their
solver; the two builders and the rendering helpers; the uncertainty and
online-learning layers; the temporal planners; and the vector-state and
closed-loop layers. Within each group one subsubsection documents each
source module, and each public symbol is a reference card giving its
purpose, signature, a minimal runnable snippet, its parameters, and the
exceptions it raises. Cards are cross-linked to the tutorial chapters
(Sections~\ref{sec:theory}--\ref{sec:domains}) and the worked examples
(Section~\ref{sec:examples}); the condensed export list in
Section~\ref{sec:api} remains a one-page index. All reference labels use
the \texttt{ref:<module>:<symbol>} scheme.

\subsection{Order structures}

The two foundational modules: \texttt{posets} supplies the functionality
and resource spaces, and \texttt{antichains} the Pareto-front value type
returned by every design problem. See Section~\ref{sec:theory} for the
underlying order theory.

\subsubsection{\texttt{posets}}
\label{ref:posets}

The \texttt{codesign.posets} module supplies the partially ordered sets
that serve as functionality spaces $F$ and resource spaces $R$. Every
concrete poset derives from the abstract \texttt{Poset} base and honours
its contract: a least element \texttt{bottom()} that seeds the Kleene
iteration, a greatest element \texttt{top()} that marks infeasibility,
and a \texttt{join()} that returns the least upper bound used to combine
constraints. The mathematical background is given in
Section~\ref{sec:theory}.

\paragraph{\texttt{Poset}}
\label{ref:posets:poset}
Abstract base class defining the poset contract. It is not instantiable
directly; a concrete subclass must override \texttt{leq}, \texttt{bottom}
and \texttt{top}, and inherits working defaults for the derived
operators. The four shipped subclasses (\texttt{Reals}, \texttt{Naturals},
\texttt{Ports}, \texttt{Discrete}) cover essentially every modelling need.

Signature: \texttt{class Poset(ABC)} with class attribute
\texttt{name: str}. The derived operators are exercised here through a
concrete subclass, since \texttt{Poset} itself cannot be instantiated:

\begin{lstlisting}
from codesign import Reals
P = Reals(unit="kg")
# Derived operators are provided by the Poset base class:
print(P.eq(1.0, 1.0))          # True
print(P.lt(1.0, 2.0))          # True
print(P.comparable(1.0, 2.0))  # True (a chain: all comparable)
print(P.is_bottom(0.0))        # True
\end{lstlisting}

The contract methods (the abstract three plus the inherited defaults):

\begin{center}
\begin{tabular}{@{}lll@{}}
\toprule
Method & Signature & Behaviour \\
\midrule
\texttt{leq}        & \texttt{leq(a, b) -> bool}   & abstract: is $a \leq b$? \\
\texttt{bottom}     & \texttt{bottom()}            & abstract: least element \\
\texttt{top}        & \texttt{top()}               & abstract: greatest element \\
\texttt{is\_top}    & \texttt{is\_top(x) -> bool}  & default via \texttt{leq} \\
\texttt{is\_bottom} & \texttt{is\_bottom(x) -> bool}& default via \texttt{leq} \\
\texttt{eq}         & \texttt{eq(a, b) -> bool}    & $a \leq b \wedge b \leq a$ \\
\texttt{lt}         & \texttt{lt(a, b) -> bool}    & $a \leq b \wedge a \neq b$ \\
\texttt{comparable} & \texttt{comparable(a, b) -> bool}& $a \leq b \vee b \leq a$ \\
\texttt{join}       & \texttt{join(a, b)}          & least upper bound \\
\texttt{format}     & \texttt{format(x) -> str}    & display string (\texttt{repr}) \\
\bottomrule
\end{tabular}
\end{center}

Parameters: the constructor takes no arguments; \texttt{name} is a
display label overridden by subclasses.

Errors: instantiating \texttt{Poset} directly raises \texttt{TypeError}
(abstract methods unimplemented). The default \texttt{join} raises
\texttt{NotImplementedError} when \texttt{a} and \texttt{b} are
incomparable (neither \texttt{leq} holds); chain and product subclasses
override \texttt{join} so it never fails on them.

\paragraph{\texttt{Reals}}
\label{ref:posets:reals}
The non-negative reals $\Rplus$ augmented with $+\infty$, forming a CPO.
This is the scalar chain for any real-valued port (mass, cost, energy).
Bottom is \texttt{0.0}, top is \texttt{math.inf}.

Signature: \texttt{Reals(unit: str = "", name: str = "R+")}.

\begin{lstlisting}
from codesign import Reals
R = Reals(unit="kg")
print(R.bottom(), R.top())   # 0.0 inf
print(R.leq(0.5, 1.2))       # True
print(R.join(0.5, 1.2))      # 1.2
print(R.format(0.56))        # 0.56 kg
print(R.format(R.top()))     # the top marker
print(R.name)                # R+[kg]
\end{lstlisting}

Parameters:
\begin{itemize}[nosep]
\item \texttt{unit} --- display-only unit string; travels through
  composition so \texttt{format} can print \texttt{0.56 kg}. Never used
  in any algebraic check.
\item \texttt{name} --- display name; if left at the default and a
  \texttt{unit} is given, becomes \texttt{R+[unit]} automatically.
\end{itemize}

Deviations from the \texttt{Poset} contract: overrides \texttt{is\_top}
(direct $+\infty$ test), \texttt{join} (returns \texttt{max}) and
\texttt{format} (prints $\topp$ for $+\infty$, else the value
with its unit).

Errors: none raised directly.

\paragraph{\texttt{Naturals}}
\label{ref:posets:naturals}
The non-negative integers $\Nat$ augmented with a top ($+\infty$), forming
a CPO. Mirrors \texttt{Reals} for integer-valued quantities such as part
counts or bit budgets. Bottom is \texttt{0}, top is \texttt{math.inf}.

Signature: \texttt{Naturals(unit: str = "", name: str = "N+")}.

\begin{lstlisting}
from codesign import Naturals
N = Naturals(unit="parts")
print(N.leq(3, N.top()))   # True
print(N.join(2, 5))        # 5
print(N.is_top(N.top()))   # True
print(N.bottom())          # 0
print(N.format(3))         # 3 parts
print(N.format(N.top()))   # ⊤
\end{lstlisting}

Parameters: \texttt{unit} and \texttt{name} as for \texttt{Reals}.
\texttt{format} mirrors \texttt{Reals}: it appends the \texttt{unit} (so
\texttt{N.format(3)} prints \texttt{3 parts}) and renders the top as
$\top$; finite values keep their integer form (\texttt{3.0} prints as
\texttt{3}). When \texttt{unit} is set and the default \texttt{name} is
kept, the name is decorated as \texttt{N+[unit]}.

Deviations from the \texttt{Poset} contract: overrides \texttt{leq} to
treat $+\infty$ as strictly above every finite natural (and equal only to
itself), \texttt{is\_top}, \texttt{join}, and \texttt{format}. Finite
inputs are coerced with \texttt{int()}.

Errors: none raised directly.

\paragraph{\texttt{Ports}}
\label{ref:posets:ports}
A named product poset: a typed bundle of named component posets, ordered
component-wise. This is the workhorse type; every design problem's $F$ and
$R$ is a \texttt{Ports}. Elements are Python dicts mapping each component
name to a value in its factor poset. The order is
$x \leq y \iff x[k] \leq y[k]$ for every component $k$; bottom and top are
the dicts of all component bottoms and tops.

Signature: \texttt{Ports(components: Mapping[str, Poset])}.

\begin{lstlisting}
from codesign import Ports, Reals
F = Ports({"capacity": Reals(unit="J")})
x = F.make(capacity=3.6e6)
print(F.leq(x, F.top()))    # True
print(F.is_top(F.top()))    # True
print(F.any_top({"capacity": F.top()["capacity"]}))  # True
print(F.format(x))          # (capacity=3.6e+06 J)
print(list(F.keys()))       # ['capacity']
\end{lstlisting}

Parameters:
\begin{itemize}[nosep]
\item \texttt{components} --- non-empty \texttt{\{name: poset\}} mapping.
  Copied into a fresh dict, so later mutation of the argument cannot
  corrupt the instance. A component poset may itself be a \texttt{Ports}
  (nesting is allowed).
\end{itemize}

Beyond the contract, \texttt{Ports} adds three public methods:

\begin{center}
\begin{tabular}{@{}lll@{}}
\toprule
Method & Signature & Behaviour \\
\midrule
\texttt{keys}    & \texttt{keys()}            & iterate component names \\
\texttt{any\_top}& \texttt{any\_top(x) -> bool}& \emph{any} component at top \\
\texttt{make}    & \texttt{make(**kwargs)}    & build a validated element dict \\
\bottomrule
\end{tabular}
\end{center}

Note the distinction between \texttt{is\_top} (\emph{every} component at
top --- total infeasibility) and \texttt{any\_top} (\emph{at least one}
component at top --- enough to disqualify a candidate in the antichain
machinery).

Errors: the constructor raises \texttt{ValueError} if \texttt{components}
is empty. \texttt{make} raises \texttt{ValueError} if any declared
component is missing from the keyword arguments, and \texttt{ValueError}
if any unknown component name is supplied. \texttt{leq}, \texttt{join} and
\texttt{format} may raise \texttt{KeyError} if an element dict lacks a
declared component.

\paragraph{\texttt{Discrete}}
\label{ref:posets:discrete}
A finite poset over an explicit element list with a caller-supplied order
predicate. Useful for enumerated design spaces (part numbers, operating
modes). The default order is equality: nothing is below anything but
itself.

Signature: \texttt{Discrete(elements, leq\_fn=None, name: str = "D")}.

\begin{lstlisting}
from codesign import Discrete
modes = Discrete(["idle", "cruise", "boost"], name="mode")
print(modes.leq("idle", "idle"))    # True
print(modes.leq("idle", "cruise"))  # False (equality order)
# bottom()/top() raise ValueError: no canonical extremum.
\end{lstlisting}

Parameters:
\begin{itemize}[nosep]
\item \texttt{elements} --- iterable of elements; stored as a list.
\item \texttt{leq\_fn} --- order predicate \texttt{leq\_fn(a, b) -> bool};
  defaults to \texttt{a == b} (the discrete order).
\item \texttt{name} --- display name.
\end{itemize}

Errors: \texttt{bottom()} and \texttt{top()} each raise
\texttt{ValueError}, because an arbitrary discrete order has no canonical
least or greatest element; subclass and override them, or wrap in a poset
that adds explicit extrema.

\paragraph{\texttt{NamedProduct}}
\label{ref:posets:namedproduct}
Backward-compatible alias: \texttt{NamedProduct = Ports}. Existing code
importing \texttt{NamedProduct} continues to work unchanged; new code
should prefer \texttt{Ports} (see~\ref{ref:posets:ports}). The two names
are the same object.

\begin{lstlisting}
from codesign import NamedProduct, Ports
print(NamedProduct is Ports)   # True
\end{lstlisting}

Errors: as for \texttt{Ports}.

\subsubsection{\texttt{antichains}}
\label{ref:antichains}

The \texttt{codesign.antichains} module provides the \texttt{Antichain}
type, the return value of every design problem's resource map
$h : F \to \Antichain[R]$. Antichains generalise the Pareto front to an
arbitrary partial order; their lattice structure (order, join by union
followed by $\Min$) is developed in Section~\ref{sec:theory}.

\paragraph{\texttt{Antichain}}
\label{ref:antichains:antichain}
A finite set of mutually incomparable elements of a poset, paired with the
poset they live in. Construction is self-normalising: dominated points and
duplicates are discarded, so the constructor is idempotent under $\Min$.
Instances are iterable and support \texttt{len()} and truthiness (empty is
falsey).

Signature: \texttt{Antichain(poset: Poset, points: Iterable = ())}.
The class-method constructors are the primary construction API; the bare
constructor is rarely called directly.

\begin{lstlisting}
from codesign import Antichain, Ports, Reals
R = Ports({"mass": Reals(unit="kg"), "cost": Reals(unit="$")})
pts = [{"mass": 1.0, "cost": 5.0},
       {"mass": 2.0, "cost": 3.0},
       {"mass": 1.5, "cost": 6.0}]   # dominated by the first
A = Antichain.from_set(R, pts)
print(len(A))    # 2 (dominated point dropped)
B = Antichain.singleton(R, {"mass": 0.5, "cost": 7.0})
U = Antichain.union_min(R, [A, B])
print(len(U))    # 3
\end{lstlisting}

Parameters:
\begin{itemize}[nosep]
\item \texttt{poset} --- the underlying \texttt{Poset} (typically a
  \texttt{Ports}) supplying \texttt{leq}, \texttt{eq} and, where present,
  \texttt{any\_top}.
\item \texttt{points} --- iterable of poset elements; normalised to an
  antichain on insertion.
\end{itemize}

Public methods (class-method constructors marked \emph{cm}):

\begin{center}
\begin{tabular}{@{}lll@{}}
\toprule
Method & Signature & Behaviour \\
\midrule
\texttt{of\_bottom}   & \texttt{of\_bottom(poset)} \emph{cm}      & singleton $\{\bot\}$; seeds Kleene \\
\texttt{empty}        & \texttt{empty(poset)} \emph{cm}           & the empty antichain \\
\texttt{singleton}    & \texttt{singleton(poset, point)} \emph{cm}& one-point antichain \\
\texttt{from\_set}    & \texttt{from\_set(poset, points)} \emph{cm}& $\Min$ of a subset \\
\texttt{union\_min}   & \texttt{union\_min(poset, antichains)} \emph{cm}& $\Min(\bigcup_i A_i)$ \\
\texttt{points}       & \texttt{points} (property)                & list copy of the points \\
\texttt{leq}          & \texttt{leq(other) -> bool}               & antichain order (see below) \\
\texttt{eq}           & \texttt{eq(other) -> bool}                & equal point sets under \texttt{eq} \\
\texttt{filter\_above}& \texttt{filter\_above(lower)}             & keep $x$ with $lower \leq x$ \\
\texttt{has\_any\_top}& \texttt{has\_any\_top() -> bool}          & any point at (a component) top \\
\texttt{is\_empty}    & \texttt{is\_empty() -> bool}              & no points \\
\bottomrule
\end{tabular}
\end{center}

The antichain order is $A \leq A'$ iff $\mathord{\uparrow} A \supseteq
\mathord{\uparrow} A'$: every point of \texttt{other} is dominated by some
point of \texttt{self}. \texttt{union\_min} and \texttt{filter\_above}
together implement the body of the Kleene iteration (Section~\ref{sec:theory});
\texttt{has\_any\_top} uses the poset's \texttt{any\_top} when available
(as on \texttt{Ports}), otherwise falls back to \texttt{is\_top}.

Errors: none raised directly. Operations delegate to the underlying
poset, so a malformed point may surface a \texttt{KeyError} from
\texttt{Ports.leq} rather than an antichain-level exception.

\subsection{Design problems and solving}

The six primitive design-problem types (\texttt{dp}), the scalar
\texttt{primitives} built on them, the three \texttt{composition}
operators, and the Kleene \texttt{solver} that computes fixed points. The
tutorial treatment is in Sections~\ref{sec:primitives}--\ref{sec:solver}.

\subsubsection{\texttt{dp}}
\label{ref:dp}

The \texttt{dp} module defines the central abstraction, the
\texttt{DesignProblem}, together with the concrete primitives that
realise it in closed form, from a catalogue, from a constraint, or from
an ODE. Every design problem exposes a single method \texttt{h(f)}
returning the antichain of minimal resources able to deliver
functionality \texttt{f}.

\paragraph{\texttt{DesignProblem}}
\label{ref:dp:designproblem}
The abstract base class for all design problems: a monotone relation
from a functionality poset \texttt{F} to an antichain in a resource
poset \texttt{R}. Subclasses must implement \texttt{h}; the class is
never instantiated directly.
\begin{lstlisting}
class DesignProblem(ABC):
    F: Poset            # functionality space
    R: Poset            # resource space
    name: str = "dp"
    def h(self, f) -> Antichain: ...   # abstract
\end{lstlisting}
Concrete instances expose \texttt{h}, \texttt{F}, \texttt{R} and a
readable \texttt{\_\_repr\_\_}:
\begin{lstlisting}
from codesign import AlgebraicDP, Ports, Reals
dp = AlgebraicDP(
    F=Ports({"speed": Reals()}),
    R=Ports({"power": Reals(unit="W")}),
    equations={"power": lambda f: 2.0 * f["speed"]},
)
print(repr(dp))              # DP(algebraic: speed:R+ -> A[...])
print(dp.h({"speed": 3.0}))  # Antichain[(power=6 W)]
\end{lstlisting}
\begin{itemize}[nosep]
  \item \texttt{F} --- \texttt{Poset} --- the functionality space.
  \item \texttt{R} --- \texttt{Poset} --- the resource space.
  \item \texttt{h(f)} --- returns the minimal-resource
        \texttt{Antichain} for \texttt{f}.
\end{itemize}
Errors: being an \texttt{ABC} with an abstract \texttt{h}, direct
instantiation raises \texttt{TypeError}; subclasses raise their own
(below).

\paragraph{\texttt{AlgebraicDP}}
\label{ref:dp:algebraicdp}
A functional design problem whose every resource is a closed-form
(monotone) expression in the functionality. The relation is a
singleton antichain. See Section~\ref{sec:primitives} for the tutorial
treatment. Monotonicity is the caller's responsibility and is not
verified.
\begin{lstlisting}
AlgebraicDP(F: Poset, R: Ports,
            equations: Mapping[str, Callable | Any],
            name: str = "algebraic")
\end{lstlisting}
\begin{lstlisting}
from codesign import AlgebraicDP, Ports, Reals
motor = AlgebraicDP(
    F=Ports({"torque": Reals(unit="Nm")}),
    R=Ports({"mass": Reals(unit="kg")}),
    equations={"mass": lambda f: 0.5 * f["torque"] + 1.0},
)
print(motor.h({"torque": 4.0}))  # Antichain[(mass=3 kg)]
\end{lstlisting}
\begin{itemize}[nosep]
  \item \texttt{R} --- \texttt{Ports} --- named resource bundle
        (required).
  \item \texttt{equations} --- one entry per resource port, a callable
        \texttt{f\_dict -> value} or a constant.
\end{itemize}
Errors: raises \texttt{TypeError} if \texttt{R} is not a \texttt{Ports};
raises \texttt{ValueError} if any resource port declared in \texttt{R}
has no entry in \texttt{equations}.

\paragraph{\texttt{FunctionDP}}
\label{ref:dp:functiondp}
Wraps an arbitrary user callable \texttt{h\_fn}. The return value is
normalised to an antichain: an \texttt{Antichain} is passed through, a
\texttt{dict} or a bare point becomes a singleton, and any other
iterable is treated as a set of points.
\begin{lstlisting}
FunctionDP(F: Poset, R: Poset,
           h_fn: Callable[[Any], Any],
           name: str = "function")
\end{lstlisting}
\begin{lstlisting}
from codesign import FunctionDP, Ports, Reals
fdp = FunctionDP(
    F=Ports({"load": Reals()}),
    R=Ports({"cost": Reals()}),
    h_fn=lambda f: {"cost": f["load"] ** 2},
)
print(fdp.h({"load": 3.0}))  # Antichain[(cost=9)]
\end{lstlisting}
\begin{itemize}[nosep]
  \item \texttt{h\_fn} --- callable returning an \texttt{Antichain}, a
        \texttt{dict}, a point, or an iterable of points.
\end{itemize}
Errors: none raised directly; the caller guarantees monotonicity, and
any exception raised inside \texttt{h\_fn} propagates unchanged.

\paragraph{\texttt{CatalogEntry}}
\label{ref:dp:catalogentry}
A dataclass describing one implementation in a catalogue: the
functionality it can deliver and the resources it costs.
\begin{lstlisting}
from codesign import CatalogEntry
e = CatalogEntry(provides={"speed": 5.0},
                 costs={"cost": 10.0}, name="A")
\end{lstlisting}
\begin{center}
\begin{tabular}{@{}lll@{}}
\toprule
Field & Type & Meaning \\
\midrule
\texttt{provides} & \texttt{Mapping} & functionality this entry fulfils
  (any \texttt{f <= provides}) \\
\texttt{costs}    & \texttt{Mapping} & resource vector it costs \\
\texttt{name}     & \texttt{str}     & optional label (default \texttt{""}) \\
\bottomrule
\end{tabular}
\end{center}
Errors: none raised directly.

\paragraph{\texttt{CatalogDP}}
\label{ref:dp:catalogdp}
Selects the cheapest catalogue entries able to deliver \texttt{f}. An
entry is feasible when its \texttt{provides} dominates \texttt{f}
port-wise; the result is the \texttt{Min} of the feasible costs. If no
entry works, \texttt{h} returns the singleton at $\topp$ (infeasible),
rather than raising.
\begin{lstlisting}
CatalogDP(F: Ports, R: Ports, catalog: Sequence,
          name: str = "catalog")
\end{lstlisting}
\begin{lstlisting}
from codesign import CatalogDP, CatalogEntry, Ports, Reals
cat = CatalogDP(
    Ports({"speed": Reals()}), Ports({"cost": Reals()}),
    catalog=[
        CatalogEntry({"speed": 5.0}, {"cost": 10.0}, "A"),
        CatalogEntry({"speed": 9.0}, {"cost": 25.0}, "B"),
    ])
print(cat.h({"speed": 6.0}))  # only B works: (cost=25)
\end{lstlisting}
\begin{itemize}[nosep]
  \item \texttt{catalog} --- a sequence of \texttt{CatalogEntry} objects
        or plain dicts with \texttt{provides}/\texttt{costs}/\texttt{name}
        keys.
\end{itemize}
Errors: raises \texttt{TypeError} if \texttt{F} or \texttt{R} is not a
\texttt{Ports}; \texttt{ValueError} (naming the DP) if the catalogue is
empty, or if two entries share the same non-empty \texttt{name} (the
clashing name(s) are listed). Multiple entries left unnamed (the default
\texttt{name=""}) are allowed.

\paragraph{\texttt{ConstraintDP}}
\label{ref:dp:constraintdp}
Defines a problem by a \texttt{sampler} over an implementation space, a
\texttt{feasible} predicate, and a \texttt{cost} function; \texttt{h}
returns the \texttt{Min} over all feasible samples (or the singleton at
$\topp$ when none is feasible).
\begin{lstlisting}
ConstraintDP(F: Poset, R: Ports,
             sampler: Callable[[Any], Iterable],
             feasible: Callable[[Any, Any], bool],
             cost: Callable[[Any], Mapping],
             name: str = "constraint")
\end{lstlisting}
\begin{lstlisting}
from codesign import ConstraintDP, Ports, Reals
cdp = ConstraintDP(
    F=Ports({"span": Reals()}), R=Ports({"area": Reals()}),
    sampler=lambda f: [0.5, 1.0, 2.0],
    feasible=lambda impl, f: impl >= f["span"],
    cost=lambda impl: {"area": impl},
)
print(cdp.h({"span": 0.8}))  # Antichain[(area=1)]
\end{lstlisting}
\begin{itemize}[nosep]
  \item \texttt{sampler(f)} --- yields candidate implementations.
  \item \texttt{feasible(impl, f)} --- whether \texttt{impl} delivers
        \texttt{f}.
  \item \texttt{cost(impl)} --- resource vector for a feasible sample.
\end{itemize}
Errors: none raised directly (no upfront interface check); exceptions
inside the callables propagate.

\paragraph{\texttt{ODE\_DP}}
\label{ref:dp:ode_dp}
Derives a monotone relation from an ODE $dx/dt = \texttt{rhs}(x, t, f)$,
reading resources from the trajectory with \texttt{extract}. In
\texttt{final\_value} mode it integrates (explicit Euler) to
\texttt{t\_end}; in \texttt{steady\_state} mode it Newton-solves
$\texttt{rhs}(x)=0$. See Section~\ref{sec:solver} for the wider
solve pipeline.
\begin{lstlisting}
ODE_DP(F: Poset, R: Ports,
       rhs: Callable[[Any, float, Any], Any],
       extract: Callable[[Any], Mapping],
       mode: str = "final_value",
       t_end: float = 10.0, n_steps: int = 200,
       x0_fn: Callable[[Any], Any] | None = None,
       name: str = "ode")
\end{lstlisting}
\begin{lstlisting}
from codesign import ODE_DP, Ports, Reals
tank = ODE_DP(
    F=Ports({"inflow": Reals()}), R=Ports({"level": Reals()}),
    rhs=lambda x, t, f: f["inflow"] - 0.5 * x,
    extract=lambda x: {"level": x},
    mode="steady_state", x0_fn=lambda f: 0.0,
)
print(tank.h({"inflow": 2.0}))  # Antichain[(level=4)]
\end{lstlisting}
\begin{itemize}[nosep]
  \item \texttt{mode} --- \texttt{"final\_value"} or
        \texttt{"steady\_state"}.
  \item \texttt{x0\_fn(f)} --- initial state; a scalar or a positional
        sequence of floats (defaults to \texttt{0.0}).
  \item \texttt{t\_end}, \texttt{n\_steps} --- integration horizon and
        step count (\texttt{final\_value} only).
\end{itemize}
Errors: the constructor raises \texttt{ValueError} if \texttt{mode} is
not \texttt{"final\_value"} or \texttt{"steady\_state"}. On the first
\texttt{h} call the state is validated: a \texttt{dict}-valued (named)
state raises \texttt{TypeError} (name the keys onto a list and index
positionally), as does a \texttt{str}/\texttt{bytes} state.

\paragraph{\texttt{UncertainDP}}
\label{ref:dp:uncertaindp}
Brackets a nominal problem between an optimistic lower DP \texttt{lower}
and a pessimistic upper DP \texttt{upper}, returning whichever the
\texttt{mode} selects (paper Sec.~VII). See
Section~\ref{sec:uncertainty} for the uncertainty solver.
\begin{lstlisting}
UncertainDP(F: Poset, R: Poset,
            lower: DesignProblem, upper: DesignProblem,
            mode: str = "upper", name: str = "uncertain")
\end{lstlisting}
\begin{lstlisting}
from codesign import AlgebraicDP, UncertainDP, Ports, Reals
F = Ports({"x": Reals()}); R = Ports({"y": Reals()})
lo = AlgebraicDP(F, R, {"y": lambda f: f["x"]})
hi = AlgebraicDP(F, R, {"y": lambda f: 1.5 * f["x"]})
u = UncertainDP(F, R, lower=lo, upper=hi, mode="upper")
print(u.h({"x": 2.0}))                     # (y=3)
print(u.with_mode("lower").h({"x": 2.0}))  # (y=2)
\end{lstlisting}
\begin{itemize}[nosep]
  \item \texttt{lower}, \texttt{upper} --- the two bounding DPs.
  \item \texttt{mode} --- \texttt{"lower"} (optimistic) or
        \texttt{"upper"} (pessimistic).
\end{itemize}
\subparagraph{\texttt{with\_mode(mode)}}
\label{ref:dp:uncertaindp:with_mode}
Returns a fresh \texttt{UncertainDP} sharing the same bounds but with the
given \texttt{mode}; the convenient way to obtain both fronts.
\newline Errors: the constructor (and hence \texttt{with\_mode}) raises
\texttt{ValueError} if \texttt{mode} is not \texttt{"lower"} or
\texttt{"upper"}.

\subsubsection{\texttt{primitives}}
\label{ref:primitives}

The \texttt{primitives} module provides small factory functions for the
everyday plumbing DPs used to glue diagrams together (paper Fig.~26,
Fig.~35). Each returns a ready-made \texttt{AlgebraicDP}
(\ref{ref:dp:algebraicdp}) over scalar \texttt{Reals} ports, so its
result is always a singleton antichain. Every factory accepts an
optional \texttt{poset} argument; when omitted a plain \texttt{Reals()}
is used for all ports. As thin wrappers over \texttt{AlgebraicDP}, they
raise nothing directly; the ports they build are always \texttt{Ports},
so the wrapped constructor never faults.

\paragraph{\texttt{adder}}
\label{ref:primitives:adder}
Sums several scalar inputs into one scalar output.
\begin{lstlisting}
adder(in_names: list[str], out_name: str, poset=None)
\end{lstlisting}
\begin{lstlisting}
from codesign import adder
a = adder(["p1", "p2"], "total")
print(a.h({"p1": 2.0, "p2": 3.0}))  # Antichain[(total=5)]
\end{lstlisting}

\paragraph{\texttt{multiplier}}
\label{ref:primitives:multiplier}
Multiplies two scalar inputs, e.g. current times voltage gives power.
\begin{lstlisting}
multiplier(in_a: str, in_b: str, out_name: str, poset=None)
\end{lstlisting}
\begin{lstlisting}
from codesign import multiplier
m = multiplier("current", "voltage", "power")
print(m.h({"current": 2.0, "voltage": 5.0}))  # (power=10)
\end{lstlisting}

\paragraph{\texttt{scale}}
\label{ref:primitives:scale}
Multiplies one input by a constant \texttt{factor}.
\begin{lstlisting}
scale(in_name: str, out_name: str, factor: float, poset=None)
\end{lstlisting}
\begin{lstlisting}
from codesign import scale
s = scale("capacity", "mass", factor=0.02)
print(s.h({"capacity": 1000.0}))  # Antichain[(mass=20)]
\end{lstlisting}

\paragraph{\texttt{constant}}
\label{ref:primitives:constant}
Ignores its (trivial) functionality port \texttt{"\_"} and emits a fixed
\texttt{value}.
\begin{lstlisting}
constant(out_name: str, value: float, poset=None)
\end{lstlisting}
\begin{lstlisting}
from codesign import constant
c = constant("base_cost", value=42.0)
print(c.h({"_": 0.0}))  # Antichain[(base_cost=42)]
\end{lstlisting}

\paragraph{\texttt{identity}}
\label{ref:primitives:identity}
Passes a single named scalar through unchanged (same port name on both
sides).
\begin{lstlisting}
identity(name: str, poset=None)
\end{lstlisting}
\begin{lstlisting}
from codesign import identity
i = identity("speed")
print(i.h({"speed": 7.0}))  # Antichain[(speed=7)]
\end{lstlisting}

\subsubsection{\texttt{composition}}
\label{ref:composition}

The \texttt{composition} module supplies the three operators that close
design problems under composition (Section~\ref{sec:composition}):
\texttt{Series}, \texttt{Parallel} and \texttt{Loop}, each with a
lowercase functional alias (\texttt{series}, \texttt{par}, \texttt{loop})
matching the paper's notation. All three preserve monotonicity, so a
composite is itself a \texttt{DesignProblem} and solves through the same
\texttt{solve} entry point (\ref{ref:solver:solve}).

\paragraph{\texttt{Series}}
\label{ref:composition:series}
Series composition: \texttt{dp1} then \texttt{dp2}, with each resource of
\texttt{dp1} fed in as a functionality of \texttt{dp2}. The composite has
\texttt{F = dp1.F} and \texttt{R = dp2.R}. When both interfaces are
\texttt{Ports}, the connectable check
\texttt{set(dp2.F.keys()) <= set(dp1.R.keys())} is enforced upfront;
extra resource ports on \texttt{dp1} are permitted and ignored.
\begin{lstlisting}
Series(dp1: DesignProblem, dp2: DesignProblem,
       name: str | None = None)
# alias: series(dp1, dp2, name=None) -> Series
\end{lstlisting}
\begin{lstlisting}
from codesign import series, scale
s1 = scale("torque", "current", factor=2.0)
s2 = scale("current", "mass", factor=0.5)
chain = series(s1, s2)
print(chain.h({"torque": 4.0}))  # Antichain[(mass=4)]
\end{lstlisting}
\begin{itemize}[nosep]
  \item \texttt{dp1}, \texttt{dp2} --- the two stages; \texttt{dp1}'s
        resources must cover \texttt{dp2}'s functionality ports.
\end{itemize}
Errors: raises \texttt{ValueError} (naming the missing ports) when both
interfaces are \texttt{Ports} and some port \texttt{dp2} requires is not
produced by \texttt{dp1}.

\paragraph{\texttt{Parallel}}
\label{ref:composition:parallel}
Parallel composition: \texttt{dp1} and \texttt{dp2} stacked side by side.
The composite has \texttt{F = F1 x F2} and \texttt{R = R1 x R2} (Ports
concatenation), and its antichain is the Cartesian product of the two.
\begin{lstlisting}
Parallel(dp1: DesignProblem, dp2: DesignProblem,
         name: str | None = None)
# alias: par(dp1, dp2, name=None) -> Parallel
\end{lstlisting}
\begin{lstlisting}
from codesign import par, scale
p = par(scale("x", "cost_x", 1.0),
        scale("y", "cost_y", 2.0))
print(p.h({"x": 3.0, "y": 4.0}))
# Antichain[(cost_x=3, cost_y=8)]
\end{lstlisting}
\begin{itemize}[nosep]
  \item \texttt{dp1}, \texttt{dp2} --- must have disjoint \texttt{F} port
        names and disjoint \texttt{R} port names.
\end{itemize}
Errors: raises \texttt{TypeError} if either \texttt{F} or either
\texttt{R} is not a \texttt{Ports}; raises \texttt{ValueError} if the
\texttt{F} port names overlap, or if the \texttt{R} port names overlap.

\paragraph{\texttt{Loop}}
\label{ref:composition:loop}
Feedback composition: closes a named \texttt{axis} present on both
\texttt{inner.F} and \texttt{inner.R}, feeding the produced resource back
as the consumed functionality (paper Def.~16, Section~\ref{sec:loop}).
The outer \texttt{F} and \texttt{R} drop the looped axis. Evaluating
\texttt{h} defers to the solver's Kleene iteration
(\ref{ref:solver:kleene_loop}).
\begin{lstlisting}
Loop(inner: DesignProblem, axis: str,
     name: str | None = None)
# alias: loop(inner, axis, name=None) -> Loop
\end{lstlisting}
\begin{lstlisting}
from codesign import loop, AlgebraicDP, solve, Ports, Reals
inner = AlgebraicDP(
    F=Ports({"payload": Reals(), "mass": Reals()}),
    R=Ports({"mass": Reals(), "power": Reals(unit="W")}),
    equations={
        "mass": lambda f: 0.5 * f["payload"] + 0.3 * f["mass"],
        "power": lambda f: 2.0 * f["mass"],
    })
lp = loop(inner, axis="mass")     # feed R.mass back to F.mass
print(solve(lp, {"payload": 2.0}).antichain)
# Antichain[(power=2.857 W)]
\end{lstlisting}
\begin{itemize}[nosep]
  \item \texttt{inner} --- the DP to close; needs \texttt{Ports} on both
        sides.
  \item \texttt{axis} --- the port name (on both \texttt{F} and
        \texttt{R}) to feed back.
\end{itemize}
Errors: raises \texttt{TypeError} if \texttt{inner.F} or \texttt{inner.R}
is not a \texttt{Ports}; raises \texttt{ValueError} if \texttt{axis} is
not declared on both \texttt{F} and \texttt{R} of \texttt{inner}.

\subsubsection{\texttt{solver}}
\label{ref:solver}

The \texttt{solver} module drives every design problem through one entry
point, \texttt{solve}, which runs the Kleene fixed-point iteration
whenever a loop is present and returns a \texttt{SolveResult}. See
Section~\ref{sec:solver} for the tutorial treatment and
Section~\ref{sec:solver-observability} for the observability features.

\paragraph{\texttt{solve}}
\label{ref:solver:solve}
Solves a (possibly composite) design problem at a given functionality.
For a \texttt{Loop} it iterates to the least fixed point; for any other
DP it is a single evaluation of \texttt{h}. This is the single
most-used function in the library.
\begin{lstlisting}
solve(
    dp: DesignProblem,
    functionality: Mapping | None = None,
    max_iter: int = 200,
    *,
    trace: bool = False,
    verbose: int = 0,
    on_iteration: Callable[[TraceEntry], None] | None = None,
    start_from: SolveResult | Antichain | None = None,
    record_trace: bool = False,
    uncertainty: list[str] | None = None,
    n_samples: int = 1000,
    rng_seed: int | None = None,
) -> SolveResult
\end{lstlisting}
\begin{lstlisting}
from codesign import AlgebraicDP, solve, minimize_cost, Ports, Reals
dp = AlgebraicDP(
    F=Ports({"speed": Reals()}),
    R=Ports({"power": Reals(unit="W")}),
    equations={"power": lambda f: 2.0 * f["speed"]},
)
res = solve(dp, {"speed": 10.0}, trace=True)
print(res.status, res.iterations, res.feasible)
# converged 0 True
print(res.antichain)  # Antichain[(power=20 W)]
\end{lstlisting}
\begin{itemize}[nosep]
  \item \texttt{functionality} --- outer \texttt{F} values; when
        \texttt{None} the solver uses \texttt{dp.F.bottom()}.
  \item \texttt{max\_iter} --- cap on Kleene steps; reaching it sets
        \texttt{status="max\_iter"}.
  \item \texttt{trace} --- collect a per-step trace on
        \texttt{result.trace} (\ref{ref:solver:traceentry}).
  \item \texttt{verbose} --- \texttt{0} silent, \texttt{1} final summary,
        \texttt{2} per-iteration feed.
  \item \texttt{on\_iteration} --- callback invoked with each
        \texttt{TraceEntry} as it is produced.
  \item \texttt{start\_from} --- warm-start seed: a prior
        \texttt{SolveResult} (its inner antichain is reused) or an
        \texttt{Antichain}; ignored for non-loop DPs.
  \item \texttt{record\_trace} --- legacy alias for \texttt{trace}.
  \item \texttt{uncertainty} --- if given (labels
        \texttt{"worst\_case"}, \texttt{"mean"}, \texttt{"p95"},
        \texttt{"cvar95"}, \texttt{"samples"}), dispatches to the
        uncertainty solver (Section~\ref{sec:uncertainty}) and returns an
        \texttt{UncertaintyResult}.
  \item \texttt{n\_samples}, \texttt{rng\_seed} --- Monte Carlo sample
        count and optional seed for stochastic uncertainty.
\end{itemize}
Errors: raises \texttt{TypeError} if \texttt{start\_from} is not a
\texttt{SolveResult}, an \texttt{Antichain}, or \texttt{None}.
Non-convergence is \emph{not} an error: it is reported through
\texttt{status} (\texttt{"max\_iter"} or \texttt{"diverged"}) and
\texttt{feasible}, so callers must inspect those fields.

\paragraph{\texttt{kleene\_loop}}
\label{ref:solver:kleene_loop}
The raw fixed-point iterator behind \texttt{solve}: it returns the loop's
antichain directly (without the \texttt{SolveResult} wrapper), writing
metadata into an optional \texttt{info\_out} dict. Most callers should
prefer \texttt{solve}.
\begin{lstlisting}
kleene_loop(
    loop_dp: Loop,
    f_outer: Mapping | None,
    max_iter: int = 200,
    *,
    trace: bool = False, verbose: int = 0,
    on_iteration: Callable | None = None,
    info_out: dict | None = None,
    start_from: Antichain | None = None,
) -> Antichain
\end{lstlisting}
\begin{lstlisting}
from codesign import loop, kleene_loop, AlgebraicDP, Ports, Reals
inner = AlgebraicDP(
    F=Ports({"payload": Reals(), "mass": Reals()}),
    R=Ports({"mass": Reals(), "power": Reals(unit="W")}),
    equations={
        "mass": lambda f: 0.5 * f["payload"] + 0.3 * f["mass"],
        "power": lambda f: 2.0 * f["mass"],
    })
lp = loop(inner, axis="mass")
print(kleene_loop(lp, {"payload": 2.0}))
# Antichain[(power=2.857 W)]
\end{lstlisting}
\begin{itemize}[nosep]
  \item \texttt{info\_out} --- receives \texttt{iterations},
        \texttt{status}, \texttt{inner\_antichain}, and (if
        \texttt{trace}) \texttt{trace}.
  \item \texttt{start\_from} --- seed antichain in the inner \texttt{R}
        poset; a cold start (bottom singleton) is used when \texttt{None}.
\end{itemize}
Errors: raises \texttt{TypeError} if \texttt{start\_from} is not
\texttt{None} and lacks a \texttt{points} attribute (i.e.\ is not an
\texttt{Antichain}). Runaway numeric values are capped at
\texttt{DIVERGENCE\_CAP} ($10^{30}$) and reported as
\texttt{status="diverged"} rather than raised.

\paragraph{\texttt{minimize\_cost}}
\label{ref:solver:minimize_cost}
Scalarises a solved antichain: returns the single resource bundle
minimising a scalar \texttt{cost\_fn}, or \texttt{None} when the result
is infeasible or empty. This is the step that turns a Pareto front into
one engineering choice.
\begin{lstlisting}
minimize_cost(result: SolveResult,
              cost_fn: Callable[[Mapping], float]) -> Mapping | None
\end{lstlisting}
\begin{lstlisting}
from codesign import AlgebraicDP, solve, minimize_cost, Ports, Reals
dp = AlgebraicDP(Ports({"speed": Reals()}),
                 Ports({"power": Reals(unit="W")}),
                 {"power": lambda f: 2.0 * f["speed"]})
res = solve(dp, {"speed": 10.0})
print(minimize_cost(res, cost_fn=lambda r: r["power"]))
# {'power': 20.0}
\end{lstlisting}
\begin{itemize}[nosep]
  \item \texttt{result} --- a \texttt{SolveResult}, or any object with
        \texttt{feasible} and \texttt{antichain} attributes.
  \item \texttt{cost\_fn(r)} --- scalar cost of a resource bundle.
\end{itemize}
Errors: none raised directly; returns \texttt{None} on an infeasible or
empty result.

\paragraph{\texttt{SolveResult}}
\label{ref:solver:solveresult}
The dataclass returned by \texttt{solve}, carrying the front plus
convergence and feasibility metadata.
\begin{lstlisting}
from codesign import SolveResult, AlgebraicDP, solve, Ports, Reals
dp = AlgebraicDP(Ports({"x": Reals()}), Ports({"y": Reals()}),
                 {"y": lambda f: f["x"]})
ac = solve(dp, {"x": 1.0}).antichain
sr = SolveResult(antichain=ac, iterations=3, status="max_iter")
print(sr.converged)  # False
\end{lstlisting}
\begin{center}
\begin{tabular}{@{}lll@{}}
\toprule
Field & Type & Meaning \\
\midrule
\texttt{antichain} & \texttt{Antichain} & the (approximate) least fixed
  point \\
\texttt{iterations} & \texttt{int} & Kleene steps taken (0 if no loop) \\
\texttt{status} & \texttt{str} & \texttt{"converged"}, \texttt{"max\_iter"}
  or \texttt{"diverged"} \\
\texttt{feasible} & \texttt{bool} & \texttt{False} iff every point hit
  $\topp$ \\
\texttt{trace} & \texttt{list | None} & per-step trace if
  \texttt{trace=True}, else \texttt{None} \\
\bottomrule
\end{tabular}
\end{center}
The read-only property \texttt{converged} is a
backward-compatibility alias for \texttt{status == "converged"}.
Errors: none raised directly.

\paragraph{\texttt{TraceEntry}}
\label{ref:solver:traceentry}
A dataclass recording one step of a Kleene iteration; a list of these is
placed on \texttt{result.trace} when \texttt{trace=True}.
\begin{lstlisting}
from codesign import solve, AlgebraicDP, Ports, Reals
dp = AlgebraicDP(Ports({"speed": Reals()}),
                 Ports({"power": Reals()}),
                 {"power": lambda f: 2.0 * f["speed"]})
te = solve(dp, {"speed": 10.0}, trace=True).trace[-1]
print(te.iteration, te.n_points, te.delta)  # 1 1 None
\end{lstlisting}
\begin{center}
\begin{tabular}{@{}lll@{}}
\toprule
Field & Type & Meaning \\
\midrule
\texttt{iteration} & \texttt{int} & 0 = seed, 1 = after first step, \ldots \\
\texttt{antichain} & \texttt{Antichain} & snapshot at this step \\
\texttt{n\_points} & \texttt{int} & \texttt{len(antichain)} \\
\texttt{delta} & \texttt{float | None} & max port change vs.\ the previous
  step; \texttt{None} at iteration 0 \\
\texttt{elapsed\_ms} & \texttt{float} & wall time for this step, in ms \\
\bottomrule
\end{tabular}
\end{center}
Errors: none raised directly.

\subsection{Builders and rendering}

The two high-level builders --- the declarative \texttt{mcdpl} front end
and the modular \texttt{system} builder, with its \texttt{module} base
class and \texttt{sugar} constraint DSL --- together with the
\texttt{diagram} and \texttt{viz} rendering helpers. See
Sections~\ref{sec:mcdp-builder}--\ref{sec:system-builder}
and~\ref{sec:visualisation}.

\subsubsection{\texttt{mcdpl}}
\label{ref:mcdpl}

The \texttt{codesign.mcdpl} module supplies a single builder class,
\texttt{MCDP}, a thin declarative front-end that mirrors the paper's
\texttt{mcdp \{ ... \}} concrete syntax in pure Python. The tutorial in
Section~\ref{sec:mcdp-builder} walks through the surface syntax; the card
below is the reference summary.

\paragraph{\texttt{MCDP}}\label{ref:mcdpl:mcdp}
Fluent builder for a single self-contained monotone design problem. One
declares functionalities with \texttt{provides} and resources with
\texttt{requires}, attaches one equation per resource with
\texttt{constraint} (or a full antichain-valued \texttt{rule}), optionally
closes feedback axes with \texttt{loop\_on}, and emits a plain
\texttt{DesignProblem} with \texttt{build}. Every declaration method
returns \texttt{self}, so calls chain. It also works as a context manager.

Signature: \texttt{MCDP(name="mcdp")}.

\begin{lstlisting}
from codesign import MCDP, solve

with MCDP("battery") as m:
    m.provides("capacity", unit="J")
    m.requires("mass", unit="kg")
    m.constraint("mass", lambda f: f["capacity"] / 1.8e6)
battery = m.build()

res = solve(battery, {"capacity": 3.6e6})
print(res.antichain)          # Antichain[(mass=2 kg)]
\end{lstlisting}

Constructor parameters:
\begin{itemize}[nosep]
  \item \texttt{name} --- \texttt{str}, optional (default \texttt{"mcdp"}).
        Names the emitted \texttt{DesignProblem} and appears in its errors.
\end{itemize}

Public methods:
\begin{center}
\begin{tabular}{@{}p{0.30\linewidth}p{0.62\linewidth}@{}}
\toprule
\textbf{Method} & \textbf{Behaviour} \\
\midrule
\texttt{provides(name, *, unit="", poset=None)} &
  Declare a functionality port. Defaults to \texttt{Reals(unit=unit)};
  pass \texttt{poset} to override. Returns \texttt{self}. \\
\texttt{requires(name, *, unit="", poset=None)} &
  Declare a resource port, same defaulting as \texttt{provides}. Returns
  \texttt{self}. \\
\texttt{constraint(resource, fn)} &
  Closed-form clause \texttt{resource >= fn(f)} where \texttt{f} is the
  functionality dict. Repeated calls on one resource are joined (max),
  matching MCDPL's multiple-\texttt{>=} semantics. Returns \texttt{self}. \\
\texttt{rule(fn)} &
  Hand-write the full \texttt{h: f -> Antichain}. Overrides all
  \texttt{constraint} clauses; use for multi-valued or branchy relations.
  Returns \texttt{self}. \\
\texttt{loop\_on(axis)} &
  Close a feedback loop on \texttt{axis}, which must be declared on both
  the \texttt{provides} and \texttt{requires} sides; the axis is projected
  out of the closed loop's $F$/$R$. May be repeated. Returns \texttt{self}. \\
\texttt{build()} &
  Emit the plain \texttt{DesignProblem} (an \texttt{AlgebraicDP}, or a
  \texttt{FunctionDP} if \texttt{rule} was used), wrapped in a
  \texttt{Loop} per declared axis. \\
\bottomrule
\end{tabular}
\end{center}

Errors: \texttt{loop\_on} raises \texttt{ValueError} naming the missing
side and listing the current \texttt{provides}/\texttt{requires}
declarations if \texttt{axis} is absent from either side. \texttt{build}
raises \texttt{ValueError} if no functionality was declared, if no
resource was declared, or (in the closed-form path) if any required
resource has neither a \texttt{constraint} nor a \texttt{rule}. The other
methods raise nothing directly.

\subsubsection{\texttt{system}}
\label{ref:system}

The \texttt{codesign.system} module supplies the \texttt{System} builder,
the recommended surface for non-trivial models. Section~\ref{sec:system-builder}
is the full tutorial (both the operator-overloaded and legacy-lambda
syntaxes); the card below summarises the API.

\paragraph{\texttt{System}}\label{ref:system:system}
Modular MCDP composition with named subsystems. One declares outer
functionalities (\texttt{provides}) and outer resources (\texttt{requires}),
adds subsystems by name (\texttt{add}), wires them with connection
constraints (\texttt{constrain}), and emits a single \texttt{DesignProblem}
with \texttt{build}. The result solves like any other DP and can itself be
nested as a subsystem. Internally \texttt{build} closes one feedback loop
over the bundle of every subsystem's $R$ ports (the hidden
\texttt{\_\_modules\_\_} axis).

Signature: \texttt{System(name="system")}.

\begin{lstlisting}
from codesign import Module, Reals, System, solve

class Battery(Module):
    F = {"capacity": Reals(unit="J")}
    R = {"mass":     Reals(unit="kg")}
    def h(self, f):
        return {"mass": f["capacity"] / 1.8e6}

sys = System("drone")
endurance  = sys.provides("endurance", unit="s")
total_mass = sys.requires("total_mass", unit="kg")
battery    = sys.add("battery", Battery())
battery.capacity >= endurance * 5.0     # module F demand
total_mass       >= battery.mass        # outer R demand
res = solve(sys.build(), {"endurance": 300.0})
print(res.antichain)          # total_mass = 0.0008333 kg
\end{lstlisting}

Constructor parameters:
\begin{itemize}[nosep]
  \item \texttt{name} --- \texttt{str}, optional (default
        \texttt{"system"}). Names the emitted DP.
\end{itemize}

Public methods:
\begin{center}
\begin{tabular}{@{}p{0.32\linewidth}p{0.60\linewidth}@{}}
\toprule
\textbf{Method} & \textbf{Behaviour} \\
\midrule
\texttt{provides(name, *, unit="", poset=None)} &
  Declare an outer functionality. Returns a \texttt{Port} (kind
  \texttt{outer\_f}) usable in constraint expressions. \\
\texttt{requires(name, *, unit="", poset=None)} &
  Declare an outer resource. Returns a \texttt{Port} (kind
  \texttt{outer\_r}) usable as the LHS of a \texttt{>=} constraint. \\
\texttt{add(module\_name, dp)} &
  Add a subsystem under a unique name; returns a \texttt{ModuleHandle}
  (see \ref{ref:sugar:modulehandle}). Aliased as \texttt{sub}. \\
\texttt{constrain(target, demand)} &
  Register \texttt{target >= demand}. \texttt{target} is a
  \texttt{"module.port"} string, an outer-R name, or a \texttt{Port};
  \texttt{demand} is a \texttt{ctx}-dict callable or an \texttt{Expr}.
  Repeated calls on one target are joined (max). Aliased as \texttt{eq}. \\
\texttt{build()} &
  Validate and emit the composed \texttt{DesignProblem}. \\
\texttt{draw\_diagram(**kwargs)} &
  Render a GraphViz block diagram; delegates to
  \texttt{draw\_system} (see \ref{ref:diagram:draw_system}). \\
\bottomrule
\end{tabular}
\end{center}

The operator-overloaded \texttt{>=} form (as in the snippet) is sugar for
\texttt{constrain}: \texttt{Port.\_\_ge\_\_} registers the same internal
constraint, so both syntaxes yield identical solves.

Errors: \texttt{provides}/\texttt{requires} raise \texttt{ValueError} on a
duplicate name. \texttt{add} raises \texttt{ValueError} if the name is
reused, contains a dot, equals the reserved \texttt{\_\_modules\_\_}, or if
the DP's $F$/$R$ are not \texttt{Ports}. \texttt{constrain} raises
\texttt{TypeError} if \texttt{demand} is neither callable nor an
\texttt{Expr}, or if \texttt{target} is a non-constrainable port kind.
\texttt{build} (via \texttt{\_validate}) raises \texttt{ValueError} when
there are no subsystems and no outer resources, when no outer resource is
declared, when a constraint targets an unknown module or port, or when a
subsystem $F$ port or outer $R$ has no constraint. \texttt{draw\_diagram}
raises \texttt{ImportError} if \texttt{graphviz} is not installed.

\subsubsection{\texttt{module}}
\label{ref:module}

The \texttt{codesign.module} module supplies \texttt{Module}, a class-based
way to declare a design problem by subclassing rather than by wiring up
posets and an \texttt{h\_fn} by hand. It is the natural unit to hand to
\texttt{System.add} (see \ref{ref:system:system}).

\paragraph{\texttt{Module}}\label{ref:module:module}
Base class for declarative design problems. A subclass provides two
class-level \texttt{dict[str, Poset]} attributes, \texttt{F} and \texttt{R},
and an instance method \texttt{h(self, f)}; the constructor wraps these into
the underlying \texttt{DesignProblem} machinery, so an instance is itself a
\texttt{DesignProblem}. The subclass \texttt{h} may return a \texttt{dict}
(a singleton antichain), a list of dicts (multi-valued), or an
\texttt{Antichain} directly --- the return value is normalised
automatically. An optional class attribute \texttt{module\_name} overrides
the default name (the lower-cased class name).

Signature: \texttt{Module()} (subclasses may override \texttt{\_\_init\_\_}
to take parameters, calling \texttt{super().\_\_init\_\_()} last).

\begin{lstlisting}
from codesign import Module, Reals, solve

class Battery(Module):
    F = {"capacity": Reals(unit="J")}
    R = {"mass":     Reals(unit="kg")}
    def h(self, f):
        return {"mass": f["capacity"] / 1.8e6}

battery = Battery()                       # a DesignProblem
print(battery.h({"capacity": 3.6e6}))     # Antichain[(mass=2 kg)]
print(solve(battery, {"capacity": 1.8e6}).antichain)  # mass=1 kg
\end{lstlisting}

Class-level declarations (set by the subclass, not passed to
\texttt{\_\_init\_\_}):
\begin{itemize}[nosep]
  \item \texttt{F} --- \texttt{dict[str, Poset]}, functionality ports.
  \item \texttt{R} --- \texttt{dict[str, Poset]}, resource ports. Units are
        carried by the individual posets (e.g. \texttt{Reals(unit="kg")}).
  \item \texttt{module\_name} --- \texttt{str}, optional name override.
\end{itemize}

\subparagraph{\texttt{h(self, f)}}\label{ref:module:module:h}
The relation to override. Takes a functionality dict \texttt{f} keyed by
the \texttt{F} port names, returns the resource requirement as a dict, a
list of dicts, or an \texttt{Antichain}. The framework wraps non-antichain
returns via \texttt{Antichain.singleton} / \texttt{Antichain.from\_set}.

Errors: \texttt{\_\_init\_\_} raises \texttt{ValueError} if neither
\texttt{F} nor \texttt{R} is declared. The default \texttt{h} raises
\texttt{NotImplementedError} if a subclass fails to override it. The return
normaliser raises \texttt{TypeError} if \texttt{h} returns something other
than a dict, a list/tuple of dicts, or an \texttt{Antichain}.

\subsubsection{\texttt{sugar}}
\label{ref:sugar}

The \texttt{codesign.sugar} module implements the operator-overloaded
constraint DSL used by the \texttt{System} builder
(Section~\ref{sec:system-builder}). A \texttt{Port} is a typed handle onto
one of a system's ports; arithmetic on ports builds an \texttt{Expr} tree
lazily; \texttt{>=} on a constrainable port registers a constraint;
\texttt{ModuleHandle} (returned by \texttt{System.add}) yields ports by
attribute access; and \texttt{sqrt}/\texttt{exp}/\texttt{log} lift the
corresponding \texttt{math} functions to operate on \texttt{Expr} trees.

The following listing exercises all six symbols; individual cards refer
back to it.

\begin{lstlisting}
from codesign import Module, Reals, System, solve, sqrt, exp, log

class Battery(Module):
    F = {"capacity": Reals(unit="J")}
    R = {"mass":     Reals(unit="kg")}
    def h(self, f):
        return {"mass": f["capacity"] / 1.8e6}

sys = System("s")
endurance = sys.provides("endurance", unit="s")   # outer_f Port
total     = sys.requires("total_mass", unit="kg")  # outer_r Port
bat       = sys.add("battery", Battery())          # ModuleHandle
e = (bat.mass + 0.5) * sqrt(endurance)             # an Expr
print(e.pretty())   # ((battery.mass + 0.5) * sqrt(endurance))
bat.capacity >= endurance * 5.0    # >= registers a constraint
total        >= bat.mass
print(solve(sys.build(), {"endurance": 300.0}).antichain)
\end{lstlisting}

\paragraph{\texttt{Port}}\label{ref:sugar:port}
A named handle onto a port of a subsystem or an outer $F$/$R$ name. Its
\texttt{kind} is one of \texttt{"module\_f"}, \texttt{"module\_r"},
\texttt{"outer\_f"}, \texttt{"outer\_r"}, which governs where it may
appear: \texttt{module\_f} and \texttt{outer\_r} ports may be the LHS of a
\texttt{>=} constraint; \texttt{module\_r} and \texttt{outer\_f} ports may
appear on the demand (RHS) side. Ports are not built directly --- they are
returned by \texttt{System.provides}/\texttt{requires} and by
\texttt{ModuleHandle} attribute access. Arithmetic operators (inherited
from \texttt{Expr}) build expression trees; \texttt{>=} registers a
constraint with the parent system.

\begin{itemize}[nosep]
  \item \texttt{full\_name} --- property; \texttt{"module.port"} for
        module ports, the bare name otherwise.
  \item \texttt{kind} --- property; the port-kind string above.
  \item \texttt{pretty()} --- the \texttt{full\_name}, used in expression
        rendering.
\end{itemize}

Errors: \texttt{Port.\_\_ge\_\_} raises \texttt{TypeError} when the LHS is a
\texttt{module\_r} port (set by the module's \texttt{h}) or an
\texttt{outer\_f} port (a system input); both name the offending port.

\paragraph{\texttt{Expr}}\label{ref:sugar:expr}
Base class for constraint-expression AST nodes (\texttt{Add}, \texttt{Sub},
\texttt{Mul}, \texttt{Div}, \texttt{Pow}, \texttt{Neg}, \texttt{Func},
\texttt{Const}, and \texttt{Port}). Users rarely name \texttt{Expr}
directly; instances arise from arithmetic on ports and numbers. The
arithmetic operators (\texttt{+ - * / **} and unary \texttt{-}) return new
\texttt{Expr} nodes rather than evaluating numerically, so the whole demand
is captured as a tree that \texttt{System} later compiles to a callable.

\begin{itemize}[nosep]
  \item \texttt{pretty()} --- a readable, fully parenthesised string of the
        tree (e.g. \texttt{((battery.mass + 0.5) * sqrt(endurance))}); also
        the \texttt{\_\_repr\_\_}.
\end{itemize}

Errors: \texttt{Expr.\_\_ge\_\_} raises \texttt{TypeError} (only ports, not
compound expressions, may be a constraint LHS); \texttt{\_\_le\_\_} raises
\texttt{TypeError} (constraints are written \texttt{target >= demand}, not
reversed); \texttt{\_\_bool\_\_} raises \texttt{TypeError} to catch an
expression used in a conditional. Coercing a non-numeric, non-\texttt{Expr}
operand raises \texttt{TypeError}.

\paragraph{\texttt{ModuleHandle}}\label{ref:sugar:modulehandle}
The handle returned by \texttt{System.add}. Attribute access
(\texttt{handle.port\_name}) returns a \texttt{Port} of kind
\texttt{module\_f} or \texttt{module\_r} depending on which side of the
subsystem's DP the port lives on. Within one subsystem the $F$ and $R$
port names must be disjoint (enforced at construction).

\begin{itemize}[nosep]
  \item \texttt{name} --- property; the subsystem's registered name.
  \item \texttt{dp} --- property; the underlying \texttt{DesignProblem}.
\end{itemize}

Its \texttt{\_\_repr\_\_} lists the $F$ and $R$ ports, e.g.
\texttt{<ModuleHandle 'battery': F=['capacity'] R=['mass']>}.

Errors: the constructor raises \texttt{ValueError} if the subsystem's $F$
and $R$ port names overlap; attribute access raises \texttt{AttributeError}
(listing the available $F$/$R$ ports) for an unknown port name.

\paragraph{\texttt{sqrt}}\label{ref:sugar:sqrt}
Lift \texttt{math.sqrt} onto an expression: \texttt{sqrt(x) -> Expr}. The
argument may be a \texttt{Port}, an \texttt{Expr}, or a number (coerced to
a constant); the result is a \texttt{Func} node preserving the tree, as in
the listing above (\texttt{sqrt(endurance)}). Errors: raises
\texttt{TypeError} if the argument is neither a number nor an \texttt{Expr}.

\paragraph{\texttt{exp}}\label{ref:sugar:exp}
Lift \texttt{math.exp} onto an expression: \texttt{exp(x) -> Expr}. Same
argument rules as \texttt{sqrt} (see \ref{ref:sugar:sqrt}); e.g.
\texttt{exp(endurance).pretty()} is \texttt{"exp(endurance)"}. Errors:
raises \texttt{TypeError} on a non-number, non-\texttt{Expr} argument.

\paragraph{\texttt{log}}\label{ref:sugar:log}
Lift \texttt{math.log} (natural log) onto an expression:
\texttt{log(x) -> Expr}. Same argument rules as \texttt{sqrt} (see
\ref{ref:sugar:sqrt}); e.g. \texttt{log(endurance).pretty()} is
\texttt{"log(endurance)"}. Errors: raises \texttt{TypeError} on a
non-number, non-\texttt{Expr} argument.

\subsubsection{\texttt{diagram}}
\label{ref:diagram}

The \texttt{codesign.diagram} module renders a \texttt{System} as a
GraphViz block diagram --- one box per subsystem with its $F$ ports on the
left and $R$ ports on the right, outer $F$/$R$ as margin nodes, constraint
wiring as port-to-port edges, and any Kleene-iteration cycle coloured. See
Section~\ref{sec:block-diagrams} for the rendered figures.

\paragraph{\texttt{draw\_system}}\label{ref:diagram:draw_system}
Build a \texttt{graphviz.Digraph} for a system. It accepts either a live
\texttt{System} or the \texttt{DesignProblem} returned by
\texttt{System.build()} (which carries the modules and constraints on its
\texttt{\_codesign\_*} attributes). Constraints written in the
operator-overloaded form resolve to individual ports; lambda-based
constraints are drawn as a dashed edge from a \texttt{$\lambda$} marker so
none are silently dropped. Building the \texttt{Digraph} object needs only
the \texttt{graphviz} Python package; the \texttt{dot} binary is required
only to \emph{render} it to a file.

Signature:
\begin{lstlisting}
draw_system(system, *, name=None, rankdir="LR", show_ports=True,
            highlight_cycles=True, graph_attrs=None)
\end{lstlisting}

\begin{lstlisting}
from codesign import Module, Reals, System, draw_system

class Battery(Module):
    F = {"capacity": Reals(unit="J")}
    R = {"mass":     Reals(unit="kg")}
    def h(self, f):
        return {"mass": f["capacity"] / 1.8e6}

sys = System("drone")
endurance  = sys.provides("endurance", unit="s")
total_mass = sys.requires("total_mass", unit="kg")
battery    = sys.add("battery", Battery())
battery.capacity >= endurance * 5.0
total_mass       >= battery.mass

dot = draw_system(sys, rankdir="LR")   # a graphviz.Digraph
# dot.render("drone", format="svg")    # writes drone.svg (needs dot)
\end{lstlisting}

Parameters:
\begin{itemize}[nosep]
  \item \texttt{system} --- a live \texttt{System} or a built System DP.
  \item \texttt{name} --- \texttt{str} or \texttt{None}; diagram title,
        defaults to \texttt{system.name}.
  \item \texttt{rankdir} --- \texttt{"LR"} (default, horizontal) or
        \texttt{"TB"} (vertical, for tall systems).
  \item \texttt{show\_ports} --- \texttt{bool}; when \texttt{False}, only
        module names are shown and edges attach to the box body.
  \item \texttt{highlight\_cycles} --- \texttt{bool}; colour edges inside
        strongly-connected components of size $>1$.
  \item \texttt{graph\_attrs} --- \texttt{dict} of extra GraphViz graph
        attributes (e.g. \texttt{\{"ranksep": "0.6"\}}).
\end{itemize}

Errors: raises \texttt{ImportError} if the \texttt{graphviz} package is not
installed; raises \texttt{TypeError} if \texttt{system} is neither a
\texttt{System} nor a built System DP (i.e. lacks the
\texttt{\_codesign\_modules} attribute). The convenience method
\texttt{System.draw\_diagram} forwards to this function.

\subsubsection{\texttt{viz}}
\label{ref:viz}

The \texttt{codesign.viz} module collects ready-made diagnostics. It is
reached via the namespace import \texttt{from codesign import viz} (its
functions are deliberately not in the top-level \texttt{\_\_all\_\_}) and
called as \texttt{viz.plot\_antichain(...)}, and so on. The three
\texttt{plot\_*} functions require \texttt{matplotlib} (imported lazily,
so the rest of the package works without it) and return a
\texttt{matplotlib.axes.Axes}; \texttt{to\_dot} has no plotting dependency
and returns a string. Section~\ref{sec:visualisation} shows the figures in
context.

The snippet below builds a small drone system used by every card that
follows; the \texttt{plot\_*} calls run headless under
\texttt{MPLBACKEND=Agg}.

\begin{lstlisting}
from codesign import Module, Reals, System, solve, viz

# ... define Battery (F: capacity -> R: mass) and Actuator
# ... (F: lift_force -> R: power) Modules as in earlier cards
sys = System("drone")
endurance = sys.provides("endurance", unit="s")
power_out = sys.requires("power", unit="W")
mass_out  = sys.requires("total_mass", unit="kg")
b = sys.add("battery", Battery())
a = sys.add("actuator", Actuator())
b.capacity   >= a.power * endurance
a.lift_force >= 9.81 * b.mass
power_out >= a.power
mass_out  >= b.mass
dp = sys.build()
res = solve(dp, {"endurance": 300.0}, trace=True)
\end{lstlisting}

\paragraph{\texttt{viz.plot\_antichain}}\label{ref:viz:plot_antichain}
Scatter an antichain's Pareto front on two or three chosen $R$-port axes,
optionally shading the dominated region (2D only). Accepts a
\texttt{SolveResult}, an \texttt{UncertaintyResult} (uses its worst-case
antichain), or an \texttt{Antichain}.

Signature:
\begin{lstlisting}
plot_antichain(result, axes, *, ax=None, title=None,
               shade_dominated=True, point_size=60.0, label=None)
\end{lstlisting}

\begin{lstlisting}
ax = viz.plot_antichain(res, ["power", "total_mass"])
\end{lstlisting}

Parameters: \texttt{result} the object to plot; \texttt{axes} a sequence of
two or three $R$-port names; \texttt{ax} an existing axes (a new figure is
made if omitted); \texttt{title}; \texttt{shade\_dominated}; \texttt{point\_size};
\texttt{label} (a legend entry is drawn only when named).

Errors: raises \texttt{TypeError} if \texttt{result} is not one of the
three accepted types; \texttt{ValueError} if \texttt{axes} has other than
2 or 3 names, or if no point is plottable (empty, infeasible, or
non-numeric); \texttt{ImportError} if \texttt{matplotlib} is absent.

\paragraph{\texttt{viz.plot\_convergence}}\label{ref:viz:plot_convergence}
Semilog plot of the per-iteration Kleene delta from a solve trace. Accepts
a \texttt{SolveResult} whose \texttt{trace} was populated (call
\texttt{solve(..., trace=True)}) or the trace list directly; zero deltas
are clamped to \texttt{floor} for the log scale.

Signature:
\begin{lstlisting}
plot_convergence(result_or_trace, *, ax=None, floor=1e-18,
                 title="Kleene-iteration convergence", label=None)
\end{lstlisting}

\begin{lstlisting}
ax = viz.plot_convergence(res)   # res from solve(..., trace=True)
\end{lstlisting}

Parameters: \texttt{result\_or\_trace}; \texttt{ax}; \texttt{title};
\texttt{floor} (lower clamp for zero deltas); \texttt{label}.

Errors: raises \texttt{ValueError} if the trace has no numeric deltas;
\texttt{ImportError} if \texttt{matplotlib} is absent.

\paragraph{\texttt{viz.plot\_uncertainty}}\label{ref:viz:plot_uncertainty}
Histogram of Monte-Carlo samples on one $R$ port, with vertical lines for
each requested summary (mean, p95, CVaR95, worst case) and an optional
nominal marker. The \texttt{UncertaintyResult} must have been produced with
\texttt{"samples"} in the \texttt{uncertainty} list.

Signature:
\begin{lstlisting}
plot_uncertainty(result, port, *, ax=None, bins=30, title=None,
                 nominal=None, show_summaries=True, xlabel=None,
                 legend_loc="best")
\end{lstlisting}

\begin{lstlisting}
ures = solve(dp, {"endurance": 300.0},
             uncertainty=["mean", "p95", "samples"],
             n_samples=200, rng_seed=0)
ax = viz.plot_uncertainty(ures, "total_mass",
                          xlabel="total mass (kg)")
\end{lstlisting}

Parameters: \texttt{result} an \texttt{UncertaintyResult} with samples;
\texttt{port} the $R$ port to histogram; \texttt{bins}; \texttt{nominal};
\texttt{show\_summaries}; \texttt{xlabel} (defaults to the raw port name);
\texttt{legend\_loc}.

Errors: raises \texttt{ValueError} if the result carries no \texttt{samples},
or if the chosen port has no finite values; \texttt{ImportError} if
\texttt{matplotlib} is absent.

\paragraph{\texttt{viz.to\_dot}}\label{ref:viz:to_dot}
Emit a GraphViz \texttt{dot} string describing a DP's structure: a
System-built DP renders its subsystems and constraint edges, a composition
tree (\texttt{Series}/\texttt{Parallel}/\texttt{Loop}) renders as nested
boxes, and a leaf DP renders as a single box. Returning a string keeps the
module free of the \texttt{graphviz} dependency; save it and run
\texttt{dot -Tsvg}.

Signature: \texttt{to\_dot(dp, *, name="codesign") -> str}.

\begin{lstlisting}
print(viz.to_dot(dp, name="drone").splitlines()[0])
# digraph drone {
print(viz.to_dot(dp, name="my drone").splitlines()[0])
# digraph "my drone" {   (non-identifier name is quoted)
\end{lstlisting}

Parameters: \texttt{dp} the design problem to describe; \texttt{name} the
graph id --- used bare when it is a valid GraphViz identifier, otherwise
double-quoted with embedded quotes and backslashes escaped.

Constraint edges track the modules referenced on the right-hand side even
when a module port is wrapped in a function or negation --- the reference
walker descends through \texttt{.inner} (\texttt{sqrt(...)}, \texttt{exp},
a leading minus) as well as \texttt{.left}/\texttt{.right}/\texttt{.arg},
so a constraint like \texttt{bat.capacity >= sqrt(act.power)} still draws
the \texttt{act}$\rightarrow$\texttt{bat} edge.

Errors: none raised directly; the walker swallows introspection failures
and returns a best-effort graph.

\subsection{Uncertainty and online learning}

The set-based and stochastic \texttt{uncertainty} layer and the
compositional \texttt{online}-learning layer of optimistic evaluators. The
modelling discussion is in Sections~\ref{sec:uncertainty}
and~\ref{sec:online-learning}.

\subsubsection{\texttt{uncertainty}}
\label{ref:uncertainty}

Reference cards for the uncertainty layer of module
\texttt{codesign.uncertainty}. Both regimes are declared on a
\texttt{Module} through its \texttt{uncertain\_set} (set-based) and
\texttt{uncertain\_dist} (stochastic) attributes, and both are consumed by
the ordinary \texttt{solve} entry point; see Section~\ref{sec:uncertainty}
for the modelling discussion. The set-based classes
(\texttt{UncertaintySet}, \texttt{Box}, \texttt{Ellipsoid}, \texttt{Disk},
\texttt{Circle}) are grouped first, the stochastic ones
(\texttt{Independence}, \texttt{GaussianCopula}, \texttt{Stochastic}) after,
followed by the result dataclass and the solver entry point. Every
listing requires numpy; the stochastic cards additionally require scipy.

\paragraph{\texttt{UncertaintySet}}\label{ref:uncertainty:uncertaintyset}
Abstract base for deterministic, set-based parameter uncertainty. A concrete
subclass declares which parameters vary (\texttt{param\_names}) and where the
set's worst point lies (\texttt{worst\_case\_values}); the worst-case solve
then reassigns those named module attributes before re-solving.
Signature: \texttt{UncertaintySet()} (no constructor of its own).
\begin{lstlisting}
from codesign import Box
from codesign.uncertainty import UncertaintySet
s = Box(mass=(1.0, 2.0, "more_is_worse"))
print(isinstance(s, UncertaintySet))   # True
print(s.param_names())                 # ['mass']
\end{lstlisting}
Public methods (both abstract, overridden by subclasses):
\begin{itemize}[nosep]
  \item \texttt{param\_names() -> list[str]} --- names of the uncertain
    parameters carried by the set.
  \item \texttt{worst\_case\_values(h\_callable, base\_values, f\_inner) ->
    dict[str, float]} --- the parameter dict at the worst point of the set.
    \texttt{h\_callable(overrides) -> Antichain} is evaluated only for
    undeclared directions (boundary search); \texttt{base\_values} supplies
    nominal values for out-of-set parameters.
\end{itemize}
Errors: the base methods raise \texttt{NotImplementedError}; instantiate a
concrete subclass instead.

\paragraph{\texttt{Box}}\label{ref:uncertainty:box}
Axis-aligned interval product over one or more parameters --- the everyday
set-based model. When every parameter declares a direction of badness the
worst case is a single corner, read off in constant time.
Signature: \texttt{Box(**params)} where each value is \texttt{(lo, hi)} or
\texttt{(lo, hi, direction)}.
\begin{lstlisting}
from codesign import Box
b = Box(se=(1.7e6, 2.3e6, "more_is_better"),
        eff=(0.83, 0.97, "more_is_better"))
print(b.param_names())                    # ['se', 'eff']
print(b.worst_case_values(None, {}, {}))
# {'se': 1700000.0, 'eff': 0.83}
\end{lstlisting}
Parameters:
\begin{itemize}[nosep]
  \item \texttt{**params} --- one keyword per uncertain parameter. The value
    is \texttt{(lo, hi)} (direction undeclared) or \texttt{(lo, hi, direction)}
    with \texttt{direction} one of \texttt{"more\_is\_better"},
    \texttt{"more\_is\_worse"}, \texttt{"less\_is\_better"},
    \texttt{"less\_is\_worse"} (in parameter units).
\end{itemize}
With $n$ undeclared parameters all $2^n$ endpoint corners are probed
(cheap for modest $n$). Errors: raises \texttt{ValueError} if a value is
neither a 2- nor 3-tuple, if an unknown direction token is passed, or if
\texttt{lo > hi} (empty interval).

\paragraph{\texttt{Ellipsoid}}\label{ref:uncertainty:ellipsoid}
An $n$-D ellipsoid $(p-c)^\top \Sigma^{-1}(p-c) \leq 1$ in parameter space;
the covariance $\Sigma$ captures scale and correlation, so it is a more
honest set than a box when parameters vary together.
\begin{lstlisting}
from codesign import Ellipsoid
e = Ellipsoid(
    center={"se": 2.0e6, "eff": 0.9},
    cov=[[1.0e10, 0.0], [0.0, 2.5e-3]],
    params=["se", "eff"],
    directions={"se": "more_is_better",
                "eff": "more_is_better"},
)
print(e.worst_case_values(None, {}, {}))
# {'se': 1929289.32..., 'eff': 0.86464...}
\end{lstlisting}
Parameters:
\begin{itemize}[nosep]
  \item \texttt{center} --- dict of centre coordinates by parameter name.
  \item \texttt{cov} --- symmetric positive-definite shape matrix $\Sigma$,
    rows/columns ordered as \texttt{params}.
  \item \texttt{params} --- parameter names in \texttt{cov}'s order.
  \item \texttt{directions=None} --- optional per-parameter badness tokens;
    when all are given the worst case is the analytic boundary point
    $c + L\,u^\star$, else the boundary is sampled.
  \item \texttt{boundary\_samples=8} --- samples per dimension used when
    directions are not fully declared.
\end{itemize}
Errors: raises \texttt{ValueError} if \texttt{cov}'s shape does not match
\texttt{params} or if \texttt{cov} is not positive definite;
\texttt{ImportError} if numpy is absent.

\paragraph{\texttt{Disk}}\label{ref:uncertainty:disk}
Two-dimensional convenience: a filled circle, i.e.\ the isotropic ellipsoid
$\Sigma = r^2 I$. A factory function returning an \texttt{Ellipsoid}, not a
class.
Signature: \texttt{Disk(center, radius, params=None, directions=None)}.
\begin{lstlisting}
from codesign import Disk, Ellipsoid
d = Disk(center={"x": 0.0, "y": 0.0}, radius=1.0,
         directions={"x": "more_is_worse",
                     "y": "more_is_worse"})
print(isinstance(d, Ellipsoid))               # True
print(d.worst_case_values(None, {}, {}))
# {'x': 0.70710..., 'y': 0.70710...}
\end{lstlisting}
Parameters: \texttt{center} (2-key dict), \texttt{radius} (in parameter
units), \texttt{params} (inferred from \texttt{center} if \texttt{None}),
\texttt{directions} (as for \texttt{Ellipsoid}). Errors: raises
\texttt{ValueError} unless exactly two parameters are supplied.

\paragraph{\texttt{Circle}}\label{ref:uncertainty:circle}
Two-dimensional convenience for the circle boundary only. For monotone
modules with declared directions the worst case lies on the boundary, so
\texttt{Circle} and \texttt{Disk} return identical worst-case answers; the
implementation delegates to \texttt{Disk}.
Signature: \texttt{Circle(center, radius, params=None, directions=None)}.
\begin{lstlisting}
from codesign import Circle, Ellipsoid
c = Circle(center={"x": 0.0, "y": 0.0}, radius=2.0,
           directions={"x": "more_is_worse",
                       "y": "more_is_worse"})
print(isinstance(c, Ellipsoid))               # True
print(c.worst_case_values(None, {}, {}))
# {'x': 1.41421..., 'y': 1.41421...}
\end{lstlisting}
Parameters: as \texttt{Disk}. Errors: as \texttt{Disk} (raises
\texttt{ValueError} unless exactly two parameters are supplied).

\paragraph{\texttt{Independence}}\label{ref:uncertainty:independence}
The independence copula: each component drawn independently in $U(0,1)$. The
default dependence structure of \texttt{Stochastic}.
Signature: \texttt{Independence()}; method
\texttt{sample\_uniform(n, d, rng) -> ndarray} of shape \texttt{(n, d)}.
\begin{lstlisting}
from codesign.uncertainty import Independence
import numpy as np
rng = np.random.default_rng(0)
u = Independence().sample_uniform(3, 2, rng)
print(u.shape)                                # (3, 2)
\end{lstlisting}
Parameters: none. \texttt{sample\_uniform} takes the sample count \texttt{n},
the dimension \texttt{d}, and a numpy \texttt{Generator} \texttt{rng}.
Errors: none raised directly.

\paragraph{\texttt{GaussianCopula}}\label{ref:uncertainty:gaussiancopula}
Gaussian copula with a given correlation matrix: draw $z \sim N(0, R)$ and
apply $\Phi$ to each coordinate to obtain correlated uniforms.
Signature: \texttt{GaussianCopula(correlation)}.
\begin{lstlisting}
from codesign import GaussianCopula
import numpy as np
gc = GaussianCopula(correlation=[[1.0, 0.4],
                                 [0.4, 1.0]])
rng = np.random.default_rng(0)
u = gc.sample_uniform(4, 2, rng)
print(u.shape, repr(gc))
# (4, 2) GaussianCopula(d=2)
\end{lstlisting}
Parameters: \texttt{correlation} --- a $(d,d)$ symmetric positive-definite
matrix with unit diagonal. Errors: the constructor raises \texttt{ValueError}
if \texttt{correlation} is not square or not positive definite;
\texttt{sample\_uniform} raises \texttt{ValueError} if \texttt{d} differs from
the copula's dimension. Sampling additionally requires scipy.

\paragraph{\texttt{Stochastic}}\label{ref:uncertainty:stochastic}
Joint distribution over uncertain parameters, built from named
\texttt{scipy.stats} frozen marginals plus an optional copula. This is the
\texttt{uncertain\_dist} declaration consumed by the Monte Carlo path of
\texttt{solve}; it has a full treatment in Section~\ref{sec:uncertainty}.
Signature: \texttt{Stochastic(marginals=None, copula=None, **marginals\_kw)}.
\begin{lstlisting}
from codesign import Stochastic
from scipy import stats
import numpy as np
st = Stochastic(marginals={
    "m": stats.norm(loc=1.0, scale=0.1)})
print(st.param_names())                       # ['m']
print(st.sample(2, np.random.default_rng(0)))
# [{'m': 1.0350...}, {'m': 0.9386...}]
\end{lstlisting}
Parameters:
\begin{itemize}[nosep]
  \item \texttt{marginals} --- dict of parameter name to a frozen scipy
    distribution (e.g.\ \texttt{stats.norm(loc, scale)}). Marginals may also
    be given as keyword arguments (\texttt{**marginals\_kw}).
  \item \texttt{copula=None} --- a \texttt{Copula}; defaults to
    \texttt{Independence}.
\end{itemize}
Methods: \texttt{param\_names() -> list[str]}; \texttt{sample(n, rng) ->
list[dict]} of $n$ parameter dicts (copula uniforms mapped through each
marginal's \texttt{ppf}). Errors: the constructor raises \texttt{ValueError}
if no marginal is supplied.

\paragraph{\texttt{UncertaintyResult}}\label{ref:uncertainty:uncertaintyresult}
Dataclass returned by an uncertainty-aware solve. Which fields are populated
depends on the summaries requested via \texttt{uncertainty=[...]}; the rest
stay \texttt{None}.
\begin{lstlisting}
from codesign.uncertainty import UncertaintyResult
r = UncertaintyResult(mean={"mass": 0.55},
                      feasibility_rate=1.0,
                      n_samples_used=1000)
print(r.mean, r.feasibility_rate)
# {'mass': 0.55} 1.0
\end{lstlisting}
Fields:
\begin{center}
\begin{tabular}{@{}lll@{}}
\toprule
Field & Type & Meaning \\
\midrule
\texttt{worst\_case} & \texttt{SolveResult} or \texttt{None} & result at the
  set's worst point \\
\texttt{mean} & \texttt{dict} or \texttt{None} & per-R-port MC mean \\
\texttt{p95} & \texttt{dict} or \texttt{None} & per-R-port 95th percentile \\
\texttt{cvar95} & \texttt{dict} or \texttt{None} & per-R-port CVaR at 95\% \\
\texttt{samples} & \texttt{list} or \texttt{None} & raw MC antichains \\
\texttt{feasibility\_rate} & \texttt{float} or \texttt{None} & feasible-sample
  fraction \\
\texttt{n\_samples\_used} & \texttt{int} & MC samples drawn (0 if none) \\
\bottomrule
\end{tabular}
\end{center}
Statistics are computed over feasible samples only. Errors: none raised
directly (a plain dataclass).

\paragraph{\texttt{solve\_with\_uncertainty}}%
\label{ref:uncertainty:solve_with_uncertainty}
Entry point for an uncertainty-aware solve, returning an
\texttt{UncertaintyResult}. Ordinarily reached through
\texttt{solve(dp, f, uncertainty=[...])}, which dispatches here when the
\texttt{uncertainty} argument is present. Signature:
\begin{lstlisting}
solve_with_uncertainty(dp, functionality, uncertainty,
    n_samples=1000, rng_seed=None, max_iter=200,
    verbose=0)
\end{lstlisting}
\begin{lstlisting}
from codesign import Module, Reals, solve, Box
class Battery(Module):
    F = {"capacity": Reals(unit="J")}
    R = {"mass":     Reals(unit="kg")}
    def __init__(self, se=2.0e6, eff=0.9):
        self.se, self.eff = se, eff
        super().__init__()
    def h(self, f):
        return {"mass": f["capacity"] / (self.se * self.eff)}
bat = Battery()
bat.uncertain_set = Box(se=(1.7e6, 2.3e6, "more_is_better"),
                        eff=(0.83, 0.97, "more_is_better"))
res = solve(bat, {"capacity": 1.0e6}, uncertainty=["worst_case"])
print(round(list(res.worst_case.antichain.points)[0]["mass"], 5))
# 0.70872
\end{lstlisting}
Parameters:
\begin{itemize}[nosep]
  \item \texttt{dp}, \texttt{functionality} --- the design problem and its
    outer F vector, as for \texttt{solve}.
  \item \texttt{uncertainty} --- list drawn from \texttt{"worst\_case"},
    \texttt{"mean"}, \texttt{"p95"}, \texttt{"cvar95"}, \texttt{"samples"}.
  \item \texttt{n\_samples=1000} --- Monte Carlo draws for the stochastic
    summaries; \texttt{rng\_seed=None} seeds the numpy generator.
  \item \texttt{max\_iter=200}, \texttt{verbose=0} --- forwarded to each
    inner solve.
\end{itemize}
Errors: raises \texttt{ValueError} for an unknown summary label, when no
module on \texttt{dp} carries an \texttt{uncertain\_set}/\texttt{uncertain\_dist},
or when stochastic summaries are requested but no \texttt{uncertain\_dist} is
present; \texttt{ImportError} if numpy is absent.

\subsubsection{\texttt{online}}
\label{ref:online}

Reference cards for the compositional online learner of module
\texttt{codesign.online} --- a port of \citet{alharbi2026online}
(\texttt{arXiv:2604.22624}). The evaluators maintain history-dependent
bounds on each candidate's inner-solve output so that provably suboptimal
candidates are eliminated without ever being solved;
\texttt{solve\_online} is the budgeted driver. See
Section~\ref{sec:online-learning} for the modelling discussion and
Section~\ref{sec:online-tier1} for the warm-start and picker options. Every
listing requires numpy; \texttt{LinearParametricEvaluator} additionally
requires scipy.

\paragraph{\texttt{OptimisticEvaluator}}\label{ref:online:optimisticevaluator}
Abstract base holding the observation history and the bound contract. The
public \texttt{bound(candidate)} returns a pair of dicts
\texttt{(lower, upper)} mapping each numeric R component to its current
lower and upper bound at the queried feature point; the base implementation
is the trivial fallback $(0, +\infty)$, which subclasses tighten.
Signature: \texttt{OptimisticEvaluator(features, r\_components)}.
\begin{lstlisting}
from codesign.online import OptimisticEvaluator
ev = OptimisticEvaluator(features=["x"],
                         r_components=["cost"])
print(ev.bound({"x": 2.0}))
# ({'cost': 0.0}, {'cost': inf})
\end{lstlisting}
Parameters: \texttt{features} --- candidate keys used as the feature vector;
\texttt{r\_components} --- resource-port names to bound.
\begin{center}
\begin{tabular}{@{}ll@{}}
\toprule
Method & Behaviour \\
\midrule
\texttt{reset()} & forget every observation \\
\texttt{observe(candidate\_id, candidate, antichain)} & record an
  inner-solve result \\
\texttt{bound(candidate)} & return \texttt{(lower, upper)} bound dicts \\
\bottomrule
\end{tabular}
\end{center}
Errors: \texttt{observe}/\texttt{bound} raise \texttt{KeyError} if a
candidate is missing a declared feature.

\paragraph{\texttt{MonotonicityEvaluator}}\label{ref:online:monotonicityevaluator}
Bounds from monotonicity alone: assuming the output is component-wise
monotone in the features, an observation at $f_0$ with antichain-min $R_0$
lower-bounds every candidate with features $\geq f_0$ by $R_0$ and
upper-bounds every candidate with features $\leq f_0$ by $R_0$. The
aggressive choice when the assumption genuinely holds.
Signature: \texttt{MonotonicityEvaluator(features, r\_components)}.
\begin{lstlisting}
from codesign import (MonotonicityEvaluator, Antichain,
                      Ports, Reals)
R = Ports({"cost": Reals()})
ev = MonotonicityEvaluator(["x"], ["cost"])
ev.observe(0, {"x": 1.0}, Antichain.singleton(R, {"cost": 5.0}))
print(ev.bound({"x": 2.0}))   # ({'cost': 5.0}, {'cost': inf})
print(ev.bound({"x": 0.5}))   # ({'cost': 0.0}, {'cost': 5.0})
\end{lstlisting}
Parameters: as \texttt{OptimisticEvaluator}. Behaviour is described, not
its implementation: the bound is exact to float comparisons, computed by a
vectorised numpy rescan of the whole history when numpy is present, with an
equivalent pure-Python fallback for numpy-absent installs. Errors: raises
\texttt{KeyError} for a missing feature.

\paragraph{\texttt{LipschitzEvaluator}}\label{ref:online:lipschitzevaluator}
Bounds from a Lipschitz assumption: with
$|h(c_1)-h(c_2)| \leq L\,\|features(c_1)-features(c_2)\|$, each observation
tightens the bound everywhere by a cone of slope $L$.
Signature: \texttt{LipschitzEvaluator(features, r\_components, L)}.
\begin{lstlisting}
from codesign import (LipschitzEvaluator, Antichain,
                      Ports, Reals)
R = Ports({"cost": Reals()})
ev = LipschitzEvaluator(["x"], ["cost"], L=2.0)
ev.observe(0, {"x": 1.0}, Antichain.singleton(R, {"cost": 5.0}))
print(ev.bound({"x": 2.0}))   # ({'cost': 3.0}, {'cost': 7.0})
\end{lstlisting}
Parameters: \texttt{L} --- a single positive float (same constant for every
R component) or a dict mapping R component name to its own constant.
Behaviour, not implementation: distances use a vectorised numpy path when
numpy is present and an equivalent pure-Python fallback otherwise. Errors:
raises \texttt{KeyError} for a missing feature, or if \texttt{L} is a dict
lacking an R component.

\paragraph{\texttt{LinearParametricEvaluator}}%
\label{ref:online:linearparametricevaluator}
The \emph{certified} optimistic bound of \citet{alharbi2026online} (Section~V-C3,
eqs.\ 26--28 and Lemma~V.5). Each resource coordinate is assumed an exact
affine function of the features, $\text{req}_k(c)=\phi(c)^\top\theta_k^*$
with $\phi(c)=[1,features(c)]$; the evaluator maintains the confidence
polytope $\Theta(H)$ of parameters consistent with every observation and
returns, per coordinate, $\min_{\theta\in\Theta(H)}\phi(i)^\top\theta$ --- one
LP solved with \texttt{scipy.optimize.linprog}. Because $\theta_k^*$ is
feasible, the optimum is $\leq \text{req}_k(i)$: a \emph{guaranteed} lower
bound (the optimism guarantee), so a Pareto-optimal candidate is never
wrongly eliminated. No upper bound is certified, so the returned upper bound
stays at $+\infty$. Signature:
\begin{lstlisting}
LinearParametricEvaluator(features, r_components,
    confidence=None, min_obs=3, *, prior_box=None,
    noise_bound=0.0, solver="highs")
\end{lstlisting}
\begin{lstlisting}
from codesign import (LinearParametricEvaluator, Antichain,
                      Ports, Reals)
R = Ports({"cost": Reals()})
ev = LinearParametricEvaluator(["x"], ["cost"],
                               prior_box=(-100.0, 100.0))
for x in (1.0, 2.0, 3.0):            # cost = 2x + 1, affine
    ev.observe(0, {"x": x},
               Antichain.singleton(R, {"cost": 2*x + 1}))
lo, hi = ev.bound({"x": 4.0})        # certified lower bound
print(round(lo["cost"], 3), hi["cost"])   # 9.0 inf
\end{lstlisting}
Parameters:
\begin{itemize}[nosep]
  \item \texttt{prior\_box=None} --- the prior box $\Theta_0$ keeping the LP
    bounded while the fit is under-determined. \texttt{None} leaves every
    parameter unbounded (always safe, but trivial until $\theta$ is pinned
    down); a \texttt{(lo, hi)} float pair applies to every parameter; a
    dict gives each R component its own \texttt{(lo, hi)}. A finite box must
    contain the true $\theta^*$ or the bound may become invalid.
  \item \texttt{noise\_bound=0.0} --- half-width of the observation band
    $|\phi(j)^\top\theta - r| \leq$ \texttt{noise\_bound}; \texttt{0.0} is
    the paper's exact/noiseless model, positive values a documented
    bounded-noise extension.
  \item \texttt{min\_obs=3} --- below this many observations the trivial
    $(0,+\infty)$ bound is returned.
  \item \texttt{solver="highs"} --- \texttt{method} forwarded to
    \texttt{linprog}.
  \item \texttt{confidence=None} --- \textbf{deprecated and ignored},
    retained only for backward compatibility with the former OLS
    $\pm$ confidence-band version; passing it emits a
    \texttt{DeprecationWarning}. The certified bound has no confidence-band
    parameter.
\end{itemize}
Errors: raises \texttt{ValueError} if \texttt{noise\_bound} is negative;
\texttt{TypeError} if a \texttt{prior\_box} entry is not a \texttt{(lo, hi)}
pair; \texttt{ImportError} if numpy or scipy is absent. An unbounded,
infeasible, or failed LP degrades safely to the trivial lower bound.

\paragraph{\texttt{GaussianProcessEvaluator}}%
\label{ref:online:gaussianprocessevaluator}
Bounds from a zero-mean Gaussian-process surrogate with an RBF kernel,
fitted in pure numpy; the bound at a query is $\text{mean} \pm
\texttt{confidence}\cdot\sigma$. The right tool when the response surface has
feature interactions or local nonlinearity, where the certified linear
evaluator would degrade to the no-information bound. The bounds are
informative but \emph{not} certified. Signature:
\begin{lstlisting}
GaussianProcessEvaluator(features, r_components,
    length_scale=0.3, sigma_f=1.0, noise=1e-3,
    confidence=2.0, min_obs=3)
\end{lstlisting}
\begin{lstlisting}
from codesign import (GaussianProcessEvaluator, Antichain,
                      Ports, Reals)
R = Ports({"cost": Reals()})
ev = GaussianProcessEvaluator(["x"], ["cost"],
                              length_scale=0.5)
for x in (0.0, 0.5, 1.0):
    ev.observe(0, {"x": x},
               Antichain.singleton(R, {"cost": x}))
lo, hi = ev.bound({"x": 0.25})
print(round(lo["cost"], 3), round(hi["cost"], 3))
# 0.066 0.289
\end{lstlisting}
Parameters:
\begin{itemize}[nosep]
  \item \texttt{length\_scale=0.3} --- RBF length scale (0.3 suits features
    in roughly $[0,1]$).
  \item \texttt{sigma\_f=1.0} --- kernel amplitude, rescaled per output by
    the empirical standard deviation.
  \item \texttt{noise=1e-3} --- observation-noise jitter stabilising the
    Cholesky factor.
  \item \texttt{confidence=2.0} --- multiplier on the predictive standard
    deviation ($\approx 95\%$ band under Gaussian residuals).
  \item \texttt{min\_obs=3} --- below this the fallback $(0,+\infty)$ bound
    is returned.
\end{itemize}
\texttt{observe} additionally invalidates the fit cache (a lazy refit runs
on the next \texttt{bound}). Errors: raises \texttt{KeyError} for a missing
feature; \texttt{ImportError} if numpy is absent.

\paragraph{\texttt{OnlineResult}}\label{ref:online:onlineresult}
Dataclass returned by \texttt{solve\_online}, carrying the surviving
antichain plus a full audit of the elimination cascade.
\begin{lstlisting}
from codesign import OnlineResult, Antichain, Ports, Reals
R = Ports({"cost": Reals()})
r = OnlineResult(
    antichain=Antichain.singleton(R, {"cost": 5.0}),
    n_evaluated=3, n_eliminated=7, n_candidates=10)
print(r.n_evaluated, r.n_eliminated, r.n_candidates)
# 3 7 10
\end{lstlisting}
Fields:
\begin{center}
\begin{tabular}{@{}lll@{}}
\toprule
Field & Type & Meaning \\
\midrule
\texttt{antichain} & \texttt{Antichain} & Min over evaluated survivors \\
\texttt{n\_evaluated} & \texttt{int} & inner solves run \\
\texttt{n\_eliminated} & \texttt{int} & candidates pruned by bounds \\
\texttt{n\_candidates} & \texttt{int} & candidates at start \\
\texttt{history} & \texttt{list[dict]} & per-iteration log; each entry has keys \texttt{pick}, \texttt{antichain}, \texttt{remaining}, \texttt{evaluated}, \texttt{eliminated}, \texttt{phase} \\
\texttt{evaluated\_ids} & \texttt{list[int]} & indices actually evaluated \\
\texttt{eliminated\_ids} & \texttt{list[int]} & indices pruned \\
\texttt{incumbent\_ids} & \texttt{list[int]} & indices in the final antichain \\
\bottomrule
\end{tabular}
\end{center}
Errors: none raised directly (a plain dataclass).

\paragraph{\texttt{solve\_online}}\label{ref:online:solve_online}
The budgeted online driver: \texttt{solve}'s elimination-based cousin. For
each candidate, \texttt{candidate\_fn} builds a fresh DP; the evaluator's
bounds prune provably suboptimal candidates before their inner solve is run.
Signature:
\begin{lstlisting}
solve_online(candidate_fn, functionality, *,
    candidates, evaluator, budget=None, max_iter=200,
    verbose=0, warm_start=None, picker="lcb")
\end{lstlisting}
\begin{lstlisting}
from codesign import (AlgebraicDP, Ports, Reals,
                      LipschitzEvaluator, solve_online)
F, R = Ports({"load": Reals()}), Ports({"cost": Reals()})
def candidate_fn(c):
    return AlgebraicDP(F, R, {"cost":
        lambda f, u=c["price"]: f["load"] * u})
cands = [{"price": p} for p in (2.0, 3.0, 5.0, 8.0)]
ev = LipschitzEvaluator(["price"], ["cost"], L=10.0)
res = solve_online(candidate_fn, {"load": 1.0},
                   candidates=cands, evaluator=ev,
                   warm_start=2, picker="ucb")
print(res.n_evaluated, res.n_eliminated, res.antichain)
# 4 0 Antichain[(cost=2)]
\end{lstlisting}
Parameters:
\begin{itemize}[nosep]
  \item \texttt{candidate\_fn} --- \texttt{candidate\_fn(candidate) -> DP}.
  \item \texttt{functionality} --- the outer F vector, passed to every inner
    solve unchanged.
  \item \texttt{candidates} --- list of feature dicts; each must contain
    every name in \texttt{evaluator.features}.
  \item \texttt{evaluator} --- the \texttt{OptimisticEvaluator} instance
    (reset and repopulated during the solve).
  \item \texttt{budget=None} --- maximum inner solves; \texttt{None} is
    unlimited (elimination still prunes).
  \item \texttt{max\_iter=200}, \texttt{verbose=0} --- forwarded to each
    inner solve; \texttt{verbose} 0/1/2 for silent/summary/trace.
  \item \texttt{warm\_start=None} --- seed evaluations before the picker
    takes over: an \texttt{int} $n$ triggers greedy farthest-point sampling
    of $n$ candidates; a list of indices seeds exactly those.
  \item \texttt{picker="lcb"} --- selection strategy: \texttt{"lcb"}
    (exploit), \texttt{"ucb"} (add exploration bonus, tune via
    \texttt{("ucb", \{"kappa": 1.0\})}), \texttt{"random"}, or a callable
    \texttt{(lo, hi, r\_components, **kwargs) -> score}.
\end{itemize}
Errors: raises \texttt{ValueError} for an unknown picker name and
\texttt{TypeError} for a picker of the wrong type; \texttt{KeyError} if a
candidate is missing a declared feature.

\subsection{Temporal layers}

The open-loop temporal planners: architecture switching
(\texttt{temporal}), scalar dynamic programming (\texttt{dynamic}), and
antichain-valued sequential co-design (\texttt{sequential}). See
Section~\ref{sec:temporal}.

\subsubsection{\texttt{temporal}}\label{ref:temporal}

The \texttt{temporal} module implements Case~1 of the temporal co-design
layers (section~\ref{sec:temporal}): architecture \emph{switching} across a
changing environment. A sequence of epochs each poses its own functionality
and admits a set of candidate architectures; \texttt{solve\_schedule}
chooses one architecture per epoch to minimise the total of the per-epoch
co-design costs plus the switching costs, by an exact Viterbi
(shortest-path) pass over the epoch-by-architecture lattice. The per-epoch
cost is an ordinary \texttt{solve()}, so this layer sits on top of static
co-design without modifying it. The \texttt{Architecture} defined here is
the shared decision object reused by every temporal layer.

\paragraph{\texttt{Architecture}}\label{ref:temporal:architecture}
A named candidate configuration that can be active in an epoch or stage;
the common decision object of the whole temporal family. It is a
dataclass wrapping a design problem.
\begin{lstlisting}
Architecture(name: str, dp: Any, tags: Mapping = {})
\end{lstlisting}
\begin{lstlisting}
from codesign import AlgebraicDP, Architecture, Ports, Reals, System
s = System("eco")
d = s.provides("demand")
s.requires("cost")
s.add("m", AlgebraicDP(F=Ports({"demand": Reals()}),
      R=Ports({"cost": Reals()}),
      equations={"cost": lambda f: 10.0})).demand >= d
s.constrain("cost", lambda x: x["m.cost"])
eco = Architecture("eco", s.build(), tags={"mode": "eco"})
print(eco.name, eco.tags)          # -> eco {'mode': 'eco'}
\end{lstlisting}
\begin{itemize}[nosep]
  \item \texttt{name} --- \texttt{str} --- identifier used in results, plots
    and switching bookkeeping.
  \item \texttt{dp} --- \texttt{Any} --- a design problem accepted by
    \texttt{solve()} (typically a built \texttt{System}); may be shared
    across epochs or epoch-specific.
  \item \texttt{tags} --- \texttt{Mapping} --- free-form metadata carried
    through to results; not interpreted by the solver.
\end{itemize}
Errors: none raised directly.

\paragraph{\texttt{Epoch}}\label{ref:temporal:epoch}
One environment regime over which a single architecture is active. The
changing environment enters through the per-epoch \texttt{functionality}.
\begin{lstlisting}
Epoch(name: str, functionality: Mapping,
      candidates: Optional[Sequence[Architecture]] = None,
      duration: float = 1.0)
\end{lstlisting}
\begin{lstlisting}
from codesign import Architecture, Epoch
arch = Architecture("eco", None, tags={"mode": "eco"})
ep = Epoch("glucose", {"mu": 0.8}, candidates=[arch], duration=2.0)
print(ep.name, ep.functionality, ep.duration)
# -> glucose {'mu': 0.8} 2.0
\end{lstlisting}
\begin{itemize}[nosep]
  \item \texttt{name} --- \texttt{str} --- epoch identifier.
  \item \texttt{functionality} --- \texttt{Mapping} --- the outer
    functionality demanded in this epoch, passed straight to
    \texttt{solve()}.
  \item \texttt{candidates} --- \texttt{Sequence[Architecture]} or
    \texttt{None} --- architectures admissible here; when \texttt{None} the
    schedule-level default set is used, so some architectures can be made
    unavailable in some regimes.
  \item \texttt{duration} --- \texttt{float} --- weight on the epoch's
    running cost, expressing unequal epoch lengths without rescaling the
    cost function.
\end{itemize}
Errors: none raised directly (a missing default set is reported by
\texttt{solve\_schedule}; see below).

\paragraph{\texttt{EpochResult}}\label{ref:temporal:epochresult}
The per-epoch outcome recorded inside a solved schedule. A dataclass; its
\texttt{total\_cost} property sums the running and switching contributions.
\begin{lstlisting}
from codesign import EpochResult
er = EpochResult(epoch="glucose", architecture="eco",
                 running_cost=8.0, switch_cost=0.5, feasible=True,
                 point={"cost": 8.0})
print(er.total_cost)               # -> 8.5
\end{lstlisting}
\begin{itemize}[nosep]
  \item \texttt{epoch} --- \texttt{str} --- the epoch's name.
  \item \texttt{architecture} --- \texttt{str} --- name of the chosen
    architecture.
  \item \texttt{running\_cost} --- \texttt{float} --- duration-weighted
    scalar cost of the chosen design point (\texttt{inf} if infeasible).
  \item \texttt{switch\_cost} --- \texttt{float} --- switching cost charged
    at this epoch's boundary (zero when the architecture is unchanged).
  \item \texttt{feasible} --- \texttt{bool} --- whether this epoch was
    satisfiable by its chosen architecture.
  \item \texttt{point} --- \texttt{Mapping} or \texttt{None} --- the chosen
    resource point.
  \item \texttt{tags} --- \texttt{Mapping} --- the chosen architecture's
    tags.
  \item \texttt{total\_cost} (property) --- \texttt{running\_cost} $+$
    \texttt{switch\_cost}.
\end{itemize}
Errors: none raised directly.

\paragraph{\texttt{ScheduleResult}}\label{ref:temporal:scheduleresult}
The outcome of \texttt{solve\_schedule}: the chosen architecture per epoch
with decomposed costs. A dataclass with a \texttt{schedule} property and a
compact \texttt{\_\_repr\_\_}.
\begin{lstlisting}
from codesign import EpochResult, ScheduleResult
e0 = EpochResult("glucose", "eco", 8.0, 0.0, True, {"cost": 8.0})
e1 = EpochResult("acetate", "bal", 6.0, 0.5, True, {"cost": 6.0})
sr = ScheduleResult([e0, e1], total_cost=14.5, feasible=True,
                    n_switches=1)
print(sr.schedule)
print(sr)
# -> ['eco', 'bal']
# -> ScheduleResult(feasible, cost=14.5, switches=1, [eco -> bal])
\end{lstlisting}
\begin{itemize}[nosep]
  \item \texttt{epochs} --- \texttt{List[EpochResult]} --- one entry per
    epoch, in order.
  \item \texttt{total\_cost} --- \texttt{float} --- running plus switching
    cost summed over epochs.
  \item \texttt{feasible} --- \texttt{bool} --- \texttt{True} iff every
    epoch was satisfiable by its chosen architecture.
  \item \texttt{n\_switches} --- \texttt{int} --- number of epoch boundaries
    at which the active architecture changed.
  \item \texttt{schedule} (property) --- \texttt{List[str]} --- the chosen
    architecture name per epoch, in order.
\end{itemize}
Errors: none raised directly.

\paragraph{\texttt{solve\_schedule}}\label{ref:temporal:solve_schedule}
Chooses an architecture per epoch minimising running plus switching cost,
by an exact Viterbi pass over the epoch/architecture lattice. With
\texttt{switch\_cost == 0} and \texttt{hysteresis == 0} the result reduces
to the epoch-local greedy choice.
\begin{lstlisting}
solve_schedule(epochs, architectures=None, *, cost_fn,
               switch_cost=0.0, hysteresis=0.0, max_iter=200,
               solve_kwargs=None) -> ScheduleResult
\end{lstlisting}
\begin{lstlisting}
from codesign import (AlgebraicDP, Architecture, Epoch, Ports,
                      Reals, System, solve_schedule)
def mk(name, c):
    s = System(name); d = s.provides("demand"); s.requires("cost")
    s.add("m", AlgebraicDP(F=Ports({"demand": Reals()}),
          R=Ports({"cost": Reals()}),
          equations={"cost": lambda f, c=c: c})).demand >= d
    s.constrain("cost", lambda x: x["m.cost"]); return s.build()
lo = Architecture("lo", mk("lo", 5.0))
hi = Architecture("hi", mk("hi", 3.0))
eps = [Epoch("e1", {"demand": 1.0}, candidates=[lo]),
       Epoch("e2", {"demand": 1.0}, candidates=[hi])]
r = solve_schedule(eps, cost_fn=lambda p: p["cost"], switch_cost=1.0)
print(r.schedule, r.total_cost, r.n_switches)   # -> ['lo', 'hi'] 9.0 1
\end{lstlisting}
\begin{itemize}[nosep]
  \item \texttt{epochs} --- \texttt{Sequence[Epoch]} --- the ordered
    environment regimes.
  \item \texttt{architectures} --- \texttt{Sequence[Architecture]} or
    \texttt{None} --- default candidate set for any epoch that does not
    supply its own; required if any epoch leaves \texttt{candidates} as
    \texttt{None}.
  \item \texttt{cost\_fn} --- \texttt{callable} --- maps a resource point to
    a scalar; lower is better.
  \item \texttt{switch\_cost} --- \texttt{float} or \texttt{callable} ---
    cost charged when the active architecture changes; a float applies
    uniformly, a callable \texttt{f(prev, next) -> float} makes it
    transition-dependent.
  \item \texttt{hysteresis} --- \texttt{float} --- extra margin a challenger
    must beat the incumbent by before a switch is preferred; suppresses
    chattering between near-equal options.
  \item \texttt{max\_iter} --- \texttt{int} --- forwarded to \texttt{solve()}
    for each epoch solve.
  \item \texttt{solve\_kwargs} --- \texttt{Mapping} or \texttt{None} ---
    extra keyword arguments forwarded to \texttt{solve()} (for example
    \texttt{uncertainty}).
\end{itemize}
Errors: raises \texttt{ValueError} when an epoch supplies no
\texttt{candidates} and no default \texttt{architectures} set was passed; the
message reads
\begin{lstlisting}[language={},numbers=none,keywordstyle={}]
epoch <name> has no candidates and no default architecture set
was supplied. Either set candidates=[...] on the epoch, or pass
architectures=[...] to solve_schedule().
\end{lstlisting}
An infeasible epoch is not an error: the schedule is returned with
\texttt{feasible=False} and the offending running cost recorded as
\texttt{inf} for inspection.

\subsubsection{\texttt{dynamic}}\label{ref:dynamic}

The \texttt{dynamic} module implements Case~2 with a \emph{scalar} carried
state (section~\ref{sec:temporal}): a finite-horizon dynamic program whose
per-stage decision is which architecture to instantiate, whose per-stage
cost is itself a co-design \texttt{solve()}, and whose carried scalar state
(fuel, charge, cumulative wear) is advanced by a user \texttt{transition}.
The value function is scalar --- at each (stage, state) the single best
cost-to-go --- found by a backward Bellman recursion. State is discretised
onto an explicit \texttt{StateGrid}; a transition that would leave the grid
envelope is rejected \emph{before} snapping, so an over-spent resource is
never silently rescued by the nearest in-range node. It reuses the
\texttt{Architecture} of the \texttt{temporal} module
(\ref{ref:temporal:architecture}).

\paragraph{\texttt{Stage}}\label{ref:dynamic:stage}
One step of the finite-horizon dynamic program: everything needed to score
and advance every candidate architecture at this point in time given the
carried state. A dataclass.
\begin{lstlisting}
Stage(name: str,
      functionality: Callable[[float], Mapping],
      transition: Callable[[float, Mapping], float],
      candidates: Optional[Sequence[Architecture]] = None,
      admissible: Optional[Callable[[float], bool]] = None)
\end{lstlisting}
\begin{lstlisting}
from codesign import Architecture, Stage
mode = Architecture("burn", None)
st = Stage("leg1", functionality=lambda s: {"demand": 1.0},
           transition=lambda s, p: s - p["fuel"],
           candidates=[mode], admissible=lambda s: s >= 0.0)
print(st.name, st.transition(10.0, {"fuel": 3.0}))   # -> leg1 7.0
\end{lstlisting}
\begin{itemize}[nosep]
  \item \texttt{name} --- \texttt{str} --- stage identifier.
  \item \texttt{functionality} --- \texttt{callable} ---
    \texttt{functionality(state) -> Mapping}, the outer functionality
    demanded when the carried resource has value \texttt{state}. A
    state-independent stage ignores its argument.
  \item \texttt{transition} --- \texttt{callable} ---
    \texttt{transition(state, point) -> new\_state}, mapping the incoming
    state and the chosen resource point to the outgoing state; the returned
    value is snapped to the nearest grid node.
  \item \texttt{candidates} --- \texttt{Sequence[Architecture]} or
    \texttt{None} --- architectures admissible here; falls back to the
    schedule-level default when \texttt{None}.
  \item \texttt{admissible} --- \texttt{callable} or \texttt{None} ---
    \texttt{admissible(state) -> bool}; states returning \texttt{False} are
    forbidden (infinite value).
\end{itemize}
Errors: none raised directly.

\paragraph{\texttt{StageResult}}\label{ref:dynamic:stageresult}
A per-stage record along a rolled-out optimal policy. A dataclass with no
public methods; constructed by \texttt{rollout}.
\begin{lstlisting}
from codesign import StageResult
sr = StageResult(stage="leg1", architecture="burn", state_in=10.0,
                 state_out=7.0, stage_cost=4.0, feasible=True,
                 point={"fuel": 3.0})
print(sr.stage, sr.state_in, sr.state_out, sr.stage_cost)
# -> leg1 10.0 7.0 4.0
\end{lstlisting}
\begin{itemize}[nosep]
  \item \texttt{stage} --- \texttt{str} --- the stage's name.
  \item \texttt{architecture} --- \texttt{str} --- chosen architecture name.
  \item \texttt{state\_in} --- \texttt{float} --- snapped carried state
    entering the stage.
  \item \texttt{state\_out} --- \texttt{float} --- snapped carried state
    leaving the stage.
  \item \texttt{stage\_cost} --- \texttt{float} --- co-design cost incurred
    at this stage.
  \item \texttt{feasible} --- \texttt{bool} --- whether a finite-cost action
    existed here.
  \item \texttt{point} --- \texttt{Mapping} or \texttt{None} --- the chosen
    resource point.
  \item \texttt{tags} --- \texttt{Mapping} --- the chosen architecture's
    tags.
\end{itemize}
Errors: none raised directly.

\paragraph{\texttt{StateGrid}}\label{ref:dynamic:stategrid}
A one-dimensional discretisation of the carried scalar resource. Tabular
dynamic programming needs a finite state set, so a continuous resource is
bucketed onto explicit, sorted nodes; transitions landing between nodes are
snapped to the nearest. Keeping the grid an explicit object makes the
discretisation visible and lets the caller trade accuracy against cost.
\begin{lstlisting}
StateGrid(nodes: Sequence[float])
\end{lstlisting}
\begin{lstlisting}
from codesign import StateGrid
g = StateGrid.linspace(0.0, 10.0, 6)   # nodes 0,2,4,6,8,10
print(list(g), len(g), g.snap(5.4))
# -> [0.0, 2.0, 4.0, 6.0, 8.0, 10.0] 6 6.0
\end{lstlisting}
\begin{itemize}[nosep]
  \item \texttt{nodes} --- \texttt{Sequence[float]} --- the grid node values
    (stored sorted). Also iterable, with \texttt{len()} giving the node
    count.
\end{itemize}
\begin{center}
\begin{tabular}{@{}lll@{}}
\toprule
Method & Signature & Behaviour \\
\midrule
\texttt{linspace} & \texttt{linspace(lo, hi, n)} (classmethod) &
  evenly spaced grid of \texttt{n} nodes on $[lo, hi]$. \\
\texttt{snap} & \texttt{snap(value) -> float} &
  the grid node nearest \texttt{value}. \\
\bottomrule
\end{tabular}
\end{center}
Errors: the constructor raises \texttt{ValueError} on an empty node
sequence (``\texttt{StateGrid needs at least one node, got an empty
sequence.\ \ldots}''); \texttt{linspace} raises \texttt{ValueError} when
\texttt{n < 1} (``\texttt{StateGrid.linspace needs at least one node, got
n=<n>.\ Pass n >= 1 \ldots}''). With \texttt{n == 1} a single node at
\texttt{lo} is returned.

\paragraph{\texttt{DynamicPolicy}}\label{ref:dynamic:dynamicpolicy}
A solved, state-indexed policy: the best architecture per (stage, state),
produced by \texttt{solve\_dynamic}. Beyond the single rolled-out path the
full table supports closed-loop control, re-queried at the actual realised
state when it diverges from the nominal roll-out.
\begin{lstlisting}
from codesign import (AlgebraicDP, Architecture, Ports, Reals, Stage,
                      StateGrid, System, solve_dynamic)
s = System("m"); d = s.provides("demand"); s.requires("fuel")
s.add("m", AlgebraicDP(F=Ports({"demand": Reals()}),
      R=Ports({"fuel": Reals()}),
      equations={"fuel": lambda f: 2.0})).demand >= d
s.constrain("fuel", lambda x: x["m.fuel"])
burn = Architecture("burn", s.build())
st = Stage("leg", lambda s: {"demand": 1.0},
           lambda s, p: s - p["fuel"], candidates=[burn])
policy = solve_dynamic([st], StateGrid.linspace(0.0, 10.0, 11),
                       cost_fn=lambda p: p["fuel"])
print(policy.cost_to_go(0, 8.0), policy.action_at(0, 8.0))
# -> 2.0 ('burn', {'fuel': 2.0}, 2.0, 6.0)
\end{lstlisting}
\begin{center}
\begin{tabular}{@{}lp{0.62\linewidth}@{}}
\toprule
Method & Behaviour \\
\midrule
\texttt{cost\_to\_go(t, state)} & optimal remaining cost from stage
  \texttt{t} at the grid node nearest \texttt{state} (\texttt{inf} if
  none). \\
\texttt{action\_at(t, state)} & best
  \texttt{(arch, point, stage\_cost, state\_out)} at (\texttt{t},
  \texttt{state}); \texttt{None} if no finite-cost action exists. \\
\bottomrule
\end{tabular}
\end{center}
Errors: none raised directly (an out-of-table state returns \texttt{inf} /
\texttt{None} rather than raising).

\paragraph{\texttt{DynamicResult}}\label{ref:dynamic:dynamicresult}
The outcome of a policy rolled out from an initial state, returned by
\texttt{rollout} and \texttt{solve\_and\_rollout}. A dataclass with a
\texttt{schedule} property and a compact \texttt{\_\_repr\_\_}; it also
carries the full \texttt{policy} for closed-loop re-querying.
\begin{lstlisting}
from codesign import (AlgebraicDP, Architecture, Ports, Reals, Stage,
                      StateGrid, System, solve_and_rollout)
s = System("m"); d = s.provides("demand"); s.requires("fuel")
s.add("m", AlgebraicDP(F=Ports({"demand": Reals()}),
      R=Ports({"fuel": Reals()}),
      equations={"fuel": lambda f: 2.0})).demand >= d
s.constrain("fuel", lambda x: x["m.fuel"])
burn = Architecture("burn", s.build())
st = Stage("leg", lambda s: {"demand": 1.0},
           lambda s, p: s - p["fuel"], candidates=[burn])
grid = StateGrid.linspace(0.0, 10.0, 11)
res = solve_and_rollout([st, st], grid, 10.0, cost_fn=lambda p: p["fuel"])
print(res.schedule, res.total_cost, res.feasible)
# -> ['burn', 'burn'] 4.0 True
\end{lstlisting}
\begin{itemize}[nosep]
  \item \texttt{stages} --- \texttt{List[StageResult]} --- per-stage record
    of the optimal policy from the initial state.
  \item \texttt{total\_cost} --- \texttt{float} --- optimal cost-to-go from
    the initial state (running plus terminal).
  \item \texttt{feasible} --- \texttt{bool} --- \texttt{True} iff a
    finite-cost policy exists from the initial state.
  \item \texttt{policy} --- \texttt{DynamicPolicy} --- the full
    state-indexed policy table.
  \item \texttt{schedule} (property) --- \texttt{List[str]} --- architecture
    chosen at each stage along the rolled-out path.
\end{itemize}
Errors: none raised directly.

\paragraph{\texttt{solve\_dynamic}}\label{ref:dynamic:solve_dynamic}
Solves the finite-horizon architecture DP by a backward Bellman pass over
the (stage, state) lattice, returning a full state-indexed
\texttt{DynamicPolicy} (\ref{ref:dynamic:dynamicpolicy}). At each (stage,
state) it solves the co-design problem for each admissible architecture at
the stage's state-dependent functionality, reads the cheapest feasible
point, advances the state, and adds the discretised cost-to-go of the
successor.
\begin{lstlisting}
solve_dynamic(stages, grid, *, cost_fn, architectures=None,
              terminal_cost=None, max_iter=200, solve_kwargs=None,
              cache=True) -> DynamicPolicy
\end{lstlisting}
\begin{lstlisting}
from codesign import (AlgebraicDP, Architecture, Ports, Reals, Stage,
                      StateGrid, System, solve_dynamic)
s = System("m"); d = s.provides("demand"); s.requires("fuel")
s.add("m", AlgebraicDP(F=Ports({"demand": Reals()}),
      R=Ports({"fuel": Reals()}),
      equations={"fuel": lambda f: 2.0})).demand >= d
s.constrain("fuel", lambda x: x["m.fuel"])
burn = Architecture("burn", s.build())
st = Stage("leg", lambda s: {"demand": 1.0},
           lambda s, p: s - p["fuel"], candidates=[burn])
grid = StateGrid.linspace(0.0, 10.0, 11)
policy = solve_dynamic([st, st], grid, cost_fn=lambda p: p["fuel"])
print(policy.stage_names, policy.cost_to_go(0, 10.0))
# -> ['leg', 'leg'] 4.0
\end{lstlisting}
\begin{itemize}[nosep]
  \item \texttt{stages} --- \texttt{Sequence[Stage]} --- the horizon in
    forward order (solved backward internally).
  \item \texttt{grid} --- \texttt{StateGrid} --- discretisation of the
    carried scalar resource.
  \item \texttt{cost\_fn} --- \texttt{callable} --- scalar cost on a
    resource point; lower is better.
  \item \texttt{architectures} --- \texttt{Sequence[Architecture]} or
    \texttt{None} --- default candidates for stages without their own.
  \item \texttt{terminal\_cost} --- \texttt{callable} or \texttt{None} ---
    \texttt{terminal\_cost(state) -> float} scoring the state after the
    final stage; defaults to zero everywhere.
  \item \texttt{max\_iter} --- \texttt{int} --- forwarded to
    \texttt{solve()}.
  \item \texttt{solve\_kwargs} --- \texttt{Mapping} or \texttt{None} ---
    extra keyword arguments forwarded to \texttt{solve()}.
  \item \texttt{cache} --- \texttt{bool} --- cache co-design solves keyed by
    (architecture name, functionality); safe when
    \texttt{functionality(state)} is deterministic.
\end{itemize}
Errors: raises \texttt{ValueError} when a stage supplies no
\texttt{candidates} and no default \texttt{architectures} set was passed; the
message reads
\begin{lstlisting}[language={},numbers=none,keywordstyle={}]
stage <name> has no candidates and no default architecture set
was supplied. Either set candidates=[...] on the stage, or pass
architectures=[...] to solve_dynamic().
\end{lstlisting}

\paragraph{\texttt{rollout}}\label{ref:dynamic:rollout}
Rolls a solved policy forward from a concrete initial state, following the
policy's chosen architecture at each stage, advancing the snapped carried
state, and assembling a \texttt{DynamicResult}
(\ref{ref:dynamic:dynamicresult}).
\begin{lstlisting}
rollout(policy, stages, initial_state, *, arch_lookup=None)
    -> DynamicResult
\end{lstlisting}
\begin{lstlisting}
from codesign import (AlgebraicDP, Architecture, Ports, Reals, Stage,
                      StateGrid, System, rollout, solve_dynamic)
s = System("m"); d = s.provides("demand"); s.requires("fuel")
s.add("m", AlgebraicDP(F=Ports({"demand": Reals()}),
      R=Ports({"fuel": Reals()}),
      equations={"fuel": lambda f: 2.0})).demand >= d
s.constrain("fuel", lambda x: x["m.fuel"])
burn = Architecture("burn", s.build())
st = Stage("leg", lambda s: {"demand": 1.0},
           lambda s, p: s - p["fuel"], candidates=[burn])
grid = StateGrid.linspace(0.0, 10.0, 11)
policy = solve_dynamic([st, st], grid, cost_fn=lambda p: p["fuel"])
res = rollout(policy, [st, st], 10.0)
print(res.schedule, res.total_cost)   # -> ['burn', 'burn'] 4.0
\end{lstlisting}
\begin{itemize}[nosep]
  \item \texttt{policy} --- \texttt{DynamicPolicy} --- a policy from
    \texttt{solve\_dynamic}.
  \item \texttt{stages} --- \texttt{Sequence[Stage]} --- the same stage list
    passed to \texttt{solve\_dynamic}.
  \item \texttt{initial\_state} --- \texttt{float} --- the starting carried
    value.
  \item \texttt{arch\_lookup} --- \texttt{Mapping} or \texttt{None} --- maps
    architecture name to \texttt{Architecture} so tags attach to each
    record; built automatically from the stage candidate lists when
    omitted.
\end{itemize}
Errors: none raised directly. When the policy has no feasible action from
the reached state the roll-out stops early, the offending record is marked
\texttt{feasible=False}, and \texttt{DynamicResult.feasible} becomes
\texttt{False}.

\paragraph{\texttt{solve\_and\_rollout}}\label{ref:dynamic:solve_and_rollout}
A convenience wrapper equivalent to \texttt{solve\_dynamic} followed by
\texttt{rollout} from \texttt{initial\_state}. Prefer the two-step form when
the policy will be queried at several initial states or in closed loop.
\begin{lstlisting}
solve_and_rollout(stages, grid, initial_state, *, cost_fn,
                  architectures=None, terminal_cost=None, max_iter=200,
                  solve_kwargs=None, cache=True) -> DynamicResult
\end{lstlisting}
\begin{lstlisting}
from codesign import (AlgebraicDP, Architecture, Ports, Reals, Stage,
                      StateGrid, System, solve_and_rollout)
s = System("m"); d = s.provides("demand"); s.requires("fuel")
s.add("m", AlgebraicDP(F=Ports({"demand": Reals()}),
      R=Ports({"fuel": Reals()}),
      equations={"fuel": lambda f: 2.0})).demand >= d
s.constrain("fuel", lambda x: x["m.fuel"])
burn = Architecture("burn", s.build())
st = Stage("leg", lambda s: {"demand": 1.0},
           lambda s, p: s - p["fuel"], candidates=[burn])
grid = StateGrid.linspace(0.0, 10.0, 11)
res = solve_and_rollout([st, st], grid, 10.0, cost_fn=lambda p: p["fuel"])
print(res.schedule, res.total_cost, res.feasible)
# -> ['burn', 'burn'] 4.0 True
\end{lstlisting}
Parameters: as \texttt{solve\_dynamic} (\ref{ref:dynamic:solve_dynamic}),
plus the positional \texttt{initial\_state} (\texttt{float}) rolled out
from.
Errors: same \texttt{ValueError} as \texttt{solve\_dynamic} when candidates
are missing.

\subsubsection{\texttt{sequential}}\label{ref:sequential}

The \texttt{sequential} module is the \emph{antichain-valued} generalisation
of the scalar \texttt{dynamic} layer (section~\ref{sec:sequential}): a
finite-horizon decision process whose stages are co-design problems coupled
by a carried state, solved by a backward Bellman recursion whose value at
each (stage, state) is a Pareto antichain rather than a scalar. The
accumulated resource (the named cost axes on the antichain) is kept
deliberately distinct from the carried state: the transition reads any
carried quantity off the full solved point, while only the cost axes
accumulate on the front. The scalar \texttt{dynamic} layer is the width-one
special case. This module also makes three theory results operational ---
monotone value (Q1, \texttt{check\_monotonicity}), front $=$ reachable
frontier (Q2), and exact factorisation at a reset (Q3,
\texttt{detect\_resets} / \texttt{factorise\_at\_resets}) --- and exposes
the precompute-then-DP structure (section~\ref{sec:precompute}). It reuses
the \texttt{Architecture} (\ref{ref:temporal:architecture}) and
\texttt{StateGrid} (\ref{ref:dynamic:stategrid}).

\paragraph{\texttt{SeqStage}}\label{ref:sequential:seqstage}
One stage of a sequential co-design problem. Structurally the antichain-
valued analogue of \texttt{Stage} (\ref{ref:dynamic:stage}): returning the
grid's bottom node from \texttt{transition} makes the stage a reset. A
dataclass.
\begin{lstlisting}
SeqStage(name: str,
         functionality: Callable[[float], Mapping],
         transition: Callable[[float, Mapping], float],
         candidates: Optional[Sequence[Architecture]] = None,
         admissible: Optional[Callable[[float], bool]] = None)
\end{lstlisting}
\begin{lstlisting}
from codesign import Architecture, SeqStage
mode = Architecture("m", None)
st = SeqStage("leg", functionality=lambda x: {"demand": 1.0},
              transition=lambda x, p: x - p["fuel"],
              candidates=[mode], admissible=lambda x: x >= 0.0)
print(st.name, st.transition(6.0, {"fuel": 2.0}))   # -> leg 4.0
\end{lstlisting}
\begin{itemize}[nosep]
  \item \texttt{name} --- \texttt{str} --- stage identifier.
  \item \texttt{functionality} --- \texttt{callable} ---
    \texttt{functionality(state) -> Mapping}, the outer functionality
    demanded when the carried state is \texttt{state}.
  \item \texttt{transition} --- \texttt{callable} ---
    \texttt{transition(state, point) -> new\_state}, from the incoming state
    and a chosen resource point (over the resource ports) to the outgoing
    state.
  \item \texttt{candidates} --- \texttt{Sequence[Architecture]} or
    \texttt{None} --- architectures admissible here; each one's co-design
    solve supplies part of the stage antichain $h_k(x)$ (the union of their
    fronts). Falls back to the problem-level default when \texttt{None}.
  \item \texttt{admissible} --- \texttt{callable} or \texttt{None} ---
    \texttt{admissible(state) -> bool}; forbidden states carry the empty
    (infeasible) value.
\end{itemize}
Errors: none raised directly.

\paragraph{\texttt{sum\_combine}}\label{ref:sequential:sum_combine}
The consumable/accumulating monoid combination $\oplus$: component-wise sum
over the cost axes. Resources accumulate across stages (fuel burned, money
spent, wear incurred); the reachable frontier can grow polynomially with the
horizon. The default \texttt{combine} for \texttt{solve\_sequential}.
\begin{lstlisting}
from codesign import sum_combine
a = {"cost": 5.0, "co2": 1.0}
b = {"cost": 2.0, "co2": 8.0}
print(sum_combine(a, b))           # -> {'cost': 7.0, 'co2': 9.0}
\end{lstlisting}
\begin{itemize}[nosep]
  \item \texttt{a}, \texttt{b} --- \texttt{Mapping} --- two resource points
    over the same axes; returns a \texttt{dict} of their component-wise sum.
\end{itemize}
Errors: none raised directly (a \texttt{KeyError} results if \texttt{b}
lacks an axis present in \texttt{a}).

\paragraph{\texttt{join\_combine}}\label{ref:sequential:join_combine}
The renewable monoid combination $\oplus$: component-wise maximum (join).
The requirement of two stages is the peak, not the sum (a worker, an oven, a
bus reused across stages); the reachable frontier stays bounded uniformly in
the horizon, and if every stage resets the problem collapses to a static
co-design.
\begin{lstlisting}
from codesign import join_combine
a = {"cost": 5.0, "co2": 1.0}
b = {"cost": 2.0, "co2": 8.0}
print(join_combine(a, b))          # -> {'cost': 5.0, 'co2': 8.0}
\end{lstlisting}
\begin{itemize}[nosep]
  \item \texttt{a}, \texttt{b} --- \texttt{Mapping} --- two resource points
    over the same axes; returns a \texttt{dict} of their component-wise
    maximum.
\end{itemize}
Errors: none raised directly (a \texttt{KeyError} results if \texttt{b}
lacks an axis present in \texttt{a}).

\paragraph{\texttt{solve\_sequential}}\label{ref:sequential:solve_sequential}
Solves the antichain-valued sequential co-design problem by the backward
Bellman recursion in the upper-set value space, carrying a full Pareto
antichain of cumulative resources at each (stage, state). Returns a
\texttt{SeqResult} (\ref{ref:sequential:seqresult}) whose value is
$V_0(\texttt{initial\_state})$. The snippet below uses a single architecture
(a width-one front); with several incomparable candidates the value has
width greater than one --- see \texttt{dp\_over\_catalog}
(\ref{ref:sequential:dp_over_catalog}) for a multi-point front.
\begin{lstlisting}
solve_sequential(stages, grid, *, cost_axes, initial_state,
                 combine=sum_combine, architectures=None,
                 max_iter=200, solve_kwargs=None) -> SeqResult
\end{lstlisting}
\begin{lstlisting}
from codesign import (AlgebraicDP, Architecture, Ports, Reals, SeqStage,
                      StateGrid, System, solve_sequential)
s = System("m"); d = s.provides("demand"); s.requires("fuel")
s.add("m", AlgebraicDP(F=Ports({"demand": Reals()}),
      R=Ports({"fuel": Reals()}),
      equations={"fuel": lambda f: 2.0})).demand >= d
s.constrain("fuel", lambda x: x["m.fuel"])
a = Architecture("a", s.build())
st = SeqStage("leg", lambda x: {"demand": 1.0},
              lambda x, p: x - p["fuel"], candidates=[a],
              admissible=lambda x: x >= -1e-9)
res = solve_sequential([st, st], StateGrid.linspace(0., 6., 13),
                       cost_axes=["fuel"], initial_state=6.)
print(res.width, sorted(p["fuel"] for p in res.value), res.feasible)
# -> 1 [4.0] True
\end{lstlisting}
\begin{itemize}[nosep]
  \item \texttt{stages} --- \texttt{Sequence[SeqStage]} --- the horizon in
    forward order (solved backward internally).
  \item \texttt{grid} --- \texttt{StateGrid} --- discretisation of the
    carried scalar state.
  \item \texttt{cost\_axes} --- \texttt{Sequence[str]} --- names of the
    resource ports to accumulate on the antichain (for example
    \texttt{["cost", "co2"]}); these define the product poset $R$.
  \item \texttt{initial\_state} --- \texttt{float} --- starting carried
    state; the returned value is $V_0$ of it.
  \item \texttt{combine} --- \texttt{callable} --- the monoid combination
    $\oplus$: \texttt{sum\_combine} (default) for a consumable resource,
    \texttt{join\_combine} for a renewable one.
  \item \texttt{architectures} --- \texttt{Sequence[Architecture]} or
    \texttt{None} --- default candidates for stages without their own.
  \item \texttt{max\_iter} --- \texttt{int} --- forwarded to
    \texttt{solve()}.
  \item \texttt{solve\_kwargs} --- \texttt{Mapping} or \texttt{None} ---
    extra keyword arguments forwarded to \texttt{solve()}.
\end{itemize}
Errors: raises \texttt{ValueError} when a stage supplies no
\texttt{candidates} and no default \texttt{architectures} set was passed; the
message reads
\begin{lstlisting}[language={},numbers=none,keywordstyle={}]
stage <name> has no candidates and no default architecture set
was supplied. Either set candidates=[...] on the stage, or pass
architectures=[...] to the sequential solver (solve_sequential /
check_monotonicity).
\end{lstlisting}

\paragraph{\texttt{SeqResult}}\label{ref:sequential:seqresult}
The outcome of \texttt{solve\_sequential} (and \texttt{dp\_over\_catalog}):
the whole-horizon Pareto front from the initial state. A dataclass with a
compact \texttt{\_\_repr\_\_}. The snippet obtains one from
\texttt{dp\_over\_catalog} to show a genuine multi-point front.
\begin{lstlisting}
from codesign import StateGrid, dp_over_catalog
cat = [("clean", {"cost": 10.0, "co2": 1.0}),
       ("cheap", {"cost": 2.0, "co2": 8.0})]
res = dp_over_catalog([cat, cat], StateGrid.linspace(0., 4., 5),
                      cost_axes=["cost", "co2"], initial_state=4.0,
                      transition=lambda s, p: s)
print(res.width, res.feasible,
      sorted((p["cost"], p["co2"]) for p in res.value))
# -> 3 True [(4.0, 16.0), (12.0, 9.0), (20.0, 2.0)]
\end{lstlisting}
\begin{itemize}[nosep]
  \item \texttt{value} --- \texttt{Antichain} --- the value antichain
    $\mathrm{Min}\,V_0(\texttt{initial\_state})$; each point is a
    \texttt{Mapping} over the cost axes.
  \item \texttt{width} --- \texttt{int} --- \texttt{len(value)}, the number
    of incomparable Pareto-optimal totals (the reachable-frontier width).
  \item \texttt{feasible} --- \texttt{bool} --- \texttt{True} iff the value
    antichain is non-empty and free of tops.
  \item \texttt{policy} --- \texttt{SeqPolicy} --- the full state-indexed
    antichain-valued policy table.
\end{itemize}
Errors: none raised directly.

\paragraph{\texttt{SeqPolicy}}\label{ref:sequential:seqpolicy}
The state-indexed antichain-valued policy from the backward pass, held on
\texttt{SeqResult.\allowbreak policy}. For each (stage, state node) it stores
the value
antichain $V_k(x)$ and, per Pareto point, the realising architecture and
chosen stage point, so an optimal choice sequence can be traced forward.
\begin{lstlisting}
from codesign import StateGrid, dp_over_catalog
cat = [("clean", {"cost": 10.0, "co2": 1.0}),
       ("cheap", {"cost": 2.0, "co2": 8.0})]
res = dp_over_catalog([cat, cat], StateGrid.linspace(0., 4., 5),
                      cost_axes=["cost", "co2"], initial_state=4.0,
                      transition=lambda s, p: s)
pol = res.policy
print(pol.width_at(0, 4.0), pol.stage_names)
# -> 3 ['stage_0', 'stage_1']
\end{lstlisting}
\begin{center}
\begin{tabular}{@{}lp{0.58\linewidth}@{}}
\toprule
Method & Behaviour \\
\midrule
\texttt{value\_at(k, state)} & the value antichain $V_k(x)$ at the grid node
  nearest \texttt{state} (empty antichain if absent). \\
\texttt{width\_at(k, state)} & the reachable-frontier width $\alpha_k(x)$,
  i.e.\ \texttt{len(value\_at(k, state))}. \\
\bottomrule
\end{tabular}
\end{center}
Errors: none raised directly.

\paragraph{\texttt{MonotonicityReport}}\label{ref:sequential:monotonicityreport}
The result of \texttt{check\_monotonicity}: whether the three hypotheses of
the monotone-value theorem --- (H1), the \emph{joint} (H2), and (H3) --- held
on the grid, with sampled witnesses when they did not. A dataclass whose
\texttt{monotone\_\allowbreak value\_\allowbreak guaranteed} property is
\texttt{True} exactly when \emph{all} hold (the Q1 theorem then guarantees a
monotone value).
\begin{lstlisting}
from codesign import MonotonicityReport
rep = MonotonicityReport(h1_ok=True, h2_ok=True, h2_joint_ok=True, h3_ok=True)
print(rep.monotone_value_guaranteed, repr(rep))
# -> True MonotonicityReport(H1=ok, H2=ok, H2joint=ok, H3=ok,
#         value_monotone=guaranteed)
\end{lstlisting}
\begin{itemize}[nosep]
  \item \texttt{h1\_ok} --- \texttt{bool} --- (H1) held: the stage map
    \texttt{h\_k} is consistently oriented at every tested grid pair.
  \item \texttt{h2\_ok} --- \texttt{bool} --- the \emph{state slice} of (H2)
    held: \texttt{phi\_k(.,\ r)} is monotone in the state at every pair.
  \item \texttt{h2\_joint\_ok} --- \texttt{bool} --- the \emph{resource
    slice} of (H2) held, so the transition is \emph{jointly} monotone on
    \texttt{X x R} (the paper's (H2-joint)): a larger resource sends the
    successor no lower in the certified orientation.
  \item \texttt{h3\_ok} --- \texttt{bool} --- (H3) held: the admissible
    region is a down-set of \texttt{X} in the certified orientation.
  \item \texttt{h1\_violations}, \texttt{h2\_violations},
    \texttt{h2\_joint\_violations}, \texttt{h3\_violations} ---
    \texttt{List[Tuple[str, float, float]]} --- sampled witnesses (the
    \texttt{h2\_joint} witness is the out-of-order successor pair
    \texttt{(stage, phi\_x\_r, phi\_x\_rp)}; the others are
    \texttt{(stage, x, x')}), empty when the condition held.
  \item \texttt{monotone\_value\_guaranteed} (property) ---
    \texttt{h1\_ok and h2\_ok and h2\_joint\_ok and h3\_ok}.
\end{itemize}
Errors: none raised directly.

\paragraph{\texttt{check\_monotonicity}}\label{ref:sequential:check_monotonicity}
Numerically verifies the three hypotheses of the monotone-value theorem on
the state grid --- (H1) consistent orientation of \texttt{h\_k}, the
\emph{joint} (H2) monotonicity of the transition on \texttt{X x R} (both its
state slice \texttt{h2\_ok} and its resource slice \texttt{h2\_joint\_ok}),
and the (H3) down-set property of admissibility --- and returns a
\texttt{MonotonicityReport} (\ref{ref:sequential:monotonicityreport}). It is
orientation-aware: it accepts the benign ``consumable but monotone''
orientation (a budget carried as remaining slack) and flags only genuinely
non-monotone stages (perishable / fatigue-as-state, with an interior
optimum), a transition that is not jointly monotone, or a non-down-set
admissible region, all of which fall outside the Q1 guarantee. Spurious
violations most often signal a too-coarse grid, since snapping is not
order-preserving at bucket boundaries.
\begin{lstlisting}
check_monotonicity(stages, grid, *, cost_axes, architectures=None,
                   max_iter=200, solve_kwargs=None,
                   max_violations=8) -> MonotonicityReport
\end{lstlisting}
\begin{lstlisting}
from codesign import (AlgebraicDP, Architecture, Ports, Reals, SeqStage,
                      StateGrid, System, check_monotonicity)
s = System("m"); d = s.provides("demand"); s.requires("fuel")
s.add("m", AlgebraicDP(F=Ports({"demand": Reals()}),
      R=Ports({"fuel": Reals()}),
      equations={"fuel": lambda f: 2.0})).demand >= d
s.constrain("fuel", lambda x: x["m.fuel"])
a = Architecture("a", s.build())
st = SeqStage("leg", lambda x: {"demand": 1.0},
              lambda x, p: x - p["fuel"], candidates=[a],
              admissible=lambda x: x >= -1e-9)
rep = check_monotonicity([st, st], StateGrid.linspace(0., 6., 13),
                         cost_axes=["fuel"])
print(rep.h1_ok, rep.h2_ok, rep.h2_joint_ok, rep.h3_ok,
      rep.monotone_value_guaranteed)
# -> True True True True True
\end{lstlisting}
\begin{itemize}[nosep]
  \item \texttt{stages}, \texttt{grid}, \texttt{cost\_axes},
    \texttt{architectures}, \texttt{max\_iter}, \texttt{solve\_kwargs} ---
    as \texttt{solve\_sequential} (\ref{ref:sequential:solve_sequential}).
  \item \texttt{max\_violations} --- \texttt{int} --- cap on the number of
    sampled witnesses recorded per hypothesis.
\end{itemize}
Errors: raises the same \texttt{ValueError} as
\texttt{solve\_sequential} when a stage has no candidates and no default set
(the message names both ``\texttt{solve\_sequential / check\_monotonicity}'').

\paragraph{\texttt{detect\_resets}}\label{ref:sequential:detect_resets}
Returns the indices of stages that are resets: a stage is a reset iff its
transition lands on the grid's bottom node for every incoming state and
every chosen resource point. At a reset the horizon factorises exactly (Q3),
in every regime, since the proof uses only distributivity.
\begin{lstlisting}
detect_resets(stages, grid, *, cost_axes, architectures=None,
              max_iter=200, solve_kwargs=None) -> List[int]
\end{lstlisting}
\begin{lstlisting}
from codesign import (AlgebraicDP, Architecture, Ports, Reals, SeqStage,
                      StateGrid, System, detect_resets)
s = System("m"); d = s.provides("demand"); s.requires("fuel")
s.add("m", AlgebraicDP(F=Ports({"demand": Reals()}),
      R=Ports({"fuel": Reals()}),
      equations={"fuel": lambda f: 2.0})).demand >= d
s.constrain("fuel", lambda x: x["m.fuel"])
a = Architecture("a", s.build())
run = SeqStage("run", lambda x: {"demand": 1.0},
               lambda x, p: x - p["fuel"], candidates=[a])
reset = SeqStage("reset", lambda x: {"demand": 1.0},
                 lambda x, p: 0.0, candidates=[a])
print(detect_resets([run, reset, run],
                    StateGrid.linspace(0., 6., 13), cost_axes=["fuel"]))
# -> [1]
\end{lstlisting}
\begin{itemize}[nosep]
  \item \texttt{stages}, \texttt{grid}, \texttt{cost\_axes},
    \texttt{architectures}, \texttt{max\_iter}, \texttt{solve\_kwargs} ---
    as \texttt{solve\_sequential} (\ref{ref:sequential:solve_sequential}).
\end{itemize}
Errors: raises \texttt{ValueError} (naming the stage) when a stage has no
\texttt{candidates} and no default \texttt{architectures} was supplied ---
the same guard as \texttt{solve\_sequential} and \texttt{check\_monotonicity},
so an under-specified stage fails loudly instead of being silently (and
vacuously) reported as a reset.

\paragraph{\texttt{factorise\_at\_resets}}\label{ref:sequential:factorise_at_resets}
Partitions the horizon into maximal quiescence-free runs given the reset
indices, returning inclusive \texttt{(start, end)} ranges. By Q3 the value
over the full horizon is the $\oplus$-product of the values of these
independent sub-problems --- the order-theoretic, deterministic analogue of
regeneration-point decomposition. Pure bookkeeping: it solves nothing.
\begin{lstlisting}
factorise_at_resets(stages, resets) -> List[Tuple[int, int]]
\end{lstlisting}
\begin{lstlisting}
from codesign import SeqStage, factorise_at_resets
f = lambda x: {"demand": 1.0}
t = lambda x, p: x
stages = [SeqStage(f"s{i}", f, t) for i in range(4)]
print(factorise_at_resets(stages, resets=[1]))
# -> [(0, 1), (2, 3)]
\end{lstlisting}
\begin{itemize}[nosep]
  \item \texttt{stages} --- \texttt{Sequence[SeqStage]} --- the horizon
    (only its length is used).
  \item \texttt{resets} --- \texttt{Sequence[int]} --- reset stage indices,
    typically from \texttt{detect\_resets}
    (\ref{ref:sequential:detect_resets}). A reset closes the run containing
    it.
\end{itemize}
Errors: none raised directly.

\paragraph{\texttt{precompute\_catalog}}\label{ref:sequential:precompute_catalog}
Solves every architecture once at a fixed functionality and returns a flat,
tagged Pareto catalog --- the co-design precomputation step of the Formula~1
seasonal framework (section~\ref{sec:precompute}). Points that look dominated
at this single-stage level are retained if non-dominated in the combined
antichain, since they may become optimal once aggregated in an outer DP.
\begin{lstlisting}
precompute_catalog(architectures, functionality, cost_axes, *,
                   max_iter=200, solve_kwargs=None)
    -> List[Tuple[str, Mapping]]
\end{lstlisting}
\begin{lstlisting}
from codesign import (AlgebraicDP, Architecture, Ports, Reals,
                      System, precompute_catalog)
s = System("m"); d = s.provides("demand"); s.requires("fuel")
s.add("m", AlgebraicDP(F=Ports({"demand": Reals()}),
      R=Ports({"fuel": Reals()}),
      equations={"fuel": lambda f: 2.0})).demand >= d
s.constrain("fuel", lambda x: x["m.fuel"])
a = Architecture("a", s.build())
cat = precompute_catalog([a], {"demand": 1.0}, cost_axes=["fuel"])
print(cat)                         # -> [('a', {'fuel': 2.0})]
\end{lstlisting}
\begin{itemize}[nosep]
  \item \texttt{architectures} --- \texttt{Sequence[Architecture]} --- the
    candidates to solve once.
  \item \texttt{functionality} --- \texttt{Mapping} --- the fixed outer
    functionality at which each architecture is solved.
  \item \texttt{cost\_axes} --- \texttt{Sequence[str]} --- resource ports
    accumulated on the catalog's antichain.
  \item \texttt{max\_iter}, \texttt{solve\_kwargs} --- forwarded to
    \texttt{solve()}.
\end{itemize}
Returns a list of \texttt{(arch\_name, point)} pairs (each \texttt{point} the
full solved resource point, so a carried-state axis remains readable);
feed it to \texttt{dp\_over\_catalog} (\ref{ref:sequential:dp_over_catalog}).
Errors: none raised directly.

\paragraph{\texttt{dp\_over\_catalog}}\label{ref:sequential:dp_over_catalog}
Runs the antichain-valued DP that \emph{selects from precomputed catalogs}:
each stage's action set is a frozen catalog rather than a live co-design
solve, so no \texttt{solve()} happens inside the Bellman sweep. Use it when
the per-stage co-design is independent of the carried state (the Formula~1
regime); otherwise use \texttt{solve\_sequential}. Returns a
\texttt{SeqResult} (\ref{ref:sequential:seqresult}).
\begin{lstlisting}
dp_over_catalog(catalogs, grid, *, cost_axes, initial_state,
                transition, combine=sum_combine,
                admissible=None) -> SeqResult
\end{lstlisting}
\begin{lstlisting}
from codesign import StateGrid, dp_over_catalog, sum_combine
cat = [("clean", {"cost": 10.0, "co2": 1.0}),
       ("cheap", {"cost": 2.0, "co2": 8.0})]
res = dp_over_catalog([cat, cat], StateGrid.linspace(0., 4., 5),
                      cost_axes=["cost", "co2"], initial_state=4.0,
                      transition=lambda s, p: s, combine=sum_combine)
print(res.width, sorted((p["cost"], p["co2"]) for p in res.value))
# -> 3 [(4.0, 16.0), (12.0, 9.0), (20.0, 2.0)]
\end{lstlisting}
\begin{itemize}[nosep]
  \item \texttt{catalogs} --- \texttt{Sequence[Sequence[Tuple[str,
    Mapping]]]} --- one catalog per stage, each a sequence of
    \texttt{(arch\_name, point)} pairs from \texttt{precompute\_catalog}.
  \item \texttt{grid} --- \texttt{StateGrid} --- discretisation of the
    carried scalar state.
  \item \texttt{cost\_axes} --- \texttt{Sequence[str]} --- resource ports
    accumulated on the antichain.
  \item \texttt{initial\_state} --- \texttt{float} --- starting carried
    state.
  \item \texttt{transition} --- \texttt{callable} ---
    \texttt{transition(state, point) -> new\_state}.
  \item \texttt{combine} --- \texttt{callable} --- monoid combination
    $\oplus$ (default \texttt{sum\_combine}).
  \item \texttt{admissible} --- \texttt{callable} or \texttt{None} --- state
    admissibility predicate.
\end{itemize}
Errors: none raised directly.

\subsection{State and closed loop}

The vector-state generalisation (\texttt{state}, \texttt{vector\_dp}) and
the closed-loop \texttt{online\_codesign} controller. See
Sections~\ref{sec:vector} and~\ref{sec:online-cd}.

\subsubsection{\texttt{state}}
\label{ref:state}

The \texttt{state} module supplies the general carried state of the temporal
layer (section~\ref{sec:vector}): a hashable state-vector value type and a
product grid of named axes carrying the component-wise order the
monotonicity results are stated over. A one-axis grid reproduces the scalar
\texttt{StateGrid} exactly, so the vector layer strictly generalises the
scalar carried state rather than parallelling it.

\paragraph{\texttt{StateVec}}\label{ref:state:statevec}
The hashable value type of a carried state vector: a tuple of
\texttt{(axis\_name, value)} pairs sorted by name, so two states with the
same contents compare and hash equal. It is a type alias, not a class ---
\texttt{StateVec = Tuple[Tuple[str, Any], ...]} --- and is used directly as
the key of the dynamic-program table. Build one with \texttt{make\_state}
and read it with \texttt{state\_get} / \texttt{state\_as\_dict} rather than
by indexing the raw tuple.

\begin{lstlisting}
from codesign import make_state
sv = make_state(fuel=12.0, flag=0)  # build a StateVec
print(sv)          # (('flag', 0), ('fuel', 12.0)): sorted, hashable
table = {sv: 0.0}  # usable directly as a DP-table key
print(table[make_state(flag=0, fuel=12.0)])  # 0.0: order-independent
\end{lstlisting}

\begin{itemize}[nosep]
  \item Contents --- one \texttt{(str, Any)} pair per axis; canonicalised
    (sorted by name) so keyword order at construction is irrelevant.
\end{itemize}
Errors: none raised directly (it is a plain tuple).

\paragraph{\texttt{make\_state}}\label{ref:state:make_state}
Build a canonical \texttt{StateVec} from keyword axis values.

\begin{lstlisting}
make_state(**axis_values) -> StateVec
\end{lstlisting}

\begin{lstlisting}
from codesign import make_state
x = make_state(fuel=12.0, flag=0)
print(x)                     # (('flag', 0), ('fuel', 12.0))
print(isinstance(x, tuple))  # True: hashable, usable as a DP-table key
\end{lstlisting}

\begin{itemize}[nosep]
  \item \texttt{**axis\_values} --- each keyword names an axis and gives its
    value; the result is normalised so equal contents hash and compare
    equal.
\end{itemize}
Errors: none raised directly.

\paragraph{\texttt{state\_get}}\label{ref:state:state_get}
Read one axis value from a state vector.
Signature: \texttt{state\_get(state: StateVec, axis: str) -> Any}.

\begin{lstlisting}
from codesign import make_state, state_get
x = make_state(fuel=12.0, flag=0)
print(state_get(x, "fuel"))  # 12.0
\end{lstlisting}

\begin{itemize}[nosep]
  \item \texttt{state} --- the state vector to read.
  \item \texttt{axis} --- the axis name whose value is wanted.
\end{itemize}
Errors: raises \texttt{KeyError} when the axis is absent, with the message
\texttt{"state vector has no axis 'charge'; it carries axes ['flag',
'fuel']"} --- the message names the missing axis and lists the axes the
vector actually carries.

\paragraph{\texttt{state\_as\_dict}}\label{ref:state:state_as_dict}
Return a plain \texttt{dict} view of a state vector, convenient inside a
transition or functionality callback.
Signature: \texttt{state\_as\_dict(state: StateVec) -> Dict[str, Any]}.

\begin{lstlisting}
from codesign import make_state, state_as_dict
x = make_state(fuel=12.0, flag=0)
print(state_as_dict(x))  # {'flag': 0, 'fuel': 12.0}
\end{lstlisting}

\begin{itemize}[nosep]
  \item \texttt{state} --- the state vector to expand.
\end{itemize}
Errors: none raised directly.

\paragraph{\texttt{Axis}}\label{ref:state:axis}
Abstract base class defining the single-axis contract that every concrete
axis satisfies and that \texttt{VectorStateGrid} composes into the product
order. An axis carries a \texttt{name} and the four operations below; the
two concrete axes (\texttt{ContinuousAxis}, \texttt{DiscreteAxis}) are
documented as deltas from this contract. The snippet exercises the contract
through a concrete axis, since \texttt{Axis} itself is not instantiable.

\begin{lstlisting}
from codesign import ContinuousAxis, Axis
ax = ContinuousAxis("charge", 0.0, 1.0, 5)  # a concrete Axis
print(isinstance(ax, Axis))   # True
print(list(ax.nodes()))       # [0.0, 0.25, 0.5, 0.75, 1.0]
print(ax.snap(0.3))           # 0.25 (nearest node)
print(ax.in_bounds(1.2))      # False
print(ax.leq(0.25, 0.5))      # True
\end{lstlisting}

\begin{center}
\begin{tabular}{@{}ll@{}}
\toprule
Method & Contract \\
\midrule
\texttt{nodes()}          & the axis's discretisation, as a sequence \\
\texttt{snap(value)}      & nearest admissible node to \texttt{value} \\
\texttt{in\_bounds(value)} & whether \texttt{value} lies within the axis \\
\texttt{leq(a, b)}        & the axis order: \texttt{a} at or below \texttt{b} \\
\bottomrule
\end{tabular}
\end{center}
Errors: the four base methods raise \texttt{NotImplementedError}; use a
concrete axis.

\paragraph{\texttt{ContinuousAxis}}\label{ref:state:continuousaxis}
A bucketed real interval, the vector analogue of the scalar grid and a
concrete \texttt{Axis} (\ref{ref:state:axis}). It places \texttt{n} evenly
spaced nodes on \texttt{[lo, hi]}; \texttt{snap} returns the nearest node,
\texttt{in\_bounds} tests \texttt{[lo, hi]} within a tolerance, and
\texttt{leq} applies the interval order, reversed when
\texttt{increasing\_is\_larger} is \texttt{False}.
Signature: \texttt{ContinuousAxis(name, lo, hi, n,
increasing\_is\_larger=True)}.

\begin{lstlisting}
from codesign import ContinuousAxis
ax = ContinuousAxis("slack", 0.0, 10.0, 3, increasing_is_larger=False)
print(list(ax.nodes()))  # [0.0, 5.0, 10.0]
print(ax.leq(10.0, 0.0))  # True: reversed, so less slack ranks higher
\end{lstlisting}

\begin{center}
\begin{tabular}{@{}llp{0.46\textwidth}@{}}
\toprule
Parameter & Type & Meaning \\
\midrule
\texttt{name} & str & axis name (e.g.\ \texttt{"charge"}, \texttt{"wear1"}) \\
\texttt{lo}, \texttt{hi} & float & interval bounds \\
\texttt{n} & int & number of evenly spaced nodes on \texttt{[lo, hi]} \\
\texttt{increasing\_is\_larger} & bool & \texttt{True}: larger value is
  higher in the order; \texttt{False} reverses it \\
\bottomrule
\end{tabular}
\end{center}
Errors: raises \texttt{ValueError} when \texttt{n < 1} (a message stating the
axis needs at least one node and that \texttt{n=1} yields a single node at
\texttt{lo}).

\paragraph{\texttt{DiscreteAxis}}\label{ref:state:discreteaxis}
A finite labelled axis with a caller-supplied order, a concrete
\texttt{Axis} (\ref{ref:state:axis}) for quantities like a regulatory flag
or an operating regime. Its \texttt{nodes} are the admissible labels;
\texttt{snap} is the identity on admissible labels (labels are used
verbatim); and \texttt{leq} follows \texttt{order} (least first), treating
labels omitted from \texttt{order}, or all labels when \texttt{order} is
\texttt{None}, as mutually incomparable --- which the monotonicity guard
reads as ``no order to exploit on this axis''.
Signature: \texttt{DiscreteAxis(name, values, order=None)}.

\begin{lstlisting}
from codesign import DiscreteAxis
flag = DiscreteAxis("penalty", values=[0, 1], order=[0, 1])
print(flag.nodes())      # [0, 1]
print(flag.leq(0, 1))    # True: 0 ranks below 1
print(flag.snap(1))      # 1 (snap is identity on admissible labels)
free = DiscreteAxis("regime", values=["a", "b"])  # no order given
print(free.leq("a", "b"))  # False: distinct labels are incomparable
\end{lstlisting}

\begin{center}
\begin{tabular}{@{}llp{0.46\textwidth}@{}}
\toprule
Parameter & Type & Meaning \\
\midrule
\texttt{name} & str & axis name \\
\texttt{values} & sequence & admissible labels (any hashable values) \\
\texttt{order} & sequence, optional & a chain giving the order, least
  first; omitted labels are incomparable \\
\bottomrule
\end{tabular}
\end{center}
Errors: raises \texttt{ValueError} when \texttt{values} is empty (a message
asking for the admissible labels).

\paragraph{\texttt{VectorStateGrid}}\label{ref:state:vectorstategrid}
A product grid over several named axes, carrying a full state vector. It
enumerates the Cartesian product of its axes' nodes as the DP state set,
snaps and bounds-checks a proposed successor axis by axis, and exposes the
component-wise product order used by the monotonicity results:
\texttt{a <= b} iff \texttt{a} is at or below \texttt{b} on every axis. A
single-axis grid reproduces the scalar \texttt{StateGrid}.
Signature: \texttt{VectorStateGrid(axes: Sequence[Axis])}.

\begin{lstlisting}
from codesign import VectorStateGrid, ContinuousAxis, DiscreteAxis
grid = VectorStateGrid([
    ContinuousAxis("charge", 0.0, 1.0, 3),
    DiscreteAxis("flag", values=[0, 1], order=[0, 1]),
])
print(len(grid))         # 6 state nodes (3 x 2)
print(grid.bottom())     # (('charge', 0.0), ('flag', 0))
s = grid.snap({"charge": 0.4, "flag": 1})
print(s)                 # (('charge', 0.5), ('flag', 1))
print(grid.in_bounds({"charge": 2.0, "flag": 0}))  # False
print(grid.leq(grid.bottom(), s))  # True
\end{lstlisting}

\begin{itemize}[nosep]
  \item \texttt{axes} --- one or more \texttt{Axis} objects with distinct
    names; their order fixes the product.
\end{itemize}

\begin{center}
\begin{tabular}{@{}p{0.44\textwidth}p{0.48\textwidth}@{}}
\toprule
Method & Behaviour \\
\midrule
\texttt{scalar(name, lo, hi, n, increasing\_is\_larger=True)} &
  classmethod: a one-axis continuous grid (the scalar case) \\
\texttt{nodes()}        & iterate over every state vector in the product \\
\texttt{snap(vector)}   & snap a \texttt{\{axis: value\}} mapping to the
  nearest node \\
\texttt{in\_bounds(vector)} & whether every axis value is within its axis \\
\texttt{leq(a, b)}      & the component-wise product order \\
\texttt{bottom()}       & the least state vector (each axis at its order
  minimum) \\
\texttt{len(grid)}      & number of nodes (product of axis sizes) \\
\bottomrule
\end{tabular}
\end{center}
Errors: raises \texttt{ValueError} on an empty axis sequence, and on
duplicate axis names with the message
\texttt{"VectorStateGrid axis names must be unique, but ['wear'] appear more
than once in ['wear', 'wear']. Give each axis a distinct name."}

\subsubsection{\texttt{vector\_dp}}
\label{ref:vector_dp}

The \texttt{vector\_dp} module is the general carried-state sequential
solver (section~\ref{sec:vector}): the antichain-valued Bellman recursion of
\texttt{solve\_sequential} (\ref{ref:sequential:solve_sequential}) carrying a
full state vector on a \texttt{VectorStateGrid} rather than a single scalar.
Its card structure mirrors \texttt{sequential} throughout, with the product
order of the vector grid substituted for the scalar order.

\paragraph{\texttt{VecStage}}\label{ref:vector_dp:vecstage}
One stage of a vector-state sequential co-design problem: a dataclass
bundling the stage's state-dependent functionality, its carried-state
transition, and (optionally) its candidate architectures and admissibility
predicate. The vector-state analogue of \texttt{sequential}'s
\texttt{SeqStage}.

\begin{lstlisting}
from codesign import VecStage, state_get
stage = VecStage(
    name="leg",
    functionality=lambda s: {"demand": state_get(s, "wear")},
    transition=lambda s, p: {"wear": dict(s)["wear"] + p["d_wear"]},
)
print(stage.name, stage.candidates, stage.admissible)  # leg None None
\end{lstlisting}

\begin{center}
\begin{tabular}{@{}llp{0.46\textwidth}@{}}
\toprule
Field & Type & Meaning \\
\midrule
\texttt{name} & str & stage identifier \\
\texttt{functionality} & callable & \texttt{f(state\_vec) -> Mapping}: the
  outer functionality demanded, as a function of the incoming state \\
\texttt{transition} & callable & \texttt{f(state\_vec, point) ->
  \{axis: value\}}: the outgoing state, snapped onto the grid \\
\texttt{candidates} & seq.\ of Architecture, optional & stage architectures;
  falls back to the solver default when \texttt{None} \\
\texttt{admissible} & callable, optional & \texttt{f(state\_vec) -> bool}
  forbidding states returning \texttt{False} \\
\bottomrule
\end{tabular}
\end{center}
Errors: none raised directly (a stage with no \texttt{candidates} and no
solver-level default is reported by the solver; see below).

\paragraph{\texttt{VecResult}}\label{ref:vector_dp:vecresult}
The outcome of \texttt{solve\_vector\_sequential}: a dataclass holding the
value antichain at the initial state, its width, feasibility, and the full
policy for roll-out.

\begin{lstlisting}
from codesign import (AlgebraicDP, Architecture, Ports, Reals, VecStage,
    VectorStateGrid, ContinuousAxis, solve_vector_sequential)
dp = AlgebraicDP(F=Ports({"demand": Reals()}),
    R=Ports({"cost": Reals(), "d_wear": Reals()}),
    equations={"cost": 1.0, "d_wear": 1.0})
stage = VecStage("leg", lambda s: {"demand": 1.0},
    lambda s, p: {"wear": dict(s)["wear"] + p["d_wear"]},
    candidates=[Architecture("run", dp)])
grid = VectorStateGrid([ContinuousAxis("wear", 0.0, 4.0, 5)])
res = solve_vector_sequential([stage], grid,
    cost_axes=["cost"], initial_state={"wear": 0.0})
print(res.feasible, res.width)  # True 1
print(res.value)                # Antichain[(cost=1)]
\end{lstlisting}

\begin{center}
\begin{tabular}{@{}llp{0.46\textwidth}@{}}
\toprule
Field & Type & Meaning \\
\midrule
\texttt{value} & Antichain & the front of cumulative cost-axis totals at the
  initial state, \(V_0(\text{initial})\) \\
\texttt{width} & int & number of points on that front \\
\texttt{feasible} & bool & \texttt{True} iff \texttt{value} is non-empty and
  carries no \(\topp\) \\
\texttt{policy} & VecPolicy & the full solved policy for roll-out \\
\bottomrule
\end{tabular}
\end{center}
Errors: none raised directly.

\paragraph{\texttt{VecPolicy}}\label{ref:vector_dp:vecpolicy}
A state-vector-indexed, antichain-valued policy. It is produced by
\texttt{solve\_\allowbreak vector\_\allowbreak sequential} (read from
\texttt{VecResult.policy})
rather than constructed directly, and it answers value and best-action
queries at any \texttt{(stage, state)}, snapping the queried state onto the
grid first.

\begin{lstlisting}
from codesign import (AlgebraicDP, Architecture, Ports, Reals, VecStage,
    VectorStateGrid, ContinuousAxis, solve_vector_sequential, make_state)
dp = AlgebraicDP(F=Ports({"demand": Reals()}),
    R=Ports({"cost": Reals(), "d_wear": Reals()}),
    equations={"cost": 1.0, "d_wear": 1.0})
stage = VecStage("leg", lambda s: {"demand": 1.0},
    lambda s, p: {"wear": dict(s)["wear"] + p["d_wear"]},
    candidates=[Architecture("run", dp)])
grid = VectorStateGrid([ContinuousAxis("wear", 0.0, 4.0, 5)])
pol = solve_vector_sequential([stage], grid,
    cost_axes=["cost"], initial_state={"wear": 0.0}).policy
s0 = make_state(wear=0.0)
print(pol.width_at(0, s0))                       # 1
print([p["cost"] for p in pol.value_at(0, s0)])  # [1.0]
print(pol.best_action_at(0, s0))  # ('run', {...}, (('wear', 1.0),))
\end{lstlisting}

\subparagraph{\texttt{value\_at(k, state)}} the value antichain at stage
\texttt{k} and (snapped) \texttt{state}; an empty antichain if the node is
unreached.

\subparagraph{\texttt{width\_at(k, state)}} the width (front size) of
\texttt{value\_at(k, state)}.

\subparagraph{\texttt{best\_action\_at(k, state)}} one realising
\texttt{(arch, full\_point, succ\_state)} at \texttt{(k, state)}, picking the
choice behind the cost axis's minimal-first point (a myopic roll-out);
\texttt{None} when no finite choice exists.

Errors: none raised directly.

\paragraph{\texttt{solve\_vector\_sequential}}\label{ref:vector_dp:solve_vector_sequential}
Solve the vector-state antichain-valued sequential co-design problem.
Identical in structure to
\texttt{solve\_sequential} (\ref{ref:sequential:solve_sequential}) but
carrying a full state vector on a \texttt{VectorStateGrid}; the transition
returns a \texttt{\{axis: value\}} mapping snapped onto the product grid, and
out-of-bounds successors are rejected before snapping on every axis.

\begin{lstlisting}
solve_vector_sequential(stages, grid, *, cost_axes, initial_state,
    combine=sum_combine, architectures=None, max_iter=200,
    solve_kwargs=None)
\end{lstlisting}

\begin{lstlisting}
from codesign import (AlgebraicDP, Architecture, Ports, Reals, VecStage,
    VectorStateGrid, ContinuousAxis, solve_vector_sequential)
dp = AlgebraicDP(F=Ports({"demand": Reals()}),
    R=Ports({"cost": Reals(), "d_wear": Reals()}),
    equations={"cost": 1.0, "d_wear": 1.0})
stage = VecStage("leg", lambda s: {"demand": 1.0},
    lambda s, p: {"wear": dict(s)["wear"] + p["d_wear"]},
    candidates=[Architecture("run", dp)])
grid = VectorStateGrid([ContinuousAxis("wear", 0.0, 4.0, 5)])
res = solve_vector_sequential([stage, stage], grid,
    cost_axes=["cost"], initial_state={"wear": 0.0})
print(res)                              # VecResult(feasible, width=1, ...)
print([p["cost"] for p in res.value])   # [2.0]: two legs, unit cost each
\end{lstlisting}

\begin{center}
\begin{tabular}{@{}llp{0.46\textwidth}@{}}
\toprule
Parameter & Type & Meaning \\
\midrule
\texttt{stages} & seq.\ of VecStage & horizon in forward order (solved
  backward) \\
\texttt{grid} & VectorStateGrid & product discretisation of the carried
  state \\
\texttt{cost\_axes} & seq.\ of str & resource ports accumulated on the
  antichain (the poset \(R\)) \\
\texttt{initial\_state} & mapping & \texttt{\{axis: value\}} start; the
  value returned is \(V_0(\text{initial\_state})\) \\
\texttt{combine} & callable & the monoid \((+)\): \texttt{sum\_combine}
  (accumulating) or \texttt{join\_combine} (renewable) \\
\texttt{architectures} & seq., optional & default candidates for stages
  without their own \\
\texttt{max\_iter}, \texttt{solve\_kwargs} & & forwarded to
  \texttt{solve} \\
\bottomrule
\end{tabular}
\end{center}
Errors: raises \texttt{ValueError} when a stage has neither
\texttt{candidates} nor a solver-level \texttt{architectures} default (the
message names the offending stage and both remedies).

\paragraph{\texttt{VectorMonotonicityReport}}\label{ref:vector_dp:vectormonotonicityreport}
The result of \texttt{check\_vector\_monotonicity}: a dataclass recording
whether (H1), the \emph{joint} (H2), and (H3) hold over the product order,
with sample violations.

\begin{lstlisting}
from codesign import VectorMonotonicityReport
rep = VectorMonotonicityReport(h1_ok=True, h2_ok=True,
                               h2_joint_ok=True, h3_ok=True)
print(rep.monotone_value_guaranteed)  # True: all hypotheses hold
print(rep)  # VectorMonotonicityReport(H1=ok, H2=ok, H2joint=ok, H3=ok, ...)
\end{lstlisting}

\begin{center}
\begin{tabular}{@{}llp{0.46\textwidth}@{}}
\toprule
Field & Type & Meaning \\
\midrule
\texttt{h1\_ok} & bool & (H1) holds: the stage map is consistently oriented
  in the carried state \\
\texttt{h2\_ok} & bool & state slice of (H2): the transition is monotone in
  the state \\
\texttt{h2\_joint\_ok} & bool & resource slice of (H2): the transition is
  jointly monotone on \texttt{X x R} \\
\texttt{h3\_ok} & bool & (H3): admissibility is a down-set of the grid \\
\texttt{h1\_violations} & list & up to \texttt{max\_violations} sampled
  \texttt{(stage, x, x')} witnesses \\
\texttt{h2\_violations}, \texttt{h2\_joint\_violations},
  \texttt{h3\_violations} & list & sampled witnesses for each hypothesis \\
\bottomrule
\end{tabular}
\end{center}
The property \texttt{monotone\_value\_guaranteed} is \texttt{h1\_ok and
h2\_ok and h2\_joint\_ok and h3\_ok}. Errors: none raised directly.

\paragraph{\texttt{check\_vector\_monotonicity}}\label{ref:vector_dp:check_vector_monotonicity}
Verify (H1), the \emph{joint} (H2), and (H3) over the vector grid's product
order, the vector-state analogue of the scalar guard. For every ordered
comparable pair \texttt{x <= x'} it checks the stage antichain is
consistently oriented (a benign consumable-but-monotone axis passes; a
genuinely non-monotone perishable one is flagged), that the transition is
jointly monotone on \texttt{X x R} (both its state slice and its resource
slice), and that admissibility is a down-set. The product grid can be large,
so at most \texttt{max\_pairs} comparable pairs per stage are sampled; the
report is exact when the grid is small enough that the cap is not reached.

\begin{lstlisting}
check_vector_monotonicity(stages, grid, *, cost_axes,
    architectures=None, max_iter=200, solve_kwargs=None,
    max_violations=8, max_pairs=4000)
\end{lstlisting}

\begin{lstlisting}
from codesign import (AlgebraicDP, Architecture, Ports, Reals, VecStage,
    VectorStateGrid, ContinuousAxis, check_vector_monotonicity)
dp = AlgebraicDP(F=Ports({"demand": Reals()}),
    R=Ports({"cost": Reals(), "d_wear": Reals()}),
    equations={"cost": 1.0, "d_wear": 1.0})
stage = VecStage("leg", lambda s: {"demand": 1.0},
    lambda s, p: {"wear": dict(s)["wear"] + p["d_wear"]},
    candidates=[Architecture("run", dp)])
grid = VectorStateGrid([ContinuousAxis("wear", 0.0, 4.0, 5)])
rep = check_vector_monotonicity([stage], grid, cost_axes=["cost"])
print(rep)  # H1=ok, H2=ok, H2joint=ok, H3=ok, value_monotone=...
print(rep.monotone_value_guaranteed)  # True
\end{lstlisting}

\begin{itemize}[nosep]
  \item \texttt{stages}, \texttt{grid}, \texttt{cost\_axes},
    \texttt{architectures}, \texttt{max\_iter}, \texttt{solve\_kwargs} ---
    as for the solver above.
  \item \texttt{max\_violations} --- cap on recorded witnesses per hypothesis.
  \item \texttt{max\_pairs} --- cap on comparable pairs sampled per stage;
    below it the check is exact.
\end{itemize}
Errors: raises the same \texttt{ValueError} as the solver when a stage has no
candidates and no default architecture set.

\subsubsection{\texttt{online\_codesign}}
\label{ref:online_codesign}

The \texttt{online\_codesign} module is the closed-loop, measurement-driven
counterpart of the open-loop temporal planners (section~\ref{sec:online-cd},
worked in example~\ref{ex:22}). At each control step it senses the measured
state, reads the live requirement and environment, re-solves the co-design
over the admissible architectures, applies the cheapest feasible choice by
stepping the true plant, and logs the outcome, so divergence between plan and
execution is absorbed by the next sensing step.

\paragraph{\texttt{ControlStep}}\label{ref:online_codesign:controlstep}
The audit record of one closed-loop control step: a dataclass capturing the
step's measured state, conditions, chosen architecture and point, cost, and
feasibility. \texttt{run\_online\_codesign} emits one per step.

\begin{lstlisting}
from codesign import ControlStep
step = ControlStep(step=0, measured_state=10.0, requirement={"rate": 2.0},
    environment={}, architecture="node", point={"energy": 2.0, "ops": 3.0},
    cost=3.0, feasible=True)
print(step)  # ControlStep(t=0, arch='node', cost=3, ok)
\end{lstlisting}

\begin{center}
\begin{tabular}{@{}llp{0.46\textwidth}@{}}
\toprule
Field & Type & Meaning \\
\midrule
\texttt{step} & int & step index \\
\texttt{measured\_state} & Any & the state sensed at this step \\
\texttt{requirement} & Mapping & the live functionality demanded \\
\texttt{environment} & Mapping & exogenous parameters, recorded for the log \\
\texttt{architecture} & str & the chosen architecture (empty when
  infeasible) \\
\texttt{point} & Mapping, optional & the chosen resource point
  (\texttt{None} when infeasible) \\
\texttt{cost} & float & the selected point's cost (\(\infty\) when
  infeasible) \\
\texttt{feasible} & bool & whether any architecture was feasible \\
\bottomrule
\end{tabular}
\end{center}
Errors: none raised directly.

\paragraph{\texttt{OnlineCoDesignResult}}\label{ref:online_codesign:onlinecodesignresult}
The trajectory produced by \texttt{run\_online\_codesign}: a dataclass
holding the per-step log, the total cost of feasible steps, and whether every
step was feasible. The \texttt{schedule} property extracts the sequence of
chosen architectures.

\begin{lstlisting}
from codesign import OnlineCoDesignResult, ControlStep
s = ControlStep(0, 10.0, {"rate": 2.0}, {}, "node",
    {"energy": 2.0, "ops": 3.0}, 3.0, True)
res = OnlineCoDesignResult(steps=[s], total_cost=3.0, feasible=True)
print(res.schedule)  # ['node']
print(res)           # OnlineCoDesignResult(feasible, steps=1, ...)
\end{lstlisting}

\begin{center}
\begin{tabular}{@{}llp{0.46\textwidth}@{}}
\toprule
Field & Type & Meaning \\
\midrule
\texttt{steps} & list of ControlStep & the closed-loop audit trail \\
\texttt{total\_cost} & float & summed cost over feasible steps \\
\texttt{feasible} & bool & \texttt{True} iff every step was feasible \\
\bottomrule
\end{tabular}
\end{center}
The property \texttt{schedule} returns \texttt{[s.architecture for s in
steps]}. Errors: none raised directly.

\paragraph{\texttt{resolve\_at}}\label{ref:online_codesign:resolve_at}
Solve every architecture at a fixed functionality and pick the cheapest
feasible point by \texttt{cost\_fn}. This is the single myopic co-design
solve performed once per control step; it is exposed separately so a caller
can drive the re-solve directly. Signature: \texttt{resolve\_at(architectures,
functionality, cost\_fn, *, max\_iter=200, solve\_kwargs=None)}.

\begin{lstlisting}
from codesign import AlgebraicDP, Architecture, Ports, Reals, resolve_at
dp = AlgebraicDP(F=Ports({"rate": Reals()}),
    R=Ports({"energy": Reals(), "ops": Reals()}),
    equations={"energy": lambda f: f["rate"], "ops": 3.0})
name, point, cost = resolve_at([Architecture("node", dp)],
    {"rate": 2.0}, cost_fn=lambda p: p["ops"])
print(name, cost)       # node 3.0
print(point["energy"])  # 2.0
\end{lstlisting}

\begin{itemize}[nosep]
  \item \texttt{architectures} --- the admissible configurations to solve.
  \item \texttt{functionality} --- the live demand to solve each at.
  \item \texttt{cost\_fn} --- scalar cost on a resource point; the cheapest
    feasible point wins.
  \item \texttt{max\_iter}, \texttt{solve\_kwargs} --- forwarded to
    \texttt{solve}.
\end{itemize}
Returns \texttt{(arch\_name, point, cost)}, or \texttt{("", None, inf)} when
no architecture is feasible. Errors: none raised directly.

\paragraph{\texttt{run\_online\_codesign}}\label{ref:online_codesign:run_online_codesign}
Run the closed-loop, measurement-driven co-design loop for \texttt{n\_steps}
steps. Each step senses, reads the live requirement and environment,
re-solves via \texttt{resolve\_at}, applies the cheapest feasible choice by
stepping the \texttt{plant} (the true process, not a nominal model), and logs
a \texttt{ControlStep}; an infeasible step holds the state and is recorded
without aborting the run.

\begin{lstlisting}
run_online_codesign(architectures, *, n_steps, sensor, requirement,
    plant, cost_fn, environment=None, initial_state=None,
    max_iter=200, solve_kwargs=None)
\end{lstlisting}

\begin{lstlisting}
from codesign import (AlgebraicDP, Architecture, Ports, Reals,
    run_online_codesign)
dp = AlgebraicDP(F=Ports({"rate": Reals()}),
    R=Ports({"energy": Reals(), "ops": Reals()}),
    equations={"energy": lambda f: f["rate"], "ops": 3.0})
res = run_online_codesign([Architecture("node", dp)], n_steps=3,
    sensor=lambda t, s: s,
    requirement=lambda t, s: {"rate": 2.0},
    plant=lambda t, s, a, p: s - p["energy"],
    cost_fn=lambda p: p["ops"], initial_state=10.0)
print(res)             # OnlineCoDesignResult(feasible, steps=3, ...)
print(res.schedule)    # ['node', 'node', 'node']
print(res.total_cost)  # 9.0
\end{lstlisting}

\begin{center}
\begin{tabular}{@{}llp{0.46\textwidth}@{}}
\toprule
Parameter & Type & Meaning \\
\midrule
\texttt{architectures} & seq.\ of Architecture & re-evaluated every step \\
\texttt{n\_steps} & int & number of closed-loop steps \\
\texttt{sensor} & callable & \texttt{sensor(step, prev\_state) ->
  measured\_state} \\
\texttt{requirement} & callable & \texttt{requirement(step, measured) ->
  functionality} \\
\texttt{plant} & callable & \texttt{plant(step, measured, arch, point) ->
  next\_state}; the true process \\
\texttt{cost\_fn} & callable & scalar cost selecting the point each step \\
\texttt{environment} & callable, optional & \texttt{environment(step,
  measured) -> params}; logged, empty by default \\
\texttt{initial\_state} & Any & handed to the first \texttt{sensor} call \\
\texttt{max\_iter}, \texttt{solve\_kwargs} & & forwarded to \texttt{solve} \\
\bottomrule
\end{tabular}
\end{center}
Errors: none raised directly. Fold \texttt{environment} into
\texttt{requirement} when it should influence the solve; it is otherwise only
recorded.

\section{Worked Examples}
\label{sec:examples}

The package ships 25 runnable examples under \texttt{examples/} and the
same 25 as executed notebooks under \texttt{notebooks/}. Each is
referenced below by its filename. This section walks through examples
1--17; the temporal, vector-state, and online examples 18--25 are covered
in Section~\ref{sec:temporal}.

Examples 1 to 5 are reproductions of models published by
\citet{censi2015mathematical}, and the section headings name the figure
or section of that paper each one comes from. Example~25 replicates the
synthetic benchmarks of \citet{alharbi2026online}. The remaining
examples are the author's, though several are modelled on published
case studies: the mobility and autonomous-system studies of
\citet{zardini2020avmobility,zardini2021autonomous,zardini2022taskdriven},
the perception and compute co-design of \citet{milojevic2025codei}, and
the Formula~1 season study of \citet{neumann2025formula1}.

\subsection{Example 1: the drone (Fig.~48 of \citealp{censi2015mathematical})}\label{ex:01}

\texttt{01\_drone.py}. The canonical MCDP from
\citet{censi2015mathematical}: the
``hello world'' of monotone co-design, and the smallest complete use of
the operator API (\texttt{FunctionDP} closed by \texttt{loop};
see~\ref{ref:dp:functiondp} and~\ref{ref:composition:loop}).

\subsubsection*{Model}

A battery powers an actuator that must lift the battery itself plus a
fixed payload, against a fixed avionics power draw. The functionalities
are the mission endurance (s), an \texttt{extra\_payload} (kg) and an
\texttt{extra\_power} (W); the resource of interest is the battery mass
(kg). The physics is a lumped model,
\[
  m \;\geq\; \frac{\bigl(c\,(g(m+m_{\text{pay}}))^2 + P_{\text{ex}}\bigr)\,T}
                  {\alpha},
\]
with $\alpha = 1.8\times10^6$~J/kg the cell specific energy. The battery
mass $m$ appears on \emph{both} sides: it is the loop axis. The inner
\texttt{FunctionDP} emits the required mass under two names --
\texttt{battery\_mass} (fed back by \texttt{loop}, axis
\texttt{battery\_mass}) and \texttt{report\_mass} (mirrored onto the
outer $R$ so the converged value survives on the interface rather than
collapsing to a unit poset). All ports are on $\Rplus$
(see~\ref{ref:posets:reals}).

\subsubsection*{Code}

\begin{lstlisting}
def build_drone():
    F = Ports({
        "endurance": Reals(unit="s"),
        "extra_payload": Reals(unit="kg"),
        "extra_power": Reals(unit="W"),
        "battery_mass": Reals(unit="kg"),
    })
    R = Ports({
        "battery_mass": Reals(unit="kg"),   # loop axis
        "report_mass": Reals(unit="kg"),    # mirror, stays visible
    })
    inner = FunctionDP(F, R, drone_inner_h, name="drone-inner")
    return loop(inner, axis="battery_mass", name="drone")
\end{lstlisting}

\subsubsection*{Output}

\begin{lstlisting}[language={},numbers=none,keywordstyle={},%
  extendedchars=true,literate={⊤}{{$\top$}}1]
Short, light: endurance=60s, payload=0.10kg, extra_P=1.0W
   iters=9, feasible=True, Antichain[(report_mass=0.0003564 kg)]

Medium, modest: endurance=300s, payload=0.50kg, extra_P=5.0W
   iters=22, feasible=True, Antichain[(report_mass=0.04921 kg)]

Longer mission: endurance=600s, payload=0.50kg, extra_P=5.0W
   iters=41, feasible=True, Antichain[(report_mass=0.1283 kg)]

Marginal: endurance=600s, payload=1.00kg, extra_P=10.0W
   iters=16, feasible=False, Antichain[(report_mass=⊤)]

Infeasible: endurance=1800s, payload=1.00kg, extra_P=10.0W
   iters=8, feasible=False, Antichain[(report_mass=⊤)]
\end{lstlisting}

\subsubsection*{Reading}

The example sweeps five mission profiles, three feasible and two
infeasible. The converged battery mass grows steeply with endurance --
doubling the mission from 300~s to 600~s roughly triples the mass
(0.049~kg to 0.128~kg), because heavier batteries demand more lift power,
which demands still more battery, the feedback the loop resolves. Beyond
the feasibility boundary the fixed point is $\topp$ (printed as the top
element): no finite mass closes the loop. Iteration counts (8--41) track
distance from that boundary: cases that diverge to $\topp$ or settle near
the origin resolve fastest, while the near-feasible 600~s/0.5~kg case
takes longest to climb. Examples~\ref{ex:06} and~\ref{ex:07} rebuild this
exact model through the \texttt{MCDP} and \texttt{System} builders; the
solver itself is \texttt{solve} (\ref{ref:solver:solve}).

\subsection{Example 2: integer optimisation (Sec. VI-D)}\label{ex:02}

\texttt{02\_integer\_optimization.py}. The classic Sec.\ VI-D problem from
the paper, and the cleanest illustration of \emph{multi-point} antichains
ascending under the Kleene iteration. It exercises \texttt{Naturals}
(\ref{ref:posets:naturals}), a multi-valued \texttt{FunctionDP}
(\ref{ref:dp:functiondp}) and the antichain algebra
(\ref{ref:antichains:antichain}).

\subsubsection*{Model}

Find the minimal pair $(x,y)\in\Nat\times\Nat$ with
\[
  x + y \;\geq\; \lceil\sqrt{x}\rceil + \lceil\sqrt{y}\rceil + c .
\]
Both loop variables are bundled into a single composite axis \texttt{xy}
(a nested \texttt{Ports} of two \texttt{Naturals}), so one \texttt{Loop}
closes the pair. For each candidate the inner relation enumerates every
split of the deficit $c$ across the two coordinates -- $x_{\text{out}}$
ranging over $[\lceil\sqrt{x}\rceil,\,\text{target}-\lceil\sqrt{y}\rceil]$
-- so the antichain gains several incomparable points at once. Infinities
are mapped to the top pair, giving the iteration a well-defined
$\topp$ to rest at when no finite pair works.

\subsubsection*{Code}

\begin{lstlisting}
target = sx + sy + c          # sx,sy = ceil(sqrt(x_in)),ceil(sqrt(y_in))
pts = []
for x_out in range(sx, target - sy + 1):
    y_out = target - x_out
    if y_out < sy:
        break
    pts.append({"xy": {"x": x_out, "y": y_out},
                "xy_report": {"x": x_out, "y": y_out}})
return Antichain.from_set(R, pts)
# ... inner wrapped and closed:
inner = FunctionDP(F=F, R=R, h_fn=h, name=f"sqrt_sum(c={c_value})")
return Loop(inner, axis="xy")
\end{lstlisting}

\subsubsection*{Output}

\begin{lstlisting}[language={},numbers=none,keywordstyle={}]
c = 1: iters = 6, feasible = True
   S_0: { (0, 0) }
   S_1: { (0, 1), (1, 0) }
   S_2: { (0, 2), (1, 1), (2, 0) }
   S_3: { (0, 3), (1, 2), (2, 1), (3, 0) }
   S_4: { (0, 3), (2, 2), (3, 0) }
   S_5: { (0, 3), (3, 0) }
   S_6: { (0, 3), (3, 0) }
   M(c=1) = { (0, 3), (3, 0) }

c = 4: iters = 6, feasible = True
   ...
   M(c=4) = { (0, 7), (3, 6), (4, 4), (6, 3), (7, 0) }
\end{lstlisting}

\subsubsection*{Reading}

The trace makes the \emph{ascent} concrete: the seed $S_0=\{(0,0)\}$ grows,
peaks (four points at $S_3$ for $c=1$), then \emph{contracts} as the
\texttt{Min} operation drops points that become dominated, settling at the
antichain $M(c)$ of minimal solutions. For $c=4$ it converges in six steps
to a five-point front. The example deliberately flags a typo in the paper:
the paper's $M(1)=\{(1,0),(0,1)\}$ is wrong; the solver's
$\{(0,3),(3,0)\}$ is correct, since $(1,0)$ gives
$1 \not\geq \lceil\sqrt1\rceil + 0 + 1 = 2$. That the number of points is
non-monotone across iterations is the whole point of working in
antichains rather than scalars. Example~\ref{ex:05} plots these very
iterates; the same \texttt{make\_looped} builder is imported there.

\subsection{Example 3: AUV seabed surveying (Sec. VIII)}\label{ex:03}

\texttt{03\_auv\_seabed.py}. Example~10 of the paper: a genuinely cyclic
co-design where velocity trades against sensor width, both feeding back
through mission time. It shows a \texttt{Loop} (\ref{ref:composition:loop})
over a \texttt{FunctionDP} that emits a small Pareto spread, then collapses
it with \texttt{minimize\_cost} (\ref{ref:solver:minimize_cost}).

\subsubsection*{Model}

An autonomous underwater vehicle sweeps an area $A$ (m$^2$) at velocity
$v$ (m/s) with sensor field of view $r$ (m). The couplings are
intrinsically cyclic: a larger $v$ shortens the time
$T = kA/(vr)$ but raises actuation power $P_{\text{act}}=\psi\,v^3$; a
wider $r$ also shortens $T$ but raises sensing power
$P_{\text{sens}}=\chi\,r$ and cost. Energy is $E=(P_{\text{act}}+
P_{\text{sens}})\,T$. The design pair $(v,r)$ is bundled into the loop
axis \texttt{design}; the outer resources $T$, $E$, \texttt{cost} stay
visible. At each step the inner relation enumerates a few $(v,r)$
candidates above the current loop value (up to the physical caps
$v\le3$, $r\le5$), producing the front; hitting a cap lifts the axis to
$\topp$ to signal infeasibility.

\subsubsection*{Code}

\begin{lstlisting}
pts = []
for v_try in (v, min(v*1.3, V_MAX), min(v*1.7, V_MAX)):
    for r_try in (r, min(r*1.3, R_MAX), min(r*1.7, R_MAX)):
        if v_try < v_in or r_try < r_in:
            continue
        T_try = K_GEOM * A / (v_try * r_try)
        E_try = (PSI_A*v_try**3 + CHI_A*r_try) * T_try
        pts.append({"design": {"v": v_try, "r": r_try},
                    "T": T_try, "E": E_try,
                    "cost": SENSOR_COST_A * r_try})
return Antichain.from_set(R, pts)
\end{lstlisting}

\subsubsection*{Output}

\begin{lstlisting}[language={},numbers=none,keywordstyle={}]
Area = 100 m^2
   iters=2, feasible=True
   T=1176s, E=29.6kJ, $=100
   T=905s, E=29.5kJ, $=130
   T=692s, E=29.5kJ, $=170
   best by composite cost: T=692s, E=29.5kJ, $=170

Area = 1000 m^2
   iters=2, feasible=True
   T=11765s, E=295.9kJ, $=100
   ...
   best by composite cost: T=6920s, E=295.1kJ, $=170
\end{lstlisting}

\subsubsection*{Reading}

Each scenario yields a three-point Pareto front trading survey time
against sensor cost: the fastest sweep ($T=692$~s) also costs the most
(\$170, the widest sensor), while the cheapest (\$100) is slowest. Energy
barely moves across the front -- the drag and sensing terms nearly cancel
over these $(v,r)$ choices -- so it is time-versus-cost that drives the
decision. Time scales linearly with $A$ (100~m$^2\to$10\,000~m$^2$ scales
$T$ by exactly $100\times$), confirming the coverage relation. The
composite scalar cost ($1/\text{s}+0.05/\text{kJ}+1/\$$) is dominated by
mission time, so \texttt{minimize\_cost} always selects the fastest
design. Contrast with the multi-objective handling in
Example~\ref{ex:08}, where the front is genuinely two-dimensional.

\subsection{Example 4: \texttt{UncertainDP} and \texttt{ODE\_DP}}\label{ex:04}

\texttt{04\_uncertain\_and\_ode.py}. Two short demos of the advanced
primitives: \texttt{UncertainDP} (\ref{ref:dp:uncertaindp}) for
set-bracketed uncertainty, and \texttt{ODE\_DP} (\ref{ref:dp:ode_dp}) for
deriving a monotone relation from a differential equation. Both are built
on \texttt{AlgebraicDP} (\ref{ref:dp:algebraicdp}) internally.

\subsubsection*{Model}

\emph{Uncertainty.} A battery's specific energy is known only to lie in
$[1.6,2.0]$~MJ/kg. Two \texttt{AlgebraicDP}s bracket it:
$m=\text{cap}/1.6\text{e}6$ (pessimistic, more mass) and
$m=\text{cap}/2.0\text{e}6$ (optimistic). \texttt{UncertainDP} pairs them
as \texttt{lower}/\texttt{upper}; \texttt{with\_mode} selects which bound
to solve. \emph{ODE.} A heater holds a temperature rise $\Delta T$ above
ambient; Newton cooling gives $\dot T=(P_{\text{in}}-h\Delta T)/C$, so at
steady state $P_{\text{in}}=h\,\Delta T$. \texttt{ODE\_DP} integrates the
right-hand side $\dot x = h\,\Delta T - x$ to its steady state and
extracts the root as the required power -- a monotone $R$ relation
recovered from dynamics rather than written in closed form.

\subsubsection*{Code}

\begin{lstlisting}
uncertain = UncertainDP(F=F, R=R, lower=optimistic,
                        upper=pessimistic, mode="upper",
                        name="battery_uncertain")
for mode in ("lower", "upper"):
    result = solve(uncertain.with_mode(mode), {"capacity": 3.6e6})

heater = ODE_DP(
    F=F, R=R,
    rhs=lambda x, t, f: H_LOSS * f["delta_T"] - x,
    extract=lambda x: {"power": float(x)},
    mode="steady_state", x0_fn=lambda f: 0.0, name="heater_ode")
\end{lstlisting}

\subsubsection*{Output}

\begin{lstlisting}[language={},numbers=none,keywordstyle={}]
UncertainDP demo: battery sizing under specific-energy uncertainty

   optimistic   (lower): mass = 1.800 kg
   pessimistic  (upper): mass = 2.250 kg

ODE_DP demo: power required to hold a steady temperature rise

   delta_T =    5 K  -->  P_in =   4.0 W   (=  h_loss * delta_T = 4.0 W)
   delta_T =   20 K  -->  P_in =  16.0 W   (=  h_loss * delta_T = 16.0 W)
   delta_T =   50 K  -->  P_in =  40.0 W   (=  h_loss * delta_T = 40.0 W)
\end{lstlisting}

\subsubsection*{Reading}

For 1~kWh ($3.6$~MJ) the optimistic cell needs 1.80~kg, the pessimistic
2.25~kg; the true mass lies between, and any design that survives the
2.25~kg (upper) case is robust to the chemistry actually delivered. This
is the set-based philosophy of Sec.~VII made operational: solve both
brackets, commit to the pessimistic one. The heater demo is a sanity
check -- the numerically integrated steady state reproduces the closed
form $P_{\text{in}}=h\,\Delta T$ exactly ($h=0.8$~W/K), confirming that
\texttt{ODE\_DP} recovers the algebraic relation a modeller would
otherwise hand-derive. See~\ref{sec:uncertainty} for the theory these two
primitives implement; Examples~\ref{ex:11} and~\ref{ex:12} push the
uncertainty machinery further.

\subsection{Example 5: visualising the Kleene ascent}\label{ex:05}

\texttt{05\_visualize\_kleene.py}. A visual companion to
Example~\ref{ex:02}: it renders each iterate $S_0, S_1, \ldots$ of the
Sec.\ VI-D problem as a scatter panel, reproducing the structure of
Fig.~36 of the paper. It shows how the trace recorded by \texttt{solve}
(\ref{ref:solver:solve}, via \texttt{record\_trace=True}) feeds a plotting
routine.

\subsubsection*{Model}

The MCDP is unchanged from Example~\ref{ex:02} -- the script imports
\texttt{make\_looped} from \texttt{02\_integer\_optimization} rather than
redefining it (the module name starts with a digit, so \texttt{importlib}
is used). The novelty is entirely in consuming
\texttt{result.trace}: a list of \texttt{TraceEntry} records
(\ref{ref:solver:solve}), one per Kleene iterate, each carrying the
antichain at that step. Each panel plots the finite points of one
$S_k$ in $\Nat\times\Nat$; $\topp$ points are skipped.

\subsubsection*{Code}

\begin{lstlisting}
result = solve(looped, {"c": c_value}, max_iter=50, record_trace=True)
trace = result.trace
for k, entry in enumerate(trace):
    A = entry.antichain
    xs = [p["xy"]["x"] for p in A.points if p["xy"]["x"] != inf]
    ys = [p["xy"]["y"] for p in A.points if p["xy"]["y"] != inf]
    axes[k].scatter(xs, ys, s=60, c="C3", zorder=3)
    axes[k].set_title(f"$S_{{{k}}}$  ({len(A.points)} pts)")
fig.savefig(out_path, dpi=140, bbox_inches="tight")
\end{lstlisting}

\subsubsection*{Output}

\begin{figure}[H]
  \centering
  \includegraphics[width=0.9\textwidth]{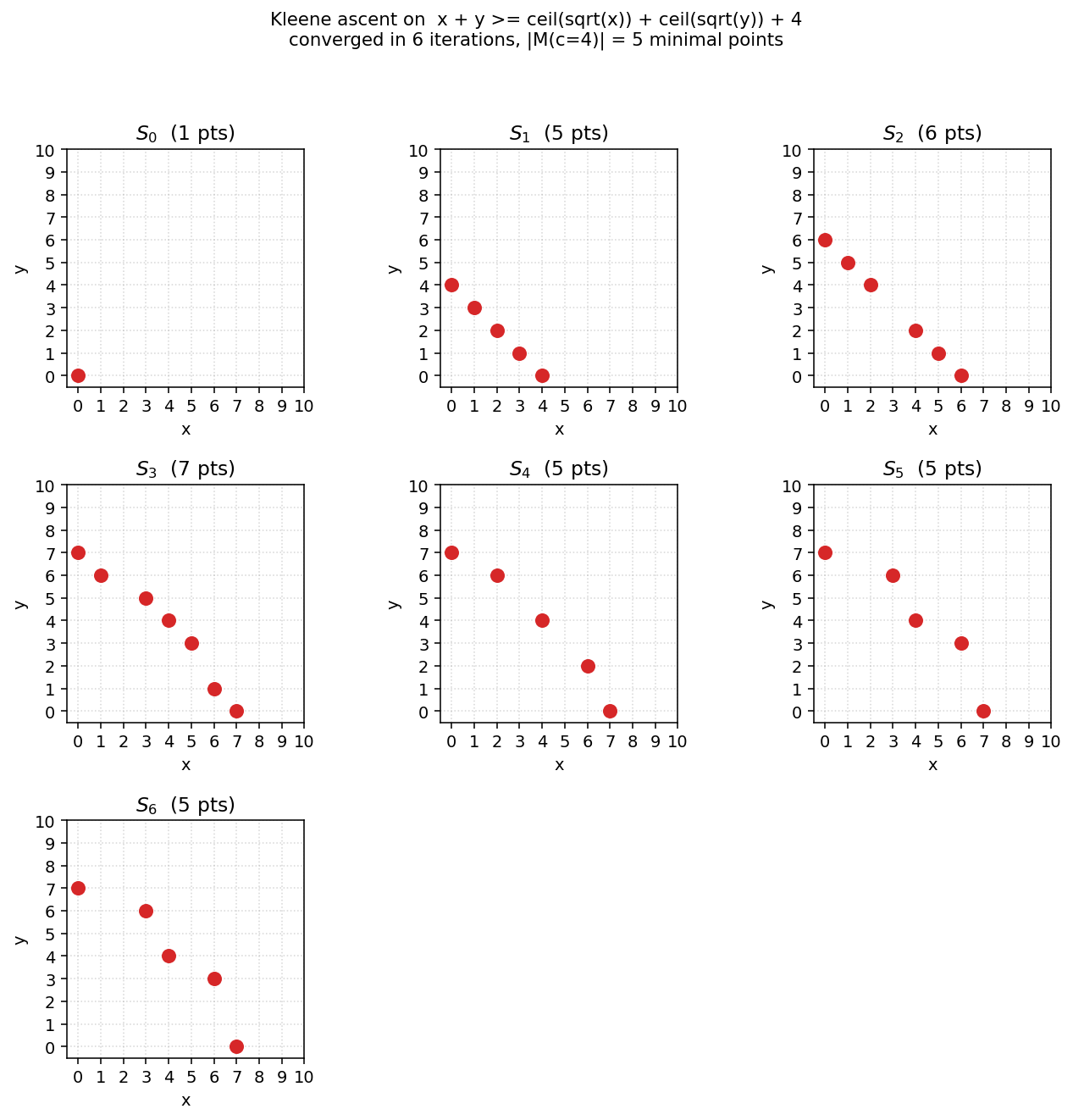}
  \caption{Kleene ascent for $c=4$: each panel is one iterate $S_k$ of the
  antichain in $\Nat\times\Nat$, from the seed $\{(0,0)\}$ to the
  five-point fixed point $M(4)$.}
  \label{fig:kleene}
\end{figure}

The script writes one PNG per value of $c\in\{1,4,8\}$ (matplotlib in the
\texttt{Agg} backend, no display required). Console output is only the
progress log:

\begin{lstlisting}[language={},numbers=none,keywordstyle={}]
Plotting Kleene iteration traces into .../outputs ...
   wrote .../outputs/kleene_trace_c1.png
   wrote .../outputs/kleene_trace_c4.png
   wrote .../outputs/kleene_trace_c8.png
\end{lstlisting}

Committed copies of the figures live under \texttt{docs/images/}
(Figure~\ref{fig:kleene}).

\subsubsection*{Reading}

The panels make the antichain algebra tangible: points appear along a
frontier $x+y=\text{const}$, the frontier advances outward, and then
dominated points fall away under \texttt{Min} until only the incomparable
minimal pairs remain. The figure is the exact geometric content of the
Sec.\ VI-D trace printed textually in Example~\ref{ex:02}; reading the two
together is the fastest way to build intuition for what a Kleene fixed
point over antichains actually looks like. This is the only example in the
first eight that touches matplotlib and therefore has a non-trivial import
cost ($\approx 1.2$~s wall-clock, essentially all of it library import).

\subsection{Example 6: the drone in MCDPL syntax}\label{ex:06}

\texttt{06\_drone\_mcdpl\_syntax.py}. The Example~\ref{ex:01} drone rebuilt
through the declarative \texttt{MCDP} builder (\ref{ref:mcdpl:mcdp}), whose
notation mirrors the paper's \texttt{mcdp\{...\}} blocks. It demonstrates
the highest-level authoring layer: \texttt{provides}/\texttt{requires}/
\texttt{constraint}/\texttt{loop\_on} inside a context manager.

\subsubsection*{Model}

The MCDP is identical to Example~\ref{ex:01}: same physics, same feedback
on \texttt{battery\_mass}. Only the \emph{spelling} changes. Inside a
\texttt{with MCDP("drone") as m:} block, functionalities are declared with
\texttt{m.provides}, resources with \texttt{m.requires}, each resource is
given its defining expression by \texttt{m.constraint}, and the recursion
is closed by \texttt{m.loop\_on("battery\_mass")}. As before,
\texttt{report\_mass} mirrors the loop value onto the visible interface.
\texttt{m.build()} compiles the block down to the same
\texttt{FunctionDP}-in-a-\texttt{Loop} the operator API produces by hand.

\subsubsection*{Code}

\begin{lstlisting}
with MCDP("drone") as m:
    m.provides("endurance", unit="s")
    m.provides("extra_payload", unit="kg")
    m.provides("extra_power", unit="W")
    m.provides("battery_mass", unit="kg")     # loop variable
    m.requires("battery_mass", unit="kg")     # loop axis
    m.requires("report_mass", unit="kg")      # mirror
    m.constraint("battery_mass", battery_mass_eq)
    m.constraint("report_mass", battery_mass_eq)
    m.loop_on("battery_mass")
return m.build()
\end{lstlisting}

\subsubsection*{Output}

\begin{lstlisting}[language={},numbers=none,keywordstyle={},%
  extendedchars=true,literate={⊤}{{$\top$}}1 {×}{{$\times$}}1]
DP(loop(drone, axis=battery_mass): endurance:R+[s]×extra_payload:R+[kg]
   ×extra_power:R+[W] -> A[report_mass:R+[kg]])

Short, light: endurance=60.0, extra_payload=0.1, extra_power=1.0
   iters=9, feasible=True, Antichain[(report_mass=0.0003564 kg)]

Medium, modest: endurance=300.0, extra_payload=0.5, extra_power=5.0
   iters=22, feasible=True, Antichain[(report_mass=0.04921 kg)]

Longer mission: endurance=600.0, extra_payload=0.5, extra_power=5.0
   iters=41, feasible=True, Antichain[(report_mass=0.1283 kg)]

Infeasible: endurance=1800.0, extra_payload=1.0, extra_power=10.0
   iters=8, feasible=False, Antichain[(report_mass=⊤)]
\end{lstlisting}

\subsubsection*{Reading}

The numbers are identical to Example~\ref{ex:01}, to the last digit
(0.0003564~kg, 0.04921~kg, 0.1283~kg) and identical iteration counts
(9, 22, 41, 8). That is the point: the \texttt{MCDP} builder is a
front-end that lowers to exactly the same internal DP as the operator API,
so the choice between them is purely one of authoring ergonomics. The
printed signature confirms the compiled object is
\texttt{loop(drone, axis=battery\_mass)} -- the same \texttt{Loop} as
before. Prefer this layer when a model is naturally a flat list of
constraint equations; drop to \texttt{FunctionDP}/\texttt{loop}
(\ref{ref:dp:functiondp}, \ref{ref:composition:loop}) when a relation
needs arbitrary Python. Example~\ref{ex:07} shows the third route, the
\texttt{System}/\texttt{Module} builder.

\subsection{Example 7: the drone, modular}\label{ex:07}

\texttt{07\_drone\_modular.py}. The Example~\ref{ex:01} drone a third way:
as a \texttt{System} (\ref{ref:system:system}) of \texttt{Module}
subclasses wired by operator-overloaded \texttt{>=} constraints. This is
the layer for genuinely modular models -- each subsystem carries its own
$F$/$R$ and relation.

\subsubsection*{Model}

Two \texttt{Module} subclasses declare class-level $F$/$R$ and an
\texttt{h}: \texttt{Battery} maps \texttt{capacity}$\to$\texttt{mass} via
specific energy, \texttt{Actuator} maps
\texttt{lift\_force}$\to$\texttt{power} as $c\,\text{lift}^2$. The
\texttt{System} exposes the outer ports, adds the two modules, and states
three inequalities that read like the textbook constraints:
battery capacity $\ge$ power$\times$endurance, actuator lift $\ge
g(m_{\text{batt}}+m_{\text{pay}})$, and total mass $\ge m_{\text{batt}}+
m_{\text{pay}}$. Each \texttt{module.port} is a handle; \texttt{>=}
records a constraint rather than evaluating a boolean.

\subsubsection*{Code}

\begin{lstlisting}
sys = System("drone")
endurance     = sys.provides("endurance", unit="s")
extra_payload = sys.provides("extra_payload", unit="kg")
extra_power   = sys.provides("extra_power", unit="W")
total_mass    = sys.requires("total_mass", unit="kg")
battery  = sys.add("battery",  Battery())
actuator = sys.add("actuator", Actuator())
battery.capacity    >= (actuator.power + extra_power) * endurance
actuator.lift_force >= G * (battery.mass + extra_payload)
total_mass          >= battery.mass + extra_payload
return sys        # sys.build() compiles to a looped DP
\end{lstlisting}

\subsubsection*{Output}

\begin{lstlisting}[language={},numbers=none,keywordstyle={},%
  extendedchars=true,literate={⊤}{{$\top$}}1]
System('drone'):
  provides: endurance, extra_payload, extra_power
  requires: total_mass
  subsystems:
    battery: (capacity) -> (mass)
    actuator: (lift_force) -> (power)
  constraints:
    battery.capacity >= ((actuator.power + extra_power) * endurance)
    actuator.lift_force >= (9.81 * (battery.mass + extra_payload))
    total_mass >= (battery.mass + extra_payload)

Short, light: endurance=60.0, ...
   iters=17, feasible=True, Antichain[(total_mass=0.1004 kg)]
Medium, modest: ...
   iters=43, feasible=True, Antichain[(total_mass=0.5492 kg)]
Longer mission: ...
   iters=81, feasible=True, Antichain[(total_mass=0.6283 kg)]
Infeasible: ...
   iters=16, feasible=False, Antichain[(total_mass=⊤)]
\end{lstlisting}

\subsubsection*{Reading}

The System reports \texttt{total\_mass} (battery + payload), so the
numbers differ from Examples~\ref{ex:01} and~\ref{ex:06} -- which reported
only the battery -- but the underlying battery masses agree exactly:
e.g.\ the medium case here is $0.5492 = 0.04921 + 0.5$~kg. Iteration
counts are roughly $2\times$ those of the single-axis model (43 vs 22, 81
vs 41) because the compiled loop now closes two coupled internal axes
(battery mass and actuator power) rather than one, so the fixed point
takes more sweeps to propagate. The declarative \texttt{>=} form and the
lambda form \texttt{sys.constrain(target, fn)} compile to the same
constraint list. This is the recommended layer for multi-subsystem
designs; Example~\ref{ex:08} scales it to three modules including a
catalog.

\subsection{Example 8: a modular vehicle with Pareto tradeoffs}\label{ex:08}

\texttt{08\_vehicle\_modular.py}. A \texttt{System}
(\ref{ref:system:system}) of three subsystems where one is a discrete
\texttt{CatalogDP} (\ref{ref:dp:catalogdp}). Because the catalog holds
Pareto-incomparable motors, the system result is a genuine multi-point
front over total mass and total cost, collapsed by
\texttt{minimize\_cost} (\ref{ref:solver:minimize_cost}).

\subsubsection*{Model}

Three modules are wired: a \texttt{CatalogDP} motor with seven entries
(each a \texttt{CatalogEntry} giving \texttt{torque} at a
(\texttt{mass}, \texttt{cost}); some share torque but trade mass against
cost, so they are incomparable), a linear \texttt{Chassis} (mass and cost
$\propto$ load) and a \texttt{Battery} (sized by energy). Five \texttt{>=}
constraints close the loop: chassis load $\ge$ payload $+$ motor $+$
battery mass, motor torque $\ge$ a load-dependent demand, battery energy
$\ge$ mission energy, and two aggregations for \texttt{total\_mass} and
\texttt{total\_cost}. The mutual dependence of chassis mass and motor
choice makes this cyclic.

\subsubsection*{Code}

\begin{lstlisting}
motor   = sys.add("motor",   make_motor_catalog())  # 7-entry CatalogDP
chassis = sys.add("chassis", Chassis())
battery = sys.add("battery", Battery())
chassis.load   >= payload + motor.mass + battery.mass
motor.torque   >= TORQUE_PER_KG * G * (payload + chassis.mass
                                       + battery.mass)
battery.energy >= mission_energy
total_mass >= payload + motor.mass + chassis.mass + battery.mass
total_cost >= motor.cost + chassis.cost + battery.cost
\end{lstlisting}

\subsubsection*{Output}

\begin{lstlisting}[language={},numbers=none,keywordstyle={}]
Small parcel: payload=2.0, mission_energy=200000.0
   iters=4, feasible=True
   Pareto front (2 points):
      total_mass=  4.82 kg,  total_cost=$ 413.00
      total_mass=  5.78 kg,  total_cost=$ 165.00
   cheapest: total_mass=5.78 kg, total_cost=$165.00

Medium load: payload=10.0, mission_energy=1000000.0
   iters=3, feasible=True
   Pareto front (2 points):
      total_mass= 22.49 kg,  total_cost=$ 475.00
      total_mass= 20.41 kg,  total_cost=$ 969.00
   cheapest: total_mass=22.49 kg, total_cost=$475.00

Heavy + long: payload=20.0, mission_energy=5000000.0
   iters=2, feasible=False
\end{lstlisting}

\subsubsection*{Reading}

Two of the three missions produce genuine two-point fronts: a lighter but
pricier motor versus a heavier but cheaper one. For the small parcel the
choice is stark -- the light design is only 0.96~kg lighter yet \$248
dearer -- so \texttt{minimize\_cost} on cost picks the 5.78~kg/\$165
option. The heavy-plus-long mission is \emph{infeasible} and the loop
detects it in just two iterations: no motor in the catalog can supply the
torque once the chassis has grown to carry a 5~MJ battery, so the fixed
point is $\topp$. This is the first-eight showcase of a discrete catalog
driving a system-level antichain; the microgrid flagship
(Example~\ref{ex:13}) and the full-vehicle study
(Example~\ref{ex:17}) build on the same \texttt{CatalogDP} pattern at much
larger scale.

\subsection{Example 9: a robotic arm with non-trivial topology}
\label{ex:09}

\texttt{09\_robotic\_arm.py}. Five subsystems wired with eight
constraints, showing where the operator-overloaded builder
(\ref{sec:system-builder}) earns its keep: the connection graph is
neither series nor parallel nor a single clean loop, and the cyclic
structure emerges from the constraint lines rather than from any
explicit loop machinery.

\subsubsection*{Model}

The arm carries two revolute \texttt{Joint} modules (shoulder, elbow),
a \texttt{Sensor}, a \texttt{Controller}, and a \texttt{Battery}. Each
joint maps functionality $(\text{torque}, \text{speed})$ to resources
$(\text{mass}, \text{electric\_power})$ with
$\text{mass} = \rho\,\tau$ and
$\text{electric\_power} = \tau\,\omega / \eta$. The topology is
genuinely coupled in three places. The controller's power draw grows
with both the sensor sample rate and the joint command rate; the
shoulder must lift the payload \emph{and} the elbow joint's own mass at
the end of its longer arm, so \texttt{shoulder.torque} reads
\texttt{elbow.mass} and closes a mechanical cycle through the elbow;
and the battery is sized by the mission-integrated power of every
electrical subsystem. The single system-level resource is
\texttt{total\_mass} over $\Rplus$ (kg).

\subsubsection*{Code}

The constraint block is the whole point: each line is one physical
relationship, and \texttt{elbow.mass} on the right of the shoulder
torque constraint is what makes the Kleene iteration iterate.

\begin{lstlisting}
elbow.torque >= G * payload_mass * elbow_arm_length
# shoulder lifts payload + the elbow joint's own mass: a cycle
shoulder.torque >= G * (payload_mass + elbow.mass) * shoulder_arm_length
sensor.sample_rate >= 2.0 * control_rate           # Nyquist
controller.input_rate   >= 2.0 * control_rate
controller.command_rate >= control_rate
battery.energy >= operating_time * (
    elbow.electric_power + shoulder.electric_power
    + controller.power + sensor.power)
total_mass >= (payload_mass + elbow.mass + shoulder.mass
    + sensor.mass + controller.mass + battery.mass)
\end{lstlisting}

\subsubsection*{Output}

Printing the assembled \texttt{System} echoes the nine constraints back
in normalised form, then three mission profiles solve to a single
minimal mass each. Every profile converges in four Kleene steps.

\begin{lstlisting}[language={},numbers=none,keywordstyle={}]
Pick-and-place light:
   iters=4, feasible=True
   total_mass = 1.97 kg
Heavier payload:
   iters=4, feasible=True
   total_mass = 7.30 kg
Long-reach precise:
   iters=4, feasible=True
   total_mass = 7.24 kg
\end{lstlisting}

No figures are produced.

\subsubsection*{Reading}

The four-step convergence for every profile shows the elbow-mass cycle
is shallow: one pass to size the elbow, one to feed its mass into the
shoulder, and the fixed point settles. The heavier-payload and
long-reach cases land at almost the same total mass by different routes
--- one is payload-dominated, the other actuator- and battery-dominated
through the higher control rate and longer operating time. This is the
cleanest small demonstration of the constraint-graph style; contrast it
with the modular vehicle of example~\ref{ex:08} and the flagship
microgrid cycle of example~\ref{ex:13}, both of which use the same
\texttt{System} builder (\ref{sec:system-builder}) at larger scale.

\subsection{Example 10: watching the solver work}
\label{ex:10}

\texttt{10\_solver\_trace.py}. The two-module drone of
example~\ref{ex:07}, instrumented to exhibit every observability
feature of \texttt{solve} in one script (see the solver observability
section, \ref{sec:solver-observability}, and the reference card
\ref{ref:solver:solve}): silent mode, the end-of-solve summary, the
per-iteration feed, the structured trace, a custom callback, and the
\texttt{status} field.

\subsubsection*{Model}

A battery ($\text{mass} = \text{capacity}/1.8\,\text{MJ/kg}$) and an
actuator ($\text{power} = 10\,F^2$) are coupled by the usual drone
cycle: capacity depends on actuator power and endurance, lift force
depends on battery mass. The model is deliberately trivial so that the
focus stays on the six calls the \texttt{\_\_main\_\_} block makes, not
on the physics.

\subsubsection*{Code}

The instrumentation is entirely in the keyword arguments to
\texttt{solve}: \texttt{verbose} controls printing, \texttt{trace}
collects a list of \texttt{TraceEntry} records
(\ref{ref:solver:traceentry}), and \texttt{on\_iteration} receives each
entry as it is produced.

\begin{lstlisting}
r = solve(drone, f, trace=True, max_iter=100)
print(f"collected {len(r.trace)} trace entries")
print(f"first 5 deltas: {[e.delta for e in r.trace[:5]]}")

def my_logger(entry):
    if entry.iteration % 5 == 0:
        print(f"iter {entry.iteration}, |A|={entry.n_points}, "
              f"delta={entry.delta}")

solve(drone, f, on_iteration=my_logger, max_iter=100)
\end{lstlisting}

\subsubsection*{Output}

The nominal solve converges in 43 iterations to the canonical
0.5492\,kg. Capping \texttt{max\_iter} short and posing an
under-resourced mission drives \texttt{status} through its three
distinct values --- these are the abbreviated forms of the final beat:

\begin{lstlisting}[language={},numbers=none,keywordstyle={}]
# 2. Final summary (verbose=1)
[solve] converged: 43 iters, |A|=1, total=0.3ms, feasible=True

# 6. Status field (max_iter vs converged vs feasible)
   max_iter=3:   status='max_iter', feasible=True,  iters=3
   max_iter=200: status='converged', feasible=True, iters=43
   infeasible:   status='diverged', feasible=False, iters=16
\end{lstlisting}

No figures are produced by the script; the companion notebook
(\texttt{10\_solver\_trace.ipynb}) adds a delta-versus-iteration plot
showing the characteristic oscillating decay of the Kleene residual to
machine precision.

\subsubsection*{Reading}

The decisive line is the last block. \texttt{status} separates two
conditions that the boolean \texttt{feasible} flag conflates: hitting
the iteration cap on a solvable problem (\texttt{'max\_iter'}, where
raising the cap would help) versus a genuinely infeasible mission
(\texttt{'diverged'}, where it would not). The trace's deltas confirm
the Kleene iteration is monotone in the fixed-point sense but not
monotone step-to-step --- the residual zig-zags downward, which is why
a single small \texttt{max\_iter} can report a non-converged but
already-feasible incumbent.

\subsection{Example 11: set-based deterministic uncertainty}
\label{ex:11}

\texttt{11\_uncertain\_drone.py}. The drone of example~\ref{ex:07} with
its battery parametrised by two internal quantities --- specific energy
and efficiency --- that are known only up to an uncertainty set. The
question posed is: under the worst-case point of that set, how heavy is
the drone? This exercises the set-based half of the uncertainty layer
(\ref{sec:uncertainty}), the \texttt{Box}
(\ref{ref:uncertainty:box}) and \texttt{Ellipsoid}
(\ref{ref:uncertainty:ellipsoid}) brackets.

\subsubsection*{Model}

The battery mass is
$\text{capacity} / (\text{specific\_energy} \cdot \text{efficiency})$.
The nominal parameters $(2.0\,\text{MJ/kg}, 0.90)$ are chosen so their
product, 1.8\,MJ/kg of \emph{delivered} energy density, reduces to the
canonical drone of examples~\ref{ex:01}, \ref{ex:06}, and~\ref{ex:07};
the nominal solve converges to the same 0.5492\,kg, and both
uncertainty sets perturb around that shared reference. Each parameter
is declared ``more is better,'' so worst-case badness is the
low-parameter direction. The \texttt{Box} admits every combination in
the two independent ranges; the \texttt{Ellipsoid} is a tilted,
correlated set that excludes the joint-extreme corner.

\subsubsection*{Code}

The uncertainty is attached to the module, not the solve; the same
system is rebuilt with a different \texttt{uncertain\_set} per query,
and \texttt{solve} is asked for the \texttt{worst\_case} summary.

\begin{lstlisting}
bat = Battery()
bat.uncertain_set = Box(
    specific_energy=(1.7e6, 2.3e6, "more_is_better"),
    efficiency=(0.83, 0.97, "more_is_better"),
)
drone = make_drone(bat)
r_box = solve(drone, f, uncertainty=["worst_case"])
wc = list(r_box.worst_case.antichain.points)[0]["total_mass"]
\end{lstlisting}

\subsubsection*{Output}

\begin{lstlisting}[language={},numbers=none,keywordstyle={}]
Nominal parameters:  total_mass = 0.5492 kg
Box uncertainty:
   worst-case total_mass = 0.5668 kg
   uncertainty penalty   = +0.0176 kg
Ellipsoid uncertainty (smaller, correlated set):
   worst-case total_mass = 0.5527 kg
   uncertainty penalty   = +0.0035 kg
\end{lstlisting}

No figures are produced.

\subsubsection*{Reading}

The two penalties, $+0.0176$\,kg for the box against $+0.0035$\,kg for
the ellipsoid, are the whole lesson. The box worst case sits at the
single corner where both parameters are simultaneously at their lowest;
the ellipsoid, being correlated and smaller, carves that implausible
corner out and returns a far milder worst case. Choosing the
uncertainty model is therefore a modelling decision about which joint
outcomes are physically credible, not a mere tolerance bound. Example
\ref{ex:12} takes the same drone stochastically and shows where the
worst case sits relative to the mean, p95, and CVaR of a Monte Carlo
distribution.

\subsection{Example 12: stochastic uncertainty with a Gaussian copula}
\label{ex:12}

\texttt{12\_stochastic\_drone.py}. The same drone battery as
example~\ref{ex:11}, but now the two internal parameters carry marginal
distributions tied by a positive correlation --- more energy-dense
cells tend also to be more efficient. This is the stochastic half of
the uncertainty layer (\ref{sec:uncertainty}): a \texttt{Stochastic}
distribution (\ref{ref:uncertainty:stochastic}) whose marginals are
glued by a \texttt{GaussianCopula}
(\ref{ref:uncertainty:gaussiancopula}).

\subsubsection*{Model}

Both parameters take a uniform marginal over the same ranges as the
box in example~\ref{ex:11}, coupled by a Gaussian copula at correlation
$0.4$. The nominal product is again 1.8\,MJ/kg delivered, so the
nominal solve is the canonical 0.5492\,kg. A single \texttt{solve} call
requests all five summaries at once --- \texttt{worst\_case} (from the
box still attached), \texttt{mean}, \texttt{p95}, \texttt{cvar95}, and
the raw \texttt{samples} --- over one thousand Monte Carlo draws at a
fixed seed.

\subsubsection*{Code}

\begin{lstlisting}
bat.uncertain_dist = Stochastic(
    marginals={
        "specific_energy": stats.uniform(loc=1.7e6, scale=0.6e6),
        "efficiency":      stats.uniform(loc=0.83, scale=0.14),
    },
    copula=GaussianCopula(correlation=[[1.0, 0.4],
                                       [0.4, 1.0]]),
)
res = solve(drone, f,
    uncertainty=["worst_case", "mean", "p95", "cvar95", "samples"],
    n_samples=1000, rng_seed=42)
\end{lstlisting}

\subsubsection*{Output}

The five summaries come out in the canonical ordering
$\text{nominal} < \text{mean} < \text{p95} < \text{cvar95} <
\text{worst\_case}$, with a feasibility rate of 1.000.

\begin{lstlisting}[language={},numbers=none,keywordstyle={}]
Nominal mass: 0.5492 kg
   mean total_mass   = 0.5506 kg
   p95 total_mass    = 0.5632 kg
   cvar95 total_mass = 0.5647 kg
   worst-case mass   = 0.5668 kg   (Box, deterministic)
Raw samples: 1000 antichains (1000 feasible)
\end{lstlisting}

The script also prints an ASCII histogram of the sampled masses; the
companion notebook renders it as a matplotlib histogram with vertical
markers at each summary. No image files are written by the script.

\subsubsection*{Reading}

The ordering is the reason to compute all five together. The mean sits
barely above nominal because the marginals are centred near the nominal
point; p95 and CVaR95 climb into the tail, and the set-based worst case
of example~\ref{ex:11} bounds them all from above. The design engineer
reads this as a spectrum of conservatism: size to the mean and half the
fleet misses the target, size to CVaR95 and only the worst 5\% of
parameter draws exceed the mass budget, size to the box worst case and
none do --- at a heavier vehicle. The copula correlation matters
because it concentrates draws along the ``both good / both bad''
diagonal, shifting the tail summaries relative to an independent model.

\subsection{Example 13: the microgrid flagship}
\label{ex:13}

\texttt{13\_microgrid.py}. An off-grid cabin must supply a daily energy
demand and a peak load without grid power. The example is the flagship
case study, exercising catalogue choice, a genuine fixed-point cycle,
warm-started sweeps, stochastic uncertainty, and the visualisation
suite in one script --- the whole \texttt{System} builder
(\ref{sec:system-builder}) and solver (\ref{ref:solver:solve}) stack
end to end.

\subsubsection*{Model}

Four subsystems contribute over $\Rplus$: a solar PV array (cheap but
sun-limited), a lithium battery whose chemistry is a discrete parameter
(LFP, NMC, LCO, NaIon, each a $(\text{Wh/kg}, \text{USD/kWh}, \text{cycle
life})$ triple), a diesel generator (reliable, carbon-heavy), and a
mounting frame. The frame is the source of the cycle: its cost and mass
scale with the total supported mass, \emph{frame included}, so the
constraint

\[
\text{frame.supported\_mass} \;\ge\;
  \text{solar.mass} + \text{battery.mass}
  + \text{diesel.mass} + \text{frame.mass}
\]

reads the frame's own output on its right-hand side and forces the
Kleene iteration to iterate to a fixed point rather than resolve in one
pass.

\subsubsection*{Code}

The sweep contrasts a cold restart against warm-starting each solve
from the previous load's fixed point via \texttt{start\_from}.

\begin{lstlisting}
warm_iters = 0
prev = None
for L in loads:                      # 50 points, 5..30 kWh
    f = {"daily_load_kwh": float(L),
         "peak_load_kw": 3.0, "backup_hours": 12.0}
    r = solve(dp, f, max_iter=400, start_from=prev)
    warm_iters += r.iterations
    prev = r                         # reuse the fixed point next step
\end{lstlisting}

\subsubsection*{Output}

At the test mission (15\,kWh/day, 3\,kW peak) all four chemistries
solve in 24 Kleene iterations; the warm-started sweep saves roughly
10\,\% of total steps.

\begin{lstlisting}[language={},numbers=none,keywordstyle={}]
   LFP     cost=$   11689  mass= 388.6kg  CO2= 416.1kg/yr  (iters=24)
   NMC     cost=$   14628  mass= 356.2kg  CO2= 416.1kg/yr  (iters=24)
   NaIon   cost=$   11328  mass= 432.8kg  CO2= 416.1kg/yr  (iters=24)

   cold-start total iters: 1205
   warm-start total iters: 1083
   speedup: 1.11x

   nominal cost:  $11689   MC mean cost: $11930
   p95 cost:      $15162   CVaR95 cost:  $16127
\end{lstlisting}

The script writes three files to \texttt{outputs/}:
\texttt{microgrid\_convergence.png} (the Kleene residual trace),
\texttt{microgrid\_uncertainty.png} (the Monte Carlo cost distribution,
x-axis labelled ``total cost (USD)''), and \texttt{microgrid.dot}
(Graphviz source for the system diagram).

\subsubsection*{Reading}

The chemistry table is a Pareto snapshot: NaIon is the cheapest but
heaviest, NMC the lightest but dearer, LFP the balanced default, and
CO$_2$ is fixed because it is set by the diesel backup, not the
chemistry. The 24-iteration count confirms the frame cycle is real ---
a acyclic model would resolve in one or two passes. The 1.11$\times$
warm-start speedup is modest but free: consecutive loads have nearby
fixed points, so seeding from the previous solution shaves off the
early Kleene steps. The stochastic pass shows the design is robust ---
feasibility stays at 1.000 as sun hours vary --- while the tail cost
climbs by roughly 40\,\% from nominal to CVaR95, the premium for
sizing against bad-sun years. Compare the discrete-catalogue mechanics
here with the online elimination of example~\ref{ex:14}, which prunes a
catalogue instead of enumerating it.

\subsection{Example 14: online elimination over a robot catalogue}
\label{ex:14}

\texttt{14\_online\_fleet.py}. A logistics service must hit a target
throughput and range by buying robots from a catalogue of 200 candidate
types, each described by four features (speed, payload, unit cost,
energy per km). Naively this is 200 inner solves; the online solver
(\ref{sec:online-cd}, \ref{ref:online:solve_online}) prunes provably
suboptimal candidates using an optimistic evaluator and runs inner
solves only for the survivors. The example reproduces the spirit of the
multi-robot fleet study of \citet{alharbi2026online}.

\subsubsection*{Model}

Each robot's inner DP is a smooth (fractional-fleet)
\texttt{AlgebraicDP} mapping the mission
$(\text{throughput}, \text{range})$ to
$(\text{total\_cost}, \text{total\_energy})$, with
$\text{total\_cost} = (\text{throughput}/(sp))\,c$ and
$\text{total\_energy} = \text{range}\cdot e \cdot 24$. Three evaluators
are compared: a \texttt{LipschitzEvaluator}
(\ref{ref:online:monotonicityevaluator} and neighbours) with per-output
constants; a \texttt{MonotonicityEvaluator} on the derived feature
$\text{cost\_per\_capacity} = c/(sp)$, under which \texttt{total\_cost}
is exactly monotone; and the certified confidence-polytope
\texttt{LinearParametricEvaluator}.

\subsubsection*{Code}

The three flavours differ only in the evaluator object handed to
\texttt{solve\_online}; the inner-DP factory and mission are shared.

\begin{lstlisting}
res = solve_online(make_dp, mission,
                   candidates=candidates,
                   evaluator=MonotonicityEvaluator(
                       features=["cost_per_capacity", "energy_per_km"],
                       r_components=["total_cost", "total_energy"]),
                   verbose=0)
print(res.n_evaluated, res.n_eliminated, len(res.antichain))
\end{lstlisting}

\subsubsection*{Output}

The exhaustive baseline finds a 5-point Pareto front. All three
evaluators recover every Pareto point; they differ only in how much of
the catalogue they had to evaluate.

\begin{lstlisting}[language={},numbers=none,keywordstyle={}]
True Pareto front: 5 non-dominated candidates

Lipschitz (L=300/30):        evaluated 192/200  Pareto recovery: OK
Monotonicity (cost_per_cap): evaluated  13/200  Pareto recovery: OK
LinearParametric (certified):evaluated 122/200  Pareto recovery: OK
\end{lstlisting}

The script writes \texttt{fleet\_online.png} (three feature-plane
panels shading each candidate by status) and
\texttt{fleet\_online\_convergence.png} (incumbent antichain size
against inner solves performed) to \texttt{outputs/}.

\subsubsection*{Reading}

The evaluated counts are the story. Monotonicity with the right derived
feature prunes hardest --- 13 of 200 --- because a single monotone
coordinate orders the whole catalogue. Lipschitz is the most general
and most conservative, evaluating 192 of 200: with a safe $L$ it never
drops a Pareto point but rarely eliminates one either. The certified
\texttt{LinearParametric} evaluator sits between, at 122 of 200, and
recovers all five: it maintains the polytope of linear parameter
vectors consistent with every observation and lower-bounds each query
by one LP over that polytope, so the bound is guaranteed and never
wrongly eliminates an optimal candidate. It prunes less aggressively
here because \texttt{total\_cost} is not affine in the raw features; on
a truly linear map it becomes exact and aggressive. This certified
bound replaces an earlier OLS $\pm$ confidence-band heuristic that
pruned to about 35 of 200 but wrongly dropped one Pareto point on this
seed --- the trade is aggressiveness for a guarantee. Example
\ref{ex:16} pushes the same three evaluators onto a deliberately
nonlinear surface and shows the certified evaluator declining to guess.

\subsection{Example 15: monoclonal antibody fed-batch co-design}
\label{ex:15}

\texttt{15\_bioprocess.py}. A biopharmaceutical company must deliver a
monoclonal antibody at a target titer (g/L) and annual demand
(kg/year). Four subsystems are coupled cyclically, and two of them are
catalogue lookups (\ref{ref:dp:catalogdp}) resolved by the solver. This
is the first biotech-upstream application in the library, calibrated to
the 2024--2026 bioprocessing literature
\citep{yoon2003low,sou2015hypothermia,trummer2006shifting,khattak2010feed,%
gagnon2011hipdog,reinhart2019bioprocessing,lao1997ammonium}, and it uses the same
\texttt{System} builder (\ref{sec:system-builder}) as the engineering
examples.

\subsubsection*{Model}

A \texttt{CellLine} module maps demanded titer to the required peak
viable-cell density and oxygen uptake via an effective integrated
specific productivity $q_P$. A \texttt{FeedStrategy} module turns the
glucose set-point into a metabolic-burden factor with a U-shaped
batch-failure penalty around 8\,mM. Two \texttt{CatalogDP} lookups pick
the smallest sufficient bioreactor (single-use 200\,L through stainless
25{,}000\,L) and the cheapest sufficient media (four commercial CD
formulations). The cycle the Kleene iteration resolves runs through the
peak cell density: higher titer needs higher density, which inflates
the metabolic burden through the feed, which inflates the density
needed for the same titer. The outer objective is the vector
$(\text{COGS}, \text{footprint}, \text{CO}_2\text{ per gram})$.

\subsubsection*{Code}

Cell line and feed are fixed at construction so the outer loop sweeps
them; the bioreactor and media are chosen by the solver. The three
outer resources are assembled from ports plus closed-over parameters
via \texttt{sys.constrain}.

\begin{lstlisting}
cell  = sys.add("cell",  CellLine(*cell_line))
feed  = sys.add("feed",  FeedStrategy(glucose_setpoint_mm))
bior  = sys.add("bior",  make_bioreactor_dp())   # CatalogDP
media = sys.add("media", make_media_dp())        # CatalogDP
# metabolically-inflated density must be supported by both catalogues
bior.peak_vcd  >= cell.peak_vcd * feed.metabolic_factor
media.peak_vcd >= cell.peak_vcd * feed.metabolic_factor
sys.constrain("cogs_per_g", cogs_eq)   # ports + annual_demand, yield
\end{lstlisting}

\subsubsection*{Output}

At each of the three mission scales the global Pareto front is a genuine
two-point tradeoff between the high-producer workhorse CHO-K1 and the
next-generation short-batch CHO-MK, both at the 8\,mM feed optimum.

\begin{lstlisting}[language={},numbers=none,keywordstyle={}]
=== titer 5.0 g/L, annual demand 100 kg/year ===
  design            COGS      footprint      CO2
  CHO-K1/glu=8mM   $ 44.04/g   39.9 m^2    214.29 g/g
  CHO-MK/glu=8mM   $ 63.65/g   30.5 m^2    214.29 g/g
  Cheapest design overall: CHO-K1/glu=8mM
\end{lstlisting}

No figures are produced.

\subsubsection*{Reading}

The two Pareto points are not comparable: CHO-K1 wins on cost of goods
(low licence fee, long batch) while CHO-MK wins on footprint (high
licence fee, short batch, so fewer parallel lines at the same annual
output). CHO-K1 is the cheapest overall at every scale, and per-gram
COGS, footprint intensity, and CO$_2$ all fall as demand rises, the
signature of amortising catalogue capex over more grams. The 8\,mM feed
wins at both cell lines because the U-shaped penalty punishes both the
starvation and the waste-accumulation extremes. All parameters sit
within published ranges (specific productivity from Reinhart \etal\
2019, metabolic limits from Khattak \etal\ 2010 and Lao \& Toth 1997);
the values are illustrative, not production-precision. Example
\ref{ex:16} fixes CHO-K1 and turns to the operating-condition
optimisation this model leaves open.

\subsection{Example 16: online DOE for the mAb fed-batch process}
\label{ex:16}

\texttt{16\_online\_doe.py}. The bioprocess of example~\ref{ex:15} is
fixed (CHO-K1, 100\,kg/year, 5\,g/L) and the question shifts to where to
operate within a $5\times5\times5\times3 = 375$-point grid of operating
conditions --- temperature, pH, glucose target, feed-start day. Each
candidate stands in for a 10--14-day bioreactor run costing
\$20{,}000--\$100{,}000, so the online elimination of
example~\ref{ex:14} (\ref{sec:online-cd},
\ref{ref:online:solve_online}) is now a Design-of-Experiments budget
question.

\subsubsection*{Model}

A closed-form effect model maps each condition vector to an effective
productivity and thence to the outer triple
$(\text{COGS}, \text{footprint}, \text{CO}_2)$, calibrated from the
process-parameter literature (temperature downshift from Yoon \etal\
2003 and Sou \etal\ 2015, pH from Trummer \etal\ 2006, glucose and feed
from Khattak \etal\ 2010 and Gagnon \etal\ 2011). The surface is
deliberately nonlinear --- a temperature/pH/glucose U-shape and a
power-law batch-failure rate. Two non-online baselines (a 75-run
factorial DOE at the pH\,=\,7.1 slice, and a 40-run uniform random
sample) are compared against the Lipschitz, Monotonicity, and certified
LinearParametric evaluators at a budget of 40 inner solves, scored by
how many of the four distinct Pareto \emph{classes} each recovers.

\subsubsection*{Code}

Recovery is counted by output value, not candidate key, because many
grid points share the same $(\text{COGS}, \text{footprint})$ outcome.

\begin{lstlisting}
true_classes = {(round(p["cogs_per_g"], 2),
                 round(p["footprint_m2"], 1)) for p in true_pareto}
for name, ev in evaluators:
    res = solve_online(make_dp, {"target_titer": 5.0},
                       candidates=candidates, evaluator=ev, budget=40)
    recovered = {(round(o["cogs_per_g"], 2),
                  round(o["footprint_m2"], 1))
                 for i in res.incumbent_ids for o in [simulate_run(...)]}
    print(name, len(true_classes & recovered), "/", len(true_classes))
\end{lstlisting}

\subsubsection*{Output}

The exhaustive baseline (375 solves) yields four Pareto classes. At
budget 40 the tuned Lipschitz evaluator matches the 75-run factorial
DOE at recovering 3 of 4; the certified LinearParametric evaluator
recovers zero, and correctly so.

\begin{lstlisting}[language={},numbers=none,keywordstyle={}]
  strategy                         evals   elim   Pareto recovered
  Factorial DOE (pH=7.1 slice)        75     -    3 / 4 classes (75%)
  Random sample (seed=42)             40     -    3 / 4 classes (75%)
  Lipschitz (normalised features)     40     0    3 / 4 classes (75%)
  Monotonicity (cogs only)            40     0    0 / 4 classes ( 0%)
  LinearParametric                    40     0    0 / 4 classes ( 0%)

  Monotonicity + 4 corner warm-start runs: recovery 1/4 classes (25%)
\end{lstlisting}

The companion notebook visualises the pick trajectory in the
(temperature, glucose) projection.

\subsubsection*{Reading}

The zero from LinearParametric is the safe-degradation property, not a
failure. On the linear catalogue of example~\ref{ex:14} the same
evaluator recovered every Pareto point; here the map is markedly
nonlinear, no single affine parameter set fits the observations, the
confidence polytope cannot support a non-trivial bound, and the
evaluator falls back to the no-information value rather than risk
eliminating a candidate it cannot certify as suboptimal. Monotonicity
alone is likewise silent because its lower bounds tighten only for
candidates above every observation, and the Pareto front lives at the
low-feature corner; seeding four hand-picked corner runs lifts it from
0 to 1 class. The generalisable lesson is that the structural prior
does the work and must match the problem: a certified global prior is
exact on an affine map and deliberately silent off it, a Gaussian
process buys reach on smooth nonlinear surfaces at the price of an
uncertified bound, and pure local bounds are always safe but need
denser observations or a warm-start to concentrate.

\subsection{Example 17: full-vehicle co-design across ICE, hybrid, and EV}
\label{ex:17}

\texttt{17\_car\_codesign.py}. The largest example in the package
(3351 lines) and the first to model three architecturally distinct
variants in one design study. A passenger car is decomposed into
per-architecture MCDP modules --- 24 for ICE, 23 for hybrid, 18 for EV
--- each with its own F and R ports, assembled with the \texttt{System}
builder (\ref{sec:system-builder}) and the catalogue DP
(\ref{ref:dp:catalogdp}) over engines, transmissions, motors, and
battery packs. Because the file is so large, only the EV mass/energy
loop and the sweep skeleton are quoted here; the three
\texttt{build\_*\_car} bodies follow the same template. Component masses,
efficiencies, and costs are calibrated from the automotive engineering
literature
\citep{genta2009chassis,bosch2018handbook,pulkrabek2003engine,heywood2018ice,%
hofmann2014hybrid,naunheimer2011transmissions,iea2023evoutlook,nemry2008cars,%
epa2024trends,larminie2012ev,mock2023pocketbook}.

\subsubsection*{Model}

ICE, hybrid, and EV share a common chassis core (body, front/rear
suspension, front/rear brakes, steering, tires, wheels) and auxiliary
core (HVAC, interior, safety, lighting/infotainment); they differ in
the powertrain modules. Two coupled cycles close inside every
architecture and are resolved by the Kleene iteration. The first is the
\emph{mass spiral}: every load-bearing module reads the design curb
weight as an F input, and the macro aggregation sums every module
weight back into \texttt{curb\_weight}, so heavier vehicle drives
heavier components drives heavier vehicle. The second is the
\emph{energy-storage loop}: fuel-tank size (ICE/HEV) or battery
capacity (EV/HEV) is sized by consumption $\times$ range, and the
storage's own mass re-enters the curb weight. For EVs the pack is
20--35\,\% of curb weight, so this loop dominates and convergence takes
30--50 iterations against 15--25 for ICE. The macro objective is the
vector $(\text{cost}, \text{curb weight}, \text{CO}_2)$, with fuel or
energy consumption, maintenance, and durability also reported.

\subsubsection*{Code}

Inside \texttt{build\_ev\_car}, one weight expression is fed to every
load-bearing module (the mass spiral) and the battery is sized from the
mission range and a consumption lambda (the energy loop).

\begin{lstlisting}
total_weight_expr = (body_m.weight + fr_m.weight + rr_m.weight
    + bf_m.weight + br_m.weight + st_m.weight + ti_m.weight
    + wh_m.weight + mot_m.weight + bat_m.weight + pe_m.weight
    + obc_m.weight + red_m.weight + btm_m.weight + DRIVER_MASS_KG)
for m in (fr_m, rr_m, bf_m, br_m, st_m, ti_m, wh_m, li_m):
    m.design_mass >= total_weight_expr          # mass spiral
sys.constrain("hv_battery.target_range_km",
              lambda x, r=mission["target_range_km"]: r)
sys.constrain("hv_battery.energy_kWh_per_100km",
              ev_energy_kWh_per_100km)           # closes energy loop
\end{lstlisting}

Each architecture has a sweep function enumerating its catalogue
combinations; \texttt{main} runs all three across four missions and
keeps the cheapest feasible design per pair.

\begin{lstlisting}
def sweep_ev(mission, *, verbose=False):
    results = []
    for body in _eligible_bodies(mission):
        for motor in EV_MOTORS:
            for battery in EV_BATTERIES:
                dp = build_ev_car(mission=mission, body=body,
                                  motor=motor, battery=battery, ...)
                res = solve(dp, dict(mission), max_iter=200)
                if res.feasible and res.antichain.points:
                    results.append((label, dict(list(
                        res.antichain.points)[0])))
    return results
\end{lstlisting}

\subsubsection*{Output}

The script prints a per-mission table and a 10-year-TCO summary. The
trimmed summary below reproduces the headline claims. For the Family
Daily mission (5 passengers, 700\,km range, 9\,s 0--100): the cheapest
ICE is \$35{,}523 at 196\,g/km, the cheapest HEV \$40{,}272 at 91\,g/km
with the best 10-year TCO at \$56{,}427, and the cheapest EV \$64{,}666
at 58\,g/km but at 2434\,kg curb weight. The 7-passenger 800\,km EV
(Suburban Utility) has zero feasible designs.

\begin{lstlisting}[language={},numbers=none,keywordstyle={}]
Mission             Arch    Cost   Weight     CO2   10y TCO
Family Daily        ICE   $35,523   1710kg  196g/km $64,998
Family Daily        HEV   $40,272   1719kg   91g/km $56,427
Family Daily        EV    $64,666   2434kg   58g/km $76,052
Suburban Utility    HEV   $42,510   1910kg  113g/km $61,066
Suburban Utility    EV        -        -       -       -   (0 feasible)
\end{lstlisting}

The full module graphs are shipped as
\texttt{docs/diagrams/car\_\{ice,hev,ev\}.svg}; the \texttt{draw\_system}
renderer shades the mass-spiral strongly-connected component in amber.

\subsubsection*{Reading}

The three architectures trade cost against CO$_2$ monotonically ---
ICE cheapest to build and dirtiest, EV dearest and cleanest, HEV in
between --- but the 10-year TCO reorders them: the HEV's low fuel bill
gives it the best total cost of ownership on the Family Daily mission,
undercutting both the cheaper-to-buy ICE and the far dearer EV. The
2434\,kg EV curb weight is the mass spiral and energy loop compounding:
a 100\,kWh 800V pack is needed for 700\,km, and the pack's own mass
inflates every load-bearing component. The Suburban Utility EV
infeasibility --- 7 passengers over 800\,km --- is the model correctly
reproducing the real-world absence of such a product against 2024
pack-level energy density (160--180\,Wh/kg): the required pack mass
never reaches a fixed point below the chassis limit. Calibration cites
Genta \& Morello, the Bosch handbook, Heywood, Hofmann, and the IEA EV
Outlook among others; the numbers are ``within published ranges for
2020--2024 production vehicles,'' and the contribution is the
compositional analysis, not the point values. This example is the
capstone of the static-solve examples; the temporal, sequential, and
online-co-design regimes begin with the following section
(\ref{sec:temporal}).

\section{Temporal, Vector-State, and Online Co-Design}
\label{sec:temporal}

The examples up to 17 solve a co-design problem \emph{once}: a fixed
functionality is posed and the Kleene iteration returns the minimal
resource antichain. Many engineering systems are not static. The best
architecture changes as the environment changes; a resource is spent or
worn across a sequence of stages; the design must be re-decided online as
measurements arrive. This section documents the six layers the package
adds on top of the ordinary \texttt{solve()} to handle these regimes, and
the theory that makes them exact. All six reuse a common decision object,
the \texttt{Architecture} (a named wrapper around a \texttt{System}), so
they compose with the whole rest of the framework.

Unlike the rest of the manual, this section does not document published
theory. The static co-design machinery underneath it is Censi's
\citep{censi2015mathematical}; the temporal, sequential, vector-state, and
receding-horizon layers built on top of it are the author's own and are
being prepared for separate publication. They have not been peer reviewed.
The closest published work is the hierarchical seasonal co-design of
\citet{neumann2025formula1}, which was developed independently.

\subsection{Case 1: architecture switching (\texttt{temporal})}

The simplest temporal problem is \emph{switching}: a system whose best
architecture changes because the environment changes, with a cost to
switch. The environment is a sequence of epochs, each posing its own
functionality and admitting a subset of candidate architectures.
\texttt{solve\_schedule} chooses one architecture per epoch to minimise
the total of the per-epoch co-design costs plus the switching costs, by an
exact Viterbi (shortest-path) pass over the epoch-by-architecture lattice.
The switching cost may be a uniform scalar or a per-transition callable,
and a hysteresis parameter lets a brief unfavourable epoch be ridden out
on the incumbent rather than switched into and back out of. The per-epoch
cost is a full \texttt{solve()}, so switching sits on top of static
co-design rather than replacing it.

\subsection{Case 2, scalar value: dynamic programming (\texttt{dynamic})}

When a resource is carried \emph{between} stages, the problem becomes a
dynamic program. \texttt{solve\_dynamic} runs a finite-horizon backward
Bellman recursion whose per-stage decision is which architecture to
instantiate, whose per-stage cost is itself a co-design \texttt{solve()},
and whose carried scalar state (fuel, charge, cumulative wear) is advanced
by a user \texttt{transition}. The value function is scalar: at each
(stage, state) it is the single best cost-to-go. The state is discretised
onto an explicit \texttt{StateGrid}; a transition that would leave the
grid envelope is rejected \emph{before} snapping, so an over-spent
resource (a negative fuel level) is never silently rescued by the nearest
in-range node. The solver returns a state-indexed \texttt{DynamicPolicy},
queryable in closed loop at off-nominal states, and \texttt{rollout}
threads it forward from a concrete initial state.

\subsection{Case 2, antichain value: sequential co-design (\texttt{sequential})}
\label{sec:sequential}

The multi-objective generalisation replaces the scalar value with a Pareto
antichain. Following the sequential co-design theory, fix an ordered
commutative resource monoid $(R, \le, \oplus, 0)$ and work in the value
space of upper sets of $R$ ordered by reverse inclusion (a larger feasible
set is a lower resource threshold, hence ``easier''). A sequential
co-design problem over stages $0,\dots,N$ and a state poset
$(X, \le_X, \bot_X)$ assigns to each stage $k$ a state-parametrised
minimal-resource map $h_k : X \to \mathcal{U}R$ and a transition
$\phi_k : X \times R \to X$. The value is the backward recursion
$V_{N+1} \equiv R$ and
\[
  V_k(x) \;=\; \mathrm{Min} \bigcup_{r \in h_k(x)}
      \Bigl( \uparrow\! r \;\oplus\; V_{k+1}\bigl(\phi_k(x, r)\bigr) \Bigr),
\]
the antichain-valued Bellman operator, implemented by
\texttt{Antichain.union\_min}. In the package the accumulated resource
(the named cost axes on the antichain) is kept deliberately \emph{distinct}
from the carried state $x$: the transition reads any carried quantity off
the full solved point, while only the cost axes accumulate on the front.
The scalar \texttt{dynamic} layer is the width-one special case (a single
cost axis gives a singleton antichain).

Three results are made operational. \textbf{(Q1) Monotone value.} If, for
every stage, (H1) $h_k$ is monotone in the carried state (more carried
state makes the stage no easier), (H2) $\phi_k$ is \emph{jointly} monotone
on $X\times R$ (monotone in the state \emph{and} in the chosen resource
point, so the successor built from a smaller state and a smaller resource
stays comparable), and (H3) admissibility is a \emph{down-set} of $X$ (a
state below an admissible one is admissible, so the out-of-bounds guard
never deletes the successor the proof constructs), then every $V_k$ is
monotone. (H1) is load-bearing and necessary; the joint slice of (H2) and
(H3) are the two obligations the informal statement previously omitted, and
both hold for the canonical consumable-sum and renewable-join transitions.
\texttt{check\_monotonicity} verifies all three on the state grid; it is
orientation-aware, accepting the benign ``consumable but monotone''
orientation (a budget carried as remaining slack) and flagging only
genuinely non-monotone (perishable or fatigue-as-state) stages, a transition
that is not jointly monotone, or a non-down-set admissible region, all of
which fall outside the guarantee.
\textbf{(Q2) Front $=$ reachable frontier.} $\mathrm{Min}\,V_k(x)$ equals
the antichain of minimal reachable cumulative resources over feasible
choice sequences; there is no tail-pruning gap, and the front size is the
width of the reachable set, which grows polynomially in the horizon for a
summed objective on a fixed number of axes (\texttt{sum\_combine}) and
stays bounded for a renewable one (\texttt{join\_combine}).
\textbf{(Q3) Exact factorisation.} At a reset stage (a quiescent
transition landing on $\bot_X$ regardless of input),
\texttt{detect\_resets} and \texttt{factorise\_at\_resets} split the
horizon into a $\oplus$-product of independent sub-problems, the
order-theoretic deterministic analogue of regeneration-point
decomposition; the proof uses only distributivity, so it holds in every
regime.

\subsection{The general carried state: vectors (\texttt{state}, \texttt{vector\_dp})}
\label{sec:vector}

Realistic problems carry structured state, not a lone scalar. The Formula 1
seasonal co-design of Neumann, Zardini, and colleagues
\citep{neumann2025formula1} carries a
state vector of two battery-wear levels plus a discrete regulatory-penalty
flag; a self-reconfiguring robot carries the accumulated wear of each drive
module plus a shared energy budget. The \texttt{state} module provides a
\texttt{VectorStateGrid}: a product grid over named axes, each either a
\texttt{ContinuousAxis} (a bucketed real interval, snapped and
bounds-checked like the scalar grid) or a \texttt{DiscreteAxis} (a finite
labelled set with a caller-supplied order, or mutually incomparable when no
order is given). \texttt{solve\_vector\_sequential} runs the same
antichain-valued Bellman recursion carrying the full state vector, with the
transition returning an $\{\text{axis} : \text{value}\}$ mapping snapped
onto the product grid and out-of-bounds successors rejected on every axis.
The monotonicity results carry over verbatim with the component-wise
product order substituted for the scalar order, and
\texttt{check\_vector\_monotonicity} verifies (H1), the joint (H2), and (H3)
over that order. A
single-axis grid reproduces the scalar \texttt{sequential} result exactly,
so the vector layer strictly generalises rather than parallels it.

A subtle correctness point governs the backward pass. The stage antichain
must not be Pareto-reduced on the cost axes \emph{before} the transition is
applied: a point dominated on cost may be the only feasible choice from a
constrained carried state (a higher-cost morphology that spares a worn
module), and its consequence for the carried state is invisible in the
cost projection. The implementation therefore enumerates all solved points
and lets \texttt{union\_min} prune only after the transition. For the same
reason the pre-transition deduplication keys candidate points on their
\emph{full} resource point (all ports), not on the cost projection: two
points that share a cost but differ on a carried axis (equal cost, different
fuel) are incomparable in the full resource poset and must both survive, so
the $\mathrm{Min}$ lives in $R$ rather than its projection. This is what
makes the reconfigurable-robot example (21) feasible.

\subsection{Precompute-then-DP: the Formula 1 structure}
\label{sec:precompute}

There are two distinct ways to combine co-design with dynamic programming,
and they should not be conflated. The \texttt{sequential} solver above
re-solves the co-design at every (stage, state), which is necessary when
the per-stage functionality depends on the carried state. The Formula 1
seasonal framework \citep{neumann2025formula1} uses the cheaper opposite order: the co-design
layer is run \emph{once} to produce track-dependent Pareto-optimal
performance-versus-wear mappings, which are then frozen into a catalog of
implementations that an outer finite-horizon dynamic program selects
among. No co-design solve happens inside the Bellman sweep; the DP only
indexes the precomputed catalog. The package exposes this explicitly:
\texttt{precompute\_catalog} solves the architectures once at a fixed
functionality and returns the tagged Pareto catalog (retaining points that
look dominated at the single-stage level, since they may become optimal
once aggregated in the DP, exactly the F1 paper's observation), and
\texttt{dp\_over\_catalog} runs the antichain-valued DP selecting catalog
indices. Precompute-then-DP agrees exactly with \texttt{solve\_sequential}
when the co-design is state-independent, and is preferred there for cost;
when the co-design genuinely depends on the carried state, only the
re-solving \texttt{sequential} form is correct.

\subsection{Online feedback co-design (\texttt{online\_codesign})}
\label{sec:online-cd}

The layers above are open-loop planners: the whole horizon is solved in
advance. The closed-loop counterpart re-decides the design online as
measurements arrive. \texttt{run\_online\_codesign} implements a
receding-horizon, measurement-in-the-loop controller: at each step it
(i) senses the measured state of the running process, (ii) reads the live
requirement and environmental conditions, (iii) re-solves the co-design at
those conditions over the admissible architectures (\texttt{resolve\_at}
picks the cheapest feasible point), (iv) applies the chosen configuration
by stepping the \emph{true} plant, and (v) logs a \texttt{ControlStep}.
Because the plant is stepped with the true process rather than a nominal
transition, any divergence between plan and execution is absorbed by the
next sensing step, which is what distinguishes feedback from open-loop
replay. This is the co-design instance of control co-design
(CCD, \citealp{garciasanz2019ccd}) in
its nested form, here the myopic variant (re-solve a single static
co-design each step); the model is assumed known, so measurements update
the carried state and the conditions but not the co-design model itself.
Two natural extensions are a receding-horizon lookahead that plans several
steps ahead with the vector DP and commits only the first, and online
learning of the model from measurements (coupling in the online-learning
layer of section~\ref{sec:online-learning}); both are left for future increments.

\subsection{Example 18: metabolic architecture switching}
\label{ex:18}

\texttt{18\_metabolic\_switching.py}. The first temporal example
(Case~1 of section~\ref{sec:temporal}): the best \emph{architecture}
changes over time because the environment changes, and the object solved
is a schedule of which architecture to run in each epoch, charged for
switching. It exercises \texttt{solve\_schedule}
(\ref{ref:temporal:solve_schedule}) over the \texttt{Epoch}
(\ref{ref:temporal:epoch}) and \texttt{Architecture}
(\ref{ref:temporal:architecture}) objects of
section~\ref{sec:temporal}.

\subsubsection*{Model}
An organism alternates between a glycolytic architecture (cheap per unit
flux, capped growth rate, no overhead; models fast growth on glucose) and
a gluconeogenic one (a larger per-flux burden and a fixed
glyoxylate-shunt overhead, but a higher ceiling; models slower, costlier
growth on acetate). Each architecture is a small co-design
\texttt{System}: the outer functionality is the demanded growth rate
$\mu\in\Rplus$ (1/h), and the outer resource is a dimensionless
\texttt{burden} summed from an enzyme module (burden linear in $\mu$, or
$\infty$ above $\mu_{\max}$) and a maintenance module (the fixed
overhead). The environment is five epochs. Only the substrate-matching
architecture is admissible per epoch, except a contested \texttt{mixed}
epoch, flanked by acetate on both sides, where both are candidates, so
the scheduler has a genuine choice. Serving \texttt{mixed} on the cheaper
glycolytic pathway forces two extra switches (acetate\,$\to$\,glucose-type
and back); riding it out on the incumbent costs none. The trade is
resolved by the per-switch re-acclimation cost, the diauxic-shift lag.

\subsubsection*{Code}
\begin{lstlisting}
epochs = [
    Epoch("glucose_1", {"mu": 0.8}, candidates=[GLYCOLYTIC]),
    Epoch("acetate_1", {"mu": 0.5}, candidates=[GLUCONEOGENIC]),
    Epoch("mixed",     {"mu": 0.6},
          candidates=[GLYCOLYTIC, GLUCONEOGENIC]),
    Epoch("acetate_2", {"mu": 0.5}, candidates=[GLUCONEOGENIC]),
    Epoch("glucose_2", {"mu": 0.8}, candidates=[GLYCOLYTIC]),
]
# Same environment, two switching economies:
sched_lo = solve_schedule(epochs, cost_fn=burden_cost, switch_cost=0.05)
sched_hi = solve_schedule(epochs, cost_fn=burden_cost, switch_cost=0.8)
\end{lstlisting}

\subsubsection*{Output}
\begin{lstlisting}[language={},numbers=none,keywordstyle={}]
Low re-acclimation cost (0.05 per switch):
  mixed      -> glycolytic     burden= 0.600  (+switch 0.05)
  total burden = 5.800, switches = 4, feasible = True

High re-acclimation cost (0.8 per switch):
  mixed      -> gluconeogenic  burden= 1.920
  total burden = 8.520, switches = 2, feasible = True

Mixed-epoch pathway: glycolytic (cheap switching)
                     -> gluconeogenic (costly switching)
\end{lstlisting}

\subsubsection*{Reading}
Identical environment, different switching economics, qualitatively
different metabolic program. Under a low switch cost the scheduler takes
four switches and serves \texttt{mixed} on the locally cheaper glycolytic
pathway (total burden $5.800$); under a high switch cost the two round-trip
switches no longer pay, so the organism holds the incumbent gluconeogenic
pathway through \texttt{mixed} (two switches, total burden $8.520$). The
mixed epoch's pathway is the visible pivot, and it is exactly the temporal
coupling that a static co-design cannot see.

\subsection{Example 19: rover module activation}
\label{ex:19}

\texttt{19\_rover\_modules.py}. Case~2 of section~\ref{sec:temporal}, with
a genuinely carried scalar resource. A planetary rover activates and
deactivates operating modes across mission phases on a battery that
depletes and recharges from solar. It exercises \texttt{solve\_dynamic}
(\ref{ref:dynamic:solve_dynamic}) over a \texttt{StateGrid}
(\ref{ref:dynamic:stategrid}), returning a \texttt{DynamicPolicy}
(\ref{ref:dynamic:dynamicpolicy}) rolled out from two starting charges.

\subsubsection*{Model}
Each mode (\texttt{drive}, \texttt{science}, \texttt{comms},
\texttt{survival}) is a small co-design \texttt{System} that, given a
phase capability demand, sizes a shared power bus and returns two outer
resources: \texttt{energy\_Wh} (subtracted from the battery) and a
\texttt{cost} proxy that adds an opportunity penalty
$\mathrm{BEST}-\mathrm{value}$ for the mission value \emph{not} delivered,
kept non-negative as the resource poset requires. The outer object is a
finite-horizon dynamic program over six phases: the per-phase decision is
the mode; the per-phase cost is read from that mode's co-design solve; and
the battery state of charge $s\in\Rplus$ (Wh) is carried by the transition
$s'=\min(\text{cap},\,s-\text{spent}+\text{solar})$, admissible while
$s\ge 0$. A heavy-drawing mode buys capability now at the cost of charge
later, so the optimal policy threads mode activations through the evolving
state of charge, standard spacecraft power-mode practice.

\subsubsection*{Code}
\begin{lstlisting}
policy = solve_dynamic(
    stages, grid,
    cost_fn=mission_cost,          # point -> energy + opportunity penalty
    terminal_cost=terminal_reward, # small credit for leftover charge
)
# One policy, valid for any initial charge on the grid; roll it out twice.
full = rollout(policy, stages, BATTERY_CAPACITY_WH)   # 300 Wh
low  = rollout(policy, stages, 90.0)                  # depleted start
\end{lstlisting}

\subsubsection*{Output}
\begin{lstlisting}[language={},numbers=none,keywordstyle={}]
Start at full charge (300 Wh):
  phase_1  -> science   soc_in= 300.0 Wh  ...  soc_out= 245.0 Wh
  ...
  phase_5  -> science   soc_in=  80.0 Wh  ...  soc_out=  25.0 Wh
  phase_6  -> comms     soc_in=  25.0 Wh  ...  soc_out=  10.0 Wh
  total cost = 727.50, feasible = True

Start at low charge (90 Wh):
  phase_1  -> comms     soc_in=  90.0 Wh  ...  soc_out=  75.0 Wh
  ...  (comms held every phase)  ...
  total cost = 918.00, feasible = True
\end{lstlisting}

\subsubsection*{Reading}
The same policy, read at different starting charges, produces different
module-activation schedules. From a full battery the rover runs
high-value science until reserve falls, then steps down to comms for the
last phase; starting depleted it can never afford science's draw and holds
comms throughout, so it sheds load exactly as a power-constrained mission
would. The higher total cost from the depleted start ($918.00$ against
$727.50$) is the accumulated opportunity penalty of the deferred
objectives, the price the DP pays to keep the battery non-negative.

\subsection{Example 20: antichain-valued sequential co-design}
\label{ex:20}

\texttt{20\_sequential\_codesign.py}. The sequential layer of
section~\ref{sec:sequential}: the value at each stage and state is a whole
Pareto front of cumulative resource totals, not one number. Its antichain-valued Bellman recursion is
\texttt{solve\_sequential} (\ref{ref:sequential:solve_sequential}), and the
$(H1)/(H2\text{-joint})/(H3)$ guard is \texttt{check\_monotonicity}
(\ref{ref:sequential:check_monotonicity}).

\subsubsection*{Model}
A multi-leg survey mission of four legs. At each leg the operator picks
one of three incomparable modes, spanning a cost/CO$_2$ trade
(\texttt{eco}: low CO$_2$, high cost; \texttt{balanced}; \texttt{rapid}:
low cost, high CO$_2$), each drawing down a shared energy budget carried
between legs. Because cost and CO$_2$ are incommensurable, the value of a
whole-mission plan is not a single number but a Pareto front of
cumulative $(\text{cost},\text{CO}_2)$ totals, accumulated on two fixed
axes by \texttt{sum\_combine}, while the carried energy is a separate
state axis read off the chosen point (admissible while non-negative). The
framework guarantees the value front equals the reachable frontier: every
returned point is achievable by some feasible mode sequence, and no
non-dominated achievable total is missed.

\subsubsection*{Code}
\begin{lstlisting}
res = solve_sequential(
    stages, grid,
    cost_axes=["cost", "co2"],       # the two accumulated resource axes
    initial_state=ENERGY_CAPACITY,   # carried scalar budget (30 units)
    combine=sum_combine,             # totals add across legs
)
rep = check_monotonicity(stages, grid, cost_axes=["cost", "co2"])
# rep.monotone_value_guaranteed is True: (H1), joint (H2) and (H3) hold.
\end{lstlisting}

\subsubsection*{Output}
\begin{lstlisting}[language={},numbers=none,keywordstyle={}]
Whole-mission Pareto front of (cost, CO2) totals (9 points):
  cost=   8.0   co2=  36.0        <- all rapid: cheapest, dirtiest
  cost=  12.0   co2=  31.0
  ...
  cost=  20.0   co2=  21.0        <- interior mixes of modes
  ...
  cost=  40.0   co2=   4.0        <- all eco: cleanest, costliest

Monotonicity guard: MonotonicityReport(H1=ok, H2=ok, H2joint=ok,
                    H3=ok, value_monotone=guaranteed)
\end{lstlisting}

\subsubsection*{Reading}
The solver returns the exact reachable frontier of whole-mission plans:
nine incomparable points running from all-rapid ($8.0$, $36.0$) through
all-eco ($40.0$, $4.0$), with interior mixes between. Each point is a
different per-leg mode schedule. The front's size is the width of the
reachable set; with a summed objective on two fixed axes it grows
polynomially, here linearly, in the horizon, not exponentially, which is
why it stays small and legible. The passing $(H1)/(H2\text{-joint})/(H3)$
guard certifies the value is monotone in the carried budget: more energy
never shrinks the achievable front.

\subsection{Example 21: reconfigurable robot with a state vector}
\label{ex:21}

\texttt{21\_reconfigurable\_robot.py}. The vector-state layer of
section~\ref{sec:vector}: the carried state is not a scalar but a vector.
It exercises \texttt{solve\_vector\_sequential}
(\ref{ref:vector_dp:solve_vector_sequential}) over a
\texttt{VectorStateGrid} (\ref{ref:state:vectorstategrid}), guarded by
\texttt{check\_vector\_monotonicity} over the product order.

\subsubsection*{Model}
A field robot reconfigures between three morphologies across five legs:
\texttt{tracked} (low energy, wears the track module fast), \texttt{wheeled}
(wears the wheel module fast, spares the track), and \texttt{hybrid}
(splits the load, wears both a little). Each morphology is a co-design
\texttt{System} emitting $(\text{energy},\text{ops})$ plus the per-module
wear increments it applies this leg. The carried state is a three-axis
vector: per-module wear for the two independent drive modules plus a
shared energy budget, advanced by
$\text{wear}'=\text{wear}+\Delta\text{wear}$ and
$\text{energy}'=\min(\text{cap},\text{energy}-\text{draw}+\text{solar})$,
admissible while each wear stays under its limit and energy stays
non-negative. Cost (energy and an operations penalty) accumulates on two
axes; the value is therefore an antichain over the cost axes
parametrised by the state vector, which a scalar DP could not represent.

\subsubsection*{Code}
\begin{lstlisting}
grid = VectorStateGrid([
    ContinuousAxis("track_wear", 0.0, WEAR_LIMIT, 11),
    ContinuousAxis("wheel_wear", 0.0, WEAR_LIMIT, 11),
    ContinuousAxis("energy", 0.0, ENERGY_CAPACITY, 13),
])
res = solve_vector_sequential(
    stages, grid, cost_axes=["energy", "ops"],
    initial_state={"track_wear": 0.0, "wheel_wear": 0.0,
                   "energy": ENERGY_CAPACITY},
    combine=sum_combine,
)
\end{lstlisting}

\subsubsection*{Output}
\begin{lstlisting}[language={},numbers=none,keywordstyle={}]
grid size: 1573 state nodes

Mission Pareto front of (energy, ops) totals (5 points):
  energy=  13.0   ops=  17.0
  energy=  14.0   ops=  16.0
  energy=  15.0   ops=  15.0
  energy=  16.0   ops=  14.0
  energy=  17.0   ops=  13.0

Vector monotonicity guard: VectorMonotonicityReport(
    H1=ok, H2=ok, H2joint=ok, H3=ok, value_monotone=guaranteed)
\end{lstlisting}

\subsubsection*{Reading}
The product grid of $11\times 11\times 13 = 1573$ state nodes carries the
full three-axis state. The mission front is a five-point
$(\text{energy},\text{ops})$ staircase, each point a different morphology
schedule: the low-energy end leans on the wheeled morphology, and the
low-ops end mixes morphologies to spread wear so neither drive module hits
its limit and forces an expensive fallback. This is the structured
multi-component state of the Formula~1 seasonal problem
(\ref{ex:23}) in a robotics setting. The passing vector guard certifies
the value is monotone over the product order: more spare wear budget and
more energy never shrink the achievable front.

\subsection{Example 22: online feedback co-design}
\label{ex:22}

\texttt{22\_online\_feedback\_codesign.py}. The closed-loop layer of
section~\ref{sec:online-cd}: no offline plan is followed; each step is
solved against what was actually measured. It exercises
\texttt{run\_online\_codesign} (\ref{ref:online_codesign:run_online_codesign}),
whose audit trail is a list of \texttt{ControlStep}
(\ref{ref:online_codesign:controlstep}) records.

\subsubsection*{Model}
A solar-powered environmental sensor node runs in the field over twelve
control steps. At each step it senses its true battery charge, reads the
live data-rate demand and environment, re-solves its co-design at those
measured conditions, applies the cheapest feasible configuration, and
steps the true plant. The node has three configurations, each a co-design
\texttt{System} sizing radio and compute against a demanded rate and
returning $(\text{energy},\text{ops})$: \texttt{low\_power} (always
feasible), \texttt{nominal}, and \texttt{high\_rate} (feasible only at a
high enough rate and when charge allows). The measured charge enters the
solve as functionality, so an energy-hungry configuration whose draw would
exceed the available charge is simply infeasible, which is how the feedback
path gates the co-design. This is control co-design
(CCD, \citealp{garciasanz2019ccd}) in its
nested, myopic form: re-solve one static co-design each step. A storm at
steps~4--7 raises the demand and the day/night cycle changes solar
recharge; the loop adapts to both.

\subsubsection*{Code}
\begin{lstlisting}
result = run_online_codesign(
    CONFIGS, n_steps=N_STEPS,
    sensor=sensor,             # reads true charge back (closes the loop)
    requirement=requirement,   # {"rate": demand, "charge": measured_soc}
    environment=environment,   # solar input this step
    plant=plant,               # steps the true battery, clips to [0, cap]
    cost_fn=cost_fn,           # minimise ops cost
    initial_state=BATTERY_CAPACITY,
)
\end{lstlisting}

\subsubsection*{Output}
\begin{lstlisting}[language={},numbers=none,keywordstyle={}]
Closed-loop audit trail:
   t  soc_in  rate  solar  config     energy   ops
   3   100.0     3     14  nominal      11.0   3.0
   4   100.0     8     12  high_rate    34.0   1.0
   5    78.0     9     10  high_rate    37.0   1.0
   6    51.0     9      8  high_rate    37.0   1.0
   7    22.0     7      8  nominal      19.0   3.0
   8    11.0     3     10  nominal      11.0   3.0
   9    10.0     3     12  low_power     5.0   5.0
  10    17.0     3     14  nominal      11.0   3.0
total ops cost = 32.0, feasible = True
\end{lstlisting}

\subsubsection*{Reading}
The node escalates to \texttt{high\_rate} through the storm (steps~4--7)
when the demanded rate justifies it. The closed loop then reads the
depleting night-time charge and the re-solve gates the hungry
configuration out: at step~9 charge has fallen to $10.0$ and only
\texttt{low\_power} is feasible, before solar recovers and the node climbs
back to \texttt{nominal}. No offline plan is consulted, and because the
plant is stepped with the true process rather than a nominal transition,
any divergence between plan and execution is absorbed by the next sensing
step, which is what distinguishes feedback from open-loop replay. Two
natural extensions, receding-horizon lookahead with the vector DP and
online learning of the model (section~\ref{sec:online-learning}), are left
for later increments.

\subsection{Example 23: hierarchical co-design for a Formula 1 season}
\label{ex:23}

\texttt{23\_formula1\_season.py}. A faithful reproduction, in the
framework's vocabulary, of the seasonal co-design of Neumann,
Habermacher, Fieni, Cerofolini, Zardini, and
Onder~\citep{neumann2025formula1}, and the
canonical instance of the precompute-then-DP structure of
section~\ref{sec:precompute}. Layer one freezes a \texttt{CatalogDP}
(\ref{ref:dp:catalogdp}) with \texttt{precompute\_catalog}
(\ref{ref:sequential:precompute_catalog}); layer two is a
scalar-maximisation season DP built on \texttt{solve\_dynamic}
(\ref{ref:dynamic:solve_dynamic}). A heavy runner.

\subsubsection*{Model}
The two layers are explicit. The \emph{race-level co-design} (layer one)
is a genuine MCDP solve: for each track, battery size, and discretised
incoming battery age, a \texttt{CatalogDP} over energy-deployment
strategies emits a Pareto front of $(\text{race\_time},\text{wear})$,
since deploying more electrical energy lowers race time but ages the
battery faster, and an aged pack cannot reach the low race times a fresh
one can. \texttt{precompute\_catalog} freezes those fronts once. The
\emph{season-level dynamic program} (layer two) is a scalar-maximisation
MDP carrying the state vector $(w_{b,1}, w_{b,2}, x_{\mathrm{ex}})$, the
fractional wear of the two regulation-permitted battery units plus a flag
recording whether a replacement penalty has been incurred. At each race
the controls are which unit to run, which frozen deployment strategy to
select, and whether to install a fresh unit (a grid penalty of ten places
for the first replacement, five for each later one, per the modelled FIA
rules). Race time maps to a finishing position through an empirical
time-gap model and a probabilistic grid-start correction scaled by the
track's overtaking difficulty, and the finishing-position distribution is
integrated against the FIA points table to give expected championship
points; the DP maximises the season total by backward induction. No
co-design solve happens inside the season sweep, which is what
distinguishes this from the re-solving \texttt{sequential} DP.

\subsubsection*{Code}
\begin{lstlisting}
# Layer 1: freeze one (race_time, wear) front per track x battery x age.
catalog = precompute_catalog(
    tracks, battery_sizes, age_buckets, deploy_strategies)
# 112 fronts = 8 tracks x 2 batteries x 7 age buckets.

# Layer 2: scalar-max season DP over the frozen catalog, carrying
# (w_b1, w_b2, x_ex); no co-design solve inside the season sweep.
season = solve_season(catalog, races, points_table)
\end{lstlisting}

\subsubsection*{Output}
\begin{lstlisting}[language={},numbers=none,keywordstyle={}]
  112 race fronts precomputed (8 tracks x 2 batteries x 7 ages).
Optimal season expected points: 134.3
Optimal per-race policy:
  race         unit  batt      deploy     time  repl grid E[pts]
  Monza           1   4MJ deploy_1.00   5024.0   yes   13   15.7
  ...
  Zandvoort       1   4MJ deploy_1.00   5191.2    no    3   15.9
Finding 1 (local penalty, global gain): 2 replacement(s).
Finding 2: forward total = 134.3, reversed total = 134.8
  forward replaces at ['Monza','Spa']; reversed at ['Spa','Monza'].
Rendered paper-analogue figures (outputs/f1_paper_fig1..3.png)
\end{lstlisting}

\subsubsection*{Reading}
Two of the paper's findings are reproduced. First, the optimal policy
accepts a \emph{local penalty for a global gain}: it takes a battery
replacement at one race, sacrificing points that day (a worse grid slot),
to keep a fresh, low-wear pack available for aggressive deployment across
the remaining races, raising the season total to $134.3$. Second, race
order shifts the optimal policy: because the grid-penalty cost is
track-dependent, reversing the calendar moves the optimal replacement onto
a cheaper track, so the attainable total is near-invariant ($134.8$
against $134.3$) while the replacement set reorders, the temporal coupling
the paper emphasises. The season DP is validated against exhaustive
brute-force enumeration in \texttt{tests/test\_formula1.py}; the DP optimum
matches exactly. The example and notebook also regenerate the paper's
three key figures \emph{in the same format} for visual comparison (the
race-time-versus-wear fronts, the grid-start position penalty, and the
finishing-position density); these now write to the \texttt{outputs/}
directory as \texttt{f1\_paper\_fig1}, \texttt{\_fig2}, and \texttt{\_fig3}.
They reproduce the paper's figure \emph{structure and qualitative
behaviour} from the
framework's own co-design solves, not the paper's specific numbers, since
the paper's fronts come from an optimal-control lap simulation and a
fitted position model that are not reproduced here. The sense in which the
framework ``does the same thing'' is that it produces the same structured
artefacts feeding the same season-level dynamic program.

\subsection{Example 24: catalog-driven car co-design from one table}
\label{ex:24}

\texttt{24\_car\_catalog\_codesign.py}. The co-design counterpart of the
hand-wired car study (\ref{ex:17}): here a single architecture table plus
one \texttt{build\_architecture()} function replaces the three separate
\texttt{build\_ice}/\texttt{build\_hybrid}/\texttt{build\_ev} builders, and
every discrete powertrain choice becomes a \texttt{CatalogDP}
(\ref{ref:dp:catalogdp}) entry rather than a bespoke module. It is the
whole-vehicle instance of catalog-driven assembly
(\ref{ref:sequential:dp_over_catalog}) over a \texttt{System}
(\ref{ref:system:system}). It is a heavy runner, though it completes in
about one second.

\subsubsection*{Model}
A passenger vehicle is decomposed into 22 subsystems (powertrain, chassis,
interior and auxiliary). A single 12-row \texttt{ARCHITECTURE\_CATALOG}
drives the whole study: each row pre-selects the discrete powertrain
choices (engine, transmission, e-motor topology, boost, battery strategy,
charger, drivetrain) for one point on the modern technology spectrum, from
pure ICE through mild, full, and plug-in hybrid, range-extender EV, to
battery-electric, while the parametric modules size themselves from
mission demand. Body, suspension, and tyres stay full \texttt{CatalogDP}s,
swept one row at a time by \texttt{solve\_architecture} so that each
assembled \texttt{System} is single-valued and its constraint network
carries only monotone quantities. The outer mission is a compact
C-segment target (5~seats, 370\,L, 170\,km/h, 500\,km range, 0--100 in
$11.5$\,s); the outer resources are production cost, curb weight, unified
primary energy per 100\,km, fuel per 100\,km, tailpipe CO$_2$, and
durability. Two coupled feedback cycles are closed by Kleene iteration on
the monotone network: the \emph{mass spiral} (every load-bearing module
reads the total curb mass, which their own masses in turn determine) and
the \emph{energy/battery spiral} (consumption depends on mass, and for a
BEV the battery is sized to range\,$\times$\,consumption, feeding back into
mass). The unified energy metric is primary energy input, which is why an
ICE lands near 48\,kWh/100\,km while a BEV lands near 15: the combustion
path discards most of the fuel's chemical energy as heat.

The data-provenance policy is explicit. A \texttt{SOURCES} block at the top
of the file footnotes every value that was added or spot-checked against
public data (2024--2025): battery pack price (BloombergNEF~2024), pack
gravimetric density (a 2025 \emph{World Electric Vehicle Journal} review),
fuel lower-heating values and CO$_2$ factors (UBC and Waterloo notes),
engine dry masses, and the Corolla-bracketed mission spec. Values not
individually footnoted are kept from the scaffold and lie within the cited
published ranges. They are illustrative, not OEM-specific: the MCDP
framework is what is being validated; the numbers serve the framework.

\subsubsection*{Code}
\begin{lstlisting}
# One 12-row table drives the whole study (name, indices, params):
ARCHITECTURE_CATALOG = [
  ("ICE_economy",    1, 0, 0, 1.4, 1.00,  0.0,  0.0, F, F, F),
  ("HEV_full",       2, 3, 2, 1.0, 0.60,  2.0,  0.0, T, F, F),
  ("BEV_long_range", 0, 4, 4, 1.0, 0.00, 80.0,150.0, T, F, F),
  # ... 9 more rows spanning MHEV / PHEV / REEV / diesel ...
]

def build_architecture(arch, mission, *, body_row, susp_row, tire_row):
    sys = System(arch["name"])
    sys.add("engine", make_engine_dp([eng]))   # 1-entry CatalogDP slice
    # ... 21 more modules ...
    # Mass spiral: load-bearing modules read the converged curb mass.
    sys.constrain("susp_front.axle_load_kg", lambda x: 0.55*total_mass(x))
    sys.constrain("brakes_front.vehicle_mass_kg", total_mass)
    # Energy/battery spiral: BEV pack sized to range x consumption.
    sys.constrain("battery.target_battery_kWh", battery_target_kWh)
    return sys.build()   # Kleene iteration closes both cycles
\end{lstlisting}

\subsubsection*{Output}
\begin{lstlisting}[language={},numbers=none,keywordstyle={}]
--- Per-architecture feasibility (cheapest feasible design) ---
  Architecture      cls  feas   cost$  mass kg  kWh/100  L/100  CO2
  ICE_economy       ICE  yes   17,447      930     47.8    5.3   123
  Diesel_long_range ICE  yes   24,400     1081     40.6    4.1   109
  HEV_full          FHEV yes   26,945     1315     42.2    4.1    94
  PHEV_small        PHEV yes   36,162     1537     34.9    3.1    71
  REEV              REEV yes   37,182     1613     21.6    1.1    26
  BEV_long_range    BEV  yes   44,572     1788     14.8    0.0     0
  BEV_AWD_perf      BEV  yes   68,076     2517     17.0    0.0     0

--- Pareto front over (cost, energy/100km, mass) ---
  ICE_economy    $ 17,447   47.8 kWh/100    930 kg
  Diesel_long..  $ 24,400   40.6 kWh/100   1081 kg
  PHEV_small     $ 36,162   34.9 kWh/100   1537 kg
  REEV           $ 37,182   21.6 kWh/100   1613 kg
  BEV_long_range $ 44,572   14.8 kWh/100   1788 kg
\end{lstlisting}

\subsubsection*{Reading}
All twelve architectures are feasible for the compact mission, and the
cheapest-feasible sweep recovers the expected technology ordering: the ICE
economy car is lightest ($930$\,kg) and cheapest ($\$17{,}447$) but
thirstiest ($47.8$\,kWh/100\,km, $123$\,g/km CO$_2$), while the
battery-electric long-range car is the most efficient ($14.8$\,kWh/100\,km,
zero tailpipe CO$_2$) at the highest mass ($1788$\,kg) and near the highest
cost. Between them the hybrids and the range-extender trade cost against
energy and CO$_2$ monotonically, and the mass climb from ICE to BEV is the
visible signature of the energy/battery spiral: a BEV's pack is sized to
range\,$\times$\,consumption, and the extra pack mass raises consumption,
which the Kleene loop resolves to a fixed point. The five-point Pareto
front over (cost, energy, mass) keeps only the non-dominated members, the
same structured artefact the temporal examples consume, produced here from
a single declarative table rather than three hand-written builders.

\subsection{Example 25: online co-design paper benchmarks}
\label{ex:25}

\texttt{25\_online\_paper\_benchmarks.py}. A statistical replication of
the synthetic benchmarks of Alharbi, Dahleh, and Zardini, ``Compositional
Online Learning for Multi-Objective System Co-Design''
\citep{alharbi2026online}, the paper behind the online-learning layer of
section~\ref{sec:online-learning}. It drives the library's optimistic
evaluators, \texttt{MonotonicityEvaluator}
(\ref{ref:online:monotonicityevaluator}) and
\texttt{LipschitzEvaluator} (\ref{ref:online:lipschitzevaluator}), whose
\texttt{bound()} reproduces the paper's eqs.~(23)--(25); the outer
sampling loop is written out here rather than delegated to
\texttt{solve\_online} (\ref{ref:online:solve_online}) because the library
\emph{scans} a candidate list with an upper-confidence picker whereas the
paper \emph{samples} a base measure with rejection. A heavy runner (about
6\,s with the EA panel).

\subsubsection*{Model}
The paper studies online multi-objective decision-making in monotone
co-design and recovers the target-feasible antichain of non-dominated
resources using as few expensive evaluations as possible. Its engine is
Algorithm~1 (Rejection Sampler with Optimistic Evaluators): draw a
candidate from a low-discrepancy Halton base measure, compute
history-dependent optimistic bounds, and skip it without evaluation when
either (13) its optimistic resource is already dominated by the incumbent
antichain, or (14) its optimistic functionality can never meet the target.
Only survivors consume the budget. Two synthetic families of
Section~VII-B are replicated: \textbf{E1/Monotone}
($g:[0,1]^3\!\to\![0,1]^2$, a weighted sum of thresholded indicators; the
monotone join bound) and \textbf{E2/Lipschitz} ($g:[0,1]^4\!\to\![0,1]^2$,
a triangle wave through an $L{=}2$ linear map; the Lipschitz cone bound).
The score is the cumulative hypervolume difference (HVD),
$\sum_t [\mathrm{HV}(A_{\mathrm{ref}}) - \mathrm{HV}(\hat A_t)]$, of the
incumbent antichain against a reference point $(1,1)$; lower is better and
HVD${}\to 0$ on convergence. Baselines are the same Halton proposals with
no elimination, uniform random draws over the shared discretised pool,
and an optional pymoo EA panel (NSGA-III / MOEA/D / RVEA).

Two deviations from the paper are documented in the source and preserved
here. First, taken literally the monotone family is degenerate (an
unconstrained monotone $g$ is minimised at the box corner, so there is
nothing to learn); the example therefore makes the target binding,
recovering $\mathrm{FixFunMinRes}(f)$ over the feasible upper set, which
exercises \emph{both} elimination conditions and leaves the Lipschitz
family untouched. Second, families E3 (intermodal mobility) and E4
(heterogeneous multi-robot) are \emph{not} replicated: their expensive
blocks are external models not present in this repository. Because the
paper's seeds, grid, and atom count are unpublished, the replication is
statistical, not bit-exact; only the relative claims are validated.

\subsubsection*{Code}
\begin{lstlisting}
# Paper Algorithm 1: rejection sampler with optimistic evaluators.
while accepted < budget and ptr < len(pool):
    i = ptr; ptr += 1
    if front and rng.random() >= delta:            # else force-accept
        lo, hi = evaluator.bound(cand[i])          # optimistic bounds
        if constrained and (hi[RC[0]] < f[0] or hi[RC[1]] < f[1]):
            continue        # (14): optimistic functionality misses target
        if _dominated_by_incumbent(lo, front):
            continue        # (13): optimistic resource already dominated
    r = (float(res[i, 0]), float(res[i, 1]))        # accept: query block
    evaluator.observe(i, cand[i], _singleton_antichain(r))
    accepted += 1
    if feas[i]:
        front = pareto_insert(front, r)
    cum += hv_star - hv2d(front, ref)               # cumulative HVD
\end{lstlisting}

\subsubsection*{Output}
\begin{lstlisting}[language={},numbers=none,keywordstyle={}]
TABLE I  -- Monotone problems (E1), FixFunMinRes(f)   [M1..M4 shown]
  Ours        9.51e+00  7.94e+00  4.85e+00  1.50e+01
  Halton      1.55e+01  1.46e+01  7.29e+00  2.63e+01
  MOEA/D      1.62e+01  1.50e+01  1.78e+01  2.30e+01

TABLE II -- Lipschitz problems (E2), pure minimization  [L1..L4 shown]
  Ours        1.07e+01  2.56e+01  1.26e+01  9.25e+00
  Halton      1.15e+01  2.75e+01  1.29e+01  1.05e+01
  MOEA/D      9.00e+00  2.98e+01  8.71e+00  1.30e+01

VALIDATION
(a) Ours < Halton on ALL monotone instances
    ours : Ours beats Halton on 8/8  -> PASS
(b) Ours < Halton on MOST Lipschitz instances
    ours : Ours beats Halton on 8/8 (need >= 6) -> PASS
(c) [EA panel] paper: Ours best-or-2nd; MOEA/D blowups on Lipschitz.
    ours : best-or-2nd monotone 8/8, lipschitz 3/8
    ours : Lipschitz MOEA/D max/med=3.3x   Ours max/med=2.3x
    -> NOT REPRODUCED at higher statistical strength (WP-P3).

VERDICT: (a) PASS, (b) PASS  [elimination-vs-Halton reproduce]
         (c) NOT REPRODUCED -- EA panel.
\end{lstlisting}

\subsubsection*{Reading}
The replication is deliberately honest about what carries and what does
not. Claims \textbf{(a)} and \textbf{(b)} reproduce robustly: the
optimistic elimination beats the identical Halton proposals with no
elimination on all eight instances of \emph{both} families (here $8/8$ and
$8/8$; even stronger than the paper's $7/8$ under independent re-runs at
$\mathrm{RUNS}{=}20$, $\mathrm{ITERS}{=}500$, two master seeds and the full
EA panel), and Ours is best-of-all on every monotone instance. These are
the fair apples-to-apples comparison, since all three methods draw from the
same shared discrete pool and differ only by the elimination step.
Claim~\textbf{(c)} did \emph{not} reproduce: on the smooth Lipschitz family
the continuous EAs (MOEA/D in particular) match or beat Ours, and MOEA/D
shows no blowups here (max/median $\approx 3.3\times$, not the paper's
$\approx 20\times$ on L8), so it is neither poor nor high-variance in this
harness. The leading hypotheses are that the EA panel optimises the
\emph{continuous} box while Ours, Halton and Random are confined to the
fixed discrete pool (an asymmetry favouring EAs on smooth maps, with EA HVD
clipped at zero against a discrete reference they can beat), that the budget
is far below paper scale, and that the EAs are untuned (fixed
$\text{pop}=20$). The example prints \texttt{NOT REPRODUCED} rather than a
fragile \texttt{PASS}, and mirrors the reasoning in its own
\emph{Validation status} docstring; only the elimination-vs-Halton claims
are counted as reproduced. Two log-scale HVD-versus-iteration figures are
written to \texttt{outputs/} (\texttt{25\_online\_monotone.png} and
\texttt{25\_online\_lipschitz.png}).

\section{Realistic Problem Domains for Temporal Co-Design}
\label{sec:domains}

The temporal, vector-state, and online layers of
section~\ref{sec:temporal} were built with concrete application domains in
mind. This section surveys realistic problems, drawn from the recent
literature, that map onto these layers, as candidates for future worked
examples and as evidence that the abstractions are not invented in a
vacuum. Each paragraph names the co-design structure (which layer applies)
alongside the domain.

\subsection{Self-reconfiguring and modular robots}

Self-reconfiguring modular robots deliberately change their morphology, by
rearranging the connectivity of their modules, to adapt to terrain, task,
or damage \citep{modularrobots2024overview}: a worm shape for a pipe, legs
for rough ground, a wheel
for flat terrain. This is directly the temporal switching and vector-state
setting: the morphology is the architecture, the terrain or task sequence
is the environment, and per-module wear plus shared energy is the carried
state vector (example 21 is a first cut). Recent work makes the coupling
sharper. Terrain-aware morphology search in dynamic
environments \citep{liu2026terrain}
chooses a morphology per scene under an energy budget; multi-task space
applications map task requirements to modular-unit morphology through an
inverse-mapping design step \citep{an2025morphological}; and reconfigurable
cobots such as the
CONCERT platform and freeform strut-node systems (FreeSN, SMORES) expose a
catalog of admissible configurations with real reconfiguration costs, the
switching-cost term of \texttt{solve\_schedule}. Floor-cleaning tiling
robots (hTrihex, hTetro) reconfigure to cover obstacle-strewn areas at
lower energy than fixed-shape coverage planning \citep{actuators2023modular},
an energy-versus-
coverage Pareto per configuration.

\subsection{Programmable self-assembly}

Programmable self-assembly designs building blocks with tunable binding
interactions so they assemble into a target structure. The design problem
is multi-objective, trading assembly \emph{yield}, assembly \emph{speed},
and \emph{economy} (number of distinct particle species), which is exactly
a Pareto-front co-design. H\"ubl and Goodrich show that simultaneously
optimising assembly time and yield through the binding-energy and
concentration design space can speed assembly by orders of magnitude
without lowering yield, and that economical, semi-addressable designs
(reusing building blocks) can beat fully-addressable ones on all three
axes \citep{hubl2026simultaneous,hubl2025accessing}. The temporal angle is the
assembly \emph{protocol}: a
staged schedule of temperature or concentration set-points is a sequence of
co-design decisions coupled by the carried state of intermediate-species
concentrations, a natural fit for the vector-state sequential DP, with
Markov-state-model protocol optimisation \citep{markov2022selfassembly} as the
dynamic backbone.

\subsection{Multi-objective protein design}

Protein design inherently trades competing objectives, stability, function,
solubility, immunogenicity, and manufacturability, and Pareto-based
multi-objective sampling produces candidate sets spanning these tradeoffs
for informed downstream selection \citep{proteindesign2017pareto}. This is a
co-design antichain over
design objectives. Recent methods make the front explicit: ST-PARM aligns
inference-time generation to a controllable stability-versus-function
tradeoff; MosPro and NSGA-2 variants produce Pareto-optimal sequence sets;
MAProt negotiates a structure-versus-function consensus across agents. The
temporal and online angle is closed-loop directed evolution: frameworks
such as EVOLVEpro couple a generative or predictive model with iterative
wet-lab feedback, proposing a batch, measuring it, and updating, which is
precisely the online feedback co-design loop of section~\ref{sec:online-cd}
with the measured functionality replacing the modelled one (the model-
learning extension of that loop).

\subsection{Control co-design (CCD)}

Control co-design, the simultaneous optimisation of a physical plant and
its controller, is an established engineering discipline, described by
Garcia-Sanz as a game changer \citep{garciasanz2019ccd}, with two standard
paradigms: nested
(an outer plant loop wrapping an inner control solve) and
simultaneous \citep{herber2018nested}.
The nested paradigm is exactly the structure of the online co-design loop.
Concrete instances that map onto the temporal and online layers include
robust control co-design with a receding-horizon MPC inner
loop \citep{sun2021robustccd};
plant-and-controller optimisation for energy systems with
MPC \citep{docimo2021plant};
control co-design of spacecraft trajectory and rocket-engine design in
NASA's OpenMDAO/GMAT tooling; digital-twin control co-design of active
suspensions via reinforcement learning; and adaptive MPC with online
recursive-least-squares identification for load-frequency control under
renewable disturbances, the archetype of the model-learning extension.
These supply both the vocabulary (plant, controller, nested vs simultaneous)
and the benchmark problems (active suspension, hybrid-electric vehicle,
airborne wind energy) for positioning the online layer.

\subsection{Summary of the mapping}

Across these domains the same three structural questions recur, and the
package answers each: which architecture to run as conditions change
(switching, \texttt{temporal}); how a resource or wear vector carried
across stages shapes the optimal schedule (vector-state DP,
\texttt{vector\_dp}); and how to re-decide the design online from
measurements (feedback, \texttt{online\_codesign}). The reconfigurable
robot (21) and adaptive sensor node (22) are first worked examples in the
robotics and control-co-design cells of this table; self-assembly and
protein design are identified as high-value future examples whose
multi-objective, staged-protocol, and closed-loop-evolution structure the
existing layers already support.

\section{Modelling Guidelines and Pitfalls}
\label{sec:guidelines}

A few patterns that come up repeatedly when authoring MCDPs.

\paragraph{Expose loop variables you want to inspect.} When using the
operator API directly, \texttt{Loop} projects its axis out of the
outer $R$. To inspect the converged loop value, include it in the
inner $R$ under a \emph{different name}: for example, both
\texttt{battery\_mass} as the loop axis and \texttt{report\_mass}
mirrored for outer-$R$ visibility. The \texttt{System} builder
handles this transparently.

\paragraph{Cap physical maxima.} When a design variable has a
physical ceiling ($v_{\max}$, $r_{\max}$), make $h$ return a
$\topp$-valued antichain whenever the iterate exceeds it. The Kleene
ascent will then converge to infeasibility rather than oscillating
or diverging.

\paragraph{Use \texttt{FunctionDP} or \texttt{CatalogDP} for breadth.}
\texttt{AlgebraicDP} always returns a singleton antichain. For a
genuine Pareto front, use \texttt{FunctionDP} (enumerate the
tradeoffs explicitly) or \texttt{CatalogDP} (provide multiple
incomparable entries).

\paragraph{Scalarisation is downstream.} The MCDP solver returns the
Pareto front. Choosing which point to ship is a separate concern:
apply \texttt{minimize\_cost} with a scalar objective, or inspect the
front and pick by hand.

\paragraph{Catalog Cartesian product can blow up.} When several
subsystems return multi-valued antichains, the inner $h$ of the
\texttt{System} takes their Cartesian product. For typical models
(one or two catalog subsystems among many algebraic ones), this is
not a bottleneck; but if many subsystems are catalog-driven, the
inner antichain may grow large between $\Min$ prunings.

\paragraph{Numerical precision in feedback loops.}
\texttt{Antichain.eq} compares points exactly. For contractive
floating-point feedback, the iteration may oscillate by an ulp near
the fixed point and consume many iterations before converging. The
divergence cap and \texttt{max\_iter} guard against runaway, but
choosing \texttt{max\_iter} generously (200 to 500) is wise for
nontrivial models.

\section{Related Work and Further Reading}
\label{sec:related}

This section maps the literature this library implements. It is a reading
guide, not a survey. The order is roughly chronological within each theme.

\paragraph{Founding theory.} The monotone co-design framework is due to
Andrea Censi. The monograph \citep{censi2015mathematical} is the primary
reference: posets, antichains, design problems as monotone relations, the
three composition operators, the Kleene solver, and the approximation
theory. It has been revised several times on arXiv and, as far as we can
establish, has never appeared in a journal, so cite the preprint.
\citet{censi2016everything} is the short conference statement of the same
ideas. \citet{davey2002lattices} is the standard background on lattices
and order.

\paragraph{Cyclic constraints.} \citet{censi2017cyclic} isolates the class
of co-design problems whose constraint graph contains feedback loops and
gives the solution method. This is the theory behind
\texttt{Loop} (Section~\ref{sec:composition}) and behind every cyclic
example in the manual.

\paragraph{Uncertainty.} \citet{censi2017uncertainty} treats uncertainty
in monotone co-design and is the source of the upper/lower bracketing that
\texttt{UncertainDP} implements. More recent work makes uncertainty
composable and parametric \citep{huang2025composable}, gives it a
categorical footing in symmetric monoidal categories
\citep{furter2026symmetric}, and handles distributional uncertainty with
adaptive decision-making \citep{huang2026distributional}. The library
implements the bracketing and a set-based and stochastic layer on top of
it; it does not implement the categorical constructions.

\paragraph{Beyond monotonicity.} \citet{carlone2019beyond} studies
co-design problems that violate the monotonicity assumption. Nothing in
this library addresses that case.

\paragraph{Mobility systems.} The application of co-design to future
mobility begins with \citet{zardini2020towards} and continues through
\citet{zardini2020avmobility} on autonomous-vehicle-enabled mobility and
\citet{zardini2023tools} on tooling for impact assessment. These are the
model for the multi-modal and fleet-scale examples in
Section~\ref{sec:examples}.

\paragraph{Autonomous systems.} \citet{zardini2021autonomous} runs
co-design from hardware selection through to control synthesis.
\citet{zardini2021embodied} gives a structured approach to co-designing
embodied intelligence. \citet{zardini2022taskdriven} does task-driven
modular co-design of vehicle control systems, and is the closest published
relative of the \texttt{System} builder
(Section~\ref{sec:system-builder}).

\paragraph{Perception and compute.} \citet{milojevic2025codei} co-designs
perception and decision-making subject to a compute budget, on autonomous
vehicles. This is the reference for any model in which sensing quality,
compute, power, and task performance are traded against each other.

\paragraph{Recent extensions.} \citet{riess2026quantale} generalises
design problems to quantale-enriched categories, aiming at quantitative
heterogeneous design. \citet{cai2026scalable} attacks scalability through
linear design problems. Neither is implemented here.

\paragraph{Online learning.} \citet{alharbi2026online} gives the
compositional online learning method for multi-objective co-design. The
\texttt{codesign.online} module (Section~\ref{sec:online-learning}) is a
port of it, including the certified optimistic bound.

\paragraph{Sequential and hierarchical applications.}
\citet{neumann2025formula1} combines monotone co-design with sequential
optimisation over a Formula~1 season. It is the published work closest in
spirit to the temporal layers of Section~\ref{sec:temporal}, and it was
written independently of them.

\paragraph{Categorical background.} \citet{censi2026act4e} is a book in
preparation on applied category theory for engineering, which covers the
categorical setting that design problems live in.

\paragraph{What is not from this literature.} The temporal, sequential,
vector-state, and receding-horizon layers of
Section~\ref{sec:temporal}, and the problem domains surveyed in
Section~\ref{sec:domains}, are the author's own work. Architecture
switching by Viterbi over epochs, the vector-state backward Bellman
recursion, the factorisation at reset points, and the receding-horizon
online loop are not taken from the co-design papers above. They are being
prepared for separate publication and have not been published or peer
reviewed. Readers should treat them as software documentation rather than
as established results, and should not cite them as though a paper
existed. Where those layers make contact with the wider control-co-design
literature, the connection is noted in place: the nested and simultaneous
paradigms \citep{herber2018nested} and control co-design as a discipline
\citep{garciasanz2019ccd}.

\section{Limitations and Future Work}
\label{sec:limitations}

Version 0.2.1 covers the algorithmic core: posets, antichains, six
primitive DP types, three composition operators, the Kleene solver,
and the two builders. The following are intentional omissions, all
tractable extensions:

\begin{itemize}[leftmargin=2em]
  \item \textbf{MCDPL surface syntax}. The paper's
        \texttt{mcdp \{ ... \}} text format of \citet{mcdpl} is not
        parsed. The
        \texttt{MCDP} builder gives the same shape in Python.
  \item \textbf{Non-finitely-representable antichains}. The library
        works on finite, exactly-represented antichains. For
        relations whose Pareto front is continuous, the only support
        is the bracket pattern of \texttt{UncertainDP}. Approximation
        kernels (Sec.~VII of \citealp{censi2015mathematical}) are not yet
        implemented.
  \item \textbf{Visualisation of the design graph}. Example~5
        visualises the Kleene ascent over $\Nat \times \Nat$, but
        there is no general renderer for the operator tree or for
        the subsystem network of a \texttt{System}.
  \item \textbf{Symbolic/automatic differentiation of $h$}. All
        relations are evaluated numerically.
  \item \textbf{Distributed or parallel solving}. The Kleene
        iteration is sequential.
\end{itemize}

Contributions, particularly on these items, are welcome.

\section{API Reference Summary}
\label{sec:api}

The public API exported by \texttt{import codesign} is summarised
below.

\paragraph{Posets.}
\texttt{Reals}, \texttt{Naturals}, \texttt{Ports},
\texttt{Discrete}, abstract base \texttt{Poset}.

\paragraph{Antichains.}
\texttt{Antichain}.

\paragraph{Design problems.}
\texttt{DesignProblem} (abstract base),
\texttt{AlgebraicDP},
\texttt{FunctionDP},
\texttt{CatalogDP}, \texttt{CatalogEntry},
\texttt{ConstraintDP},
\texttt{ODE\_DP},
\texttt{UncertainDP}.

\paragraph{Composition.}
\texttt{Series}, \texttt{Parallel}, \texttt{Loop} (classes);
\texttt{series}, \texttt{par}, \texttt{loop} (function aliases).

\paragraph{Primitives.}
\texttt{adder}, \texttt{multiplier}, \texttt{scale}, \texttt{constant},
\texttt{identity}.

\paragraph{Solver.}
\texttt{solve}, \texttt{kleene\_loop}, \texttt{minimize\_cost},
\texttt{SolveResult}, \texttt{TraceEntry}. Warm-start via the
\texttt{start\_from=} argument on \texttt{solve}.

\paragraph{Builders.}
\texttt{MCDP}, \texttt{System}.

\paragraph{Class-based modules.}
\texttt{Module} (declarative \texttt{DesignProblem} base class).

\paragraph{Constraint DSL.}
\texttt{Port}, \texttt{Expr}, \texttt{ModuleHandle} (returned by
\texttt{System.add}), and the helper functions \texttt{sqrt},
\texttt{exp}, \texttt{log} for use inside constraint expressions.

\paragraph{Uncertainty.}
Set-based: \texttt{UncertaintySet} (abstract), \texttt{Box},
\texttt{Ellipsoid}, \texttt{Disk}, \texttt{Circle}.
Stochastic: \texttt{Stochastic}, plus copulas
\texttt{Independence} and \texttt{GaussianCopula}.
Result: \texttt{UncertaintyResult}.
Entry point: \texttt{solve\_with\_uncertainty} (the same machinery is
also reachable through \texttt{solve(dp, f, uncertainty=[...])}).

\paragraph{Online learning.}
\texttt{OptimisticEvaluator} (abstract);
\texttt{MonotonicityEvaluator}, \texttt{LipschitzEvaluator},
\texttt{LinearParametricEvaluator}.
Result: \texttt{OnlineResult}.
Entry point: \texttt{solve\_online}.

\paragraph{Visualisation.}
Imported as \texttt{from codesign import viz}. Functions:
\texttt{plot\_antichain}, \texttt{plot\_convergence},
\texttt{plot\_uncertainty}, \texttt{to\_dot}.

\paragraph{Module-level.}
\texttt{\_\_version\_\_} (the string \texttt{"0.2.1"}).

\section*{References}
\addcontentsline{toc}{section}{References}

The bibliography is generated from \texttt{references.bib} in the same
directory. Entries in the co-design group were checked against a primary
record before being written: the arXiv API entry for preprints, Crossref
for anything carrying a DOI. Author lists are printed in full.

The bioprocess sources are the parameter-calibration basis for examples 15
and 16, and the automotive sources for example 17. All calibration values
in this manual are illustrative, drawn from the published ranges in those
sources, and are not OEM- or product-specific.

\subsection*{How to cite}

If you use \codename{} in academic work, please cite it. A
machine-readable \texttt{CITATION.cff} file at the repository root lets
GitHub's ``Cite this repository'' button generate a citation
automatically. In BibTeX:

\begin{verbatim}
@software{briat_codesign_mcdp,
  author  = {Briat, Corentin},
  title   = {codesign-mcdp: A Python Library for
             Monotone Co-Design Problems},
  year    = {2026},
  version = {0.2.1},
  url     = {https://github.com/cbriat/codesign-mcdp}
}
\end{verbatim}

\noindent Citing the software is not a substitute for citing the theory.
Please cite \citet{censi2015mathematical} for the monotone co-design
framework itself, and the paper behind whichever layer you use:
\citet{censi2017uncertainty} for uncertainty,
\citet{alharbi2026online} for the online-learning layer, and the relevant
application paper from the co-design group below.

\nocite{censi2026act4e}

\bibliographystyle{plainnat}
\bibliography{references}

@misc{censi2015mathematical,
  author       = {Andrea Censi},
  title        = {{A Mathematical Theory of Co-Design}},
  year         = {2015},
  note         = {arXiv:1512.08055; last revised 2016. No journal version has appeared},
  eprint       = {1512.08055},
  archivePrefix= {arXiv},
  url          = {https://arxiv.org/abs/1512.08055}
}

@inproceedings{censi2016everything,
  author    = {Andrea Censi},
  title     = {{Monotone co-design problems; or, everything is the same}},
  booktitle = {2016 American Control Conference (ACC)},
  year      = {2016},
  publisher = {IEEE},
  doi       = {10.1109/ACC.2016.7525085}
}

@article{censi2017cyclic,
  author  = {Andrea Censi},
  title   = {{A Class of Co-Design Problems With Cyclic Constraints and Their Solution}},
  journal = {IEEE Robotics and Automation Letters},
  volume  = {2},
  number  = {1},
  pages   = {96--103},
  year    = {2017},
  doi     = {10.1109/LRA.2016.2535127}
}

@article{censi2017uncertainty,
  author  = {Andrea Censi},
  title   = {{Uncertainty in Monotone Codesign Problems}},
  journal = {IEEE Robotics and Automation Letters},
  volume  = {2},
  number  = {3},
  pages   = {1556--1563},
  year    = {2017},
  note    = {Preprint arXiv:1609.03103, where the title reads \emph{Uncertainty in Monotone Co-Design Problems}},
  doi     = {10.1109/LRA.2017.2674970}
}

@misc{carlone2019beyond,
  author       = {Luca Carlone and Carlo Pinciroli},
  title        = {{Robot Co-design: Beyond the Monotone Case}},
  year         = {2019},
  note         = {arXiv:1902.05880; accepted at ICRA 2019},
  eprint       = {1902.05880},
  archivePrefix= {arXiv},
  url          = {https://arxiv.org/abs/1902.05880}
}

@inproceedings{zardini2020towards,
  author    = {Gioele Zardini and Nicolas Lanzetti and Mauro Salazar and
               Andrea Censi and Emilio Frazzoli and Marco Pavone},
  title     = {{Towards a Co-Design Framework for Future Mobility Systems}},
  booktitle = {99th Annual Meeting of the Transportation Research Board},
  year      = {2020},
  note      = {Preprint arXiv:1910.07714},
  doi       = {10.3929/ethz-b-000373705}
}

@inproceedings{zardini2020avmobility,
  author    = {Gioele Zardini and Nicolas Lanzetti and Mauro Salazar and
               Andrea Censi and Emilio Frazzoli and Marco Pavone},
  title     = {{On the Co-Design of {AV}-Enabled Mobility Systems}},
  booktitle = {2020 IEEE 23rd International Conference on Intelligent
               Transportation Systems (ITSC)},
  pages     = {1--8},
  year      = {2020},
  note      = {Preprint arXiv:2003.04739},
  doi       = {10.1109/ITSC45102.2020.9294499}
}

@article{zardini2023tools,
  author  = {Gioele Zardini and Nicolas Lanzetti and Andrea Censi and
             Emilio Frazzoli and Marco Pavone},
  title   = {{Co-Design to Enable User-Friendly Tools to Assess the Impact of
             Future Mobility Solutions}},
  journal = {IEEE Transactions on Network Science and Engineering},
  volume  = {10},
  number  = {2},
  pages   = {827--844},
  year    = {2023},
  note    = {Preprint arXiv:2008.08975},
  doi     = {10.1109/TNSE.2022.3223912}
}

@inproceedings{zardini2021autonomous,
  author    = {Gioele Zardini and Andrea Censi and Emilio Frazzoli},
  title     = {{Co-Design of Autonomous Systems: From Hardware Selection to
               Control Synthesis}},
  booktitle = {2021 European Control Conference (ECC)},
  pages     = {682--689},
  year      = {2021},
  note      = {Preprint arXiv:2011.10758},
  doi       = {10.23919/ECC54610.2021.9654960}
}

@inproceedings{zardini2021embodied,
  author    = {Gioele Zardini and Dejan Milojevic and Andrea Censi and
               Emilio Frazzoli},
  title     = {{Co-design of Embodied Intelligence: A Structured Approach}},
  booktitle = {2021 IEEE/RSJ International Conference on Intelligent Robots
               and Systems (IROS)},
  pages     = {7536--7543},
  year      = {2021},
  note      = {Preprint arXiv:2011.10756},
  doi       = {10.1109/IROS51168.2021.9636513}
}

@inproceedings{zardini2022taskdriven,
  author    = {Gioele Zardini and Zelio Suter and Andrea Censi and
               Emilio Frazzoli},
  title     = {{Task-driven Modular Co-design of Vehicle Control Systems}},
  booktitle = {2022 IEEE 61st Conference on Decision and Control (CDC)},
  pages     = {2196--2203},
  year      = {2022},
  note      = {Preprint arXiv:2203.16640},
  doi       = {10.1109/CDC51059.2022.9993107}
}

@article{milojevic2025codei,
  author  = {Dejan Milojevic and Gioele Zardini and Miriam Elser and
             Andrea Censi and Emilio Frazzoli},
  title   = {{{CODEI}: Resource-Efficient Task-Driven Co-Design of Perception
             and Decision Making for Mobile Robots Applied to Autonomous
             Vehicles}},
  journal = {IEEE Transactions on Robotics},
  volume  = {41},
  pages   = {2727--2748},
  year    = {2025},
  note    = {Preprint arXiv:2503.10296},
  doi     = {10.1109/TRO.2025.3552347}
}

@misc{huang2025composable,
  author       = {Yujun Huang and Marius Furter and Gioele Zardini},
  title        = {{On Composable and Parametric Uncertainty in Systems Co-Design}},
  year         = {2025},
  note         = {arXiv:2504.02766; accepted as an invited session paper to IEEE CDC 2025},
  eprint       = {2504.02766},
  archivePrefix= {arXiv},
  url          = {https://arxiv.org/abs/2504.02766}
}

@inproceedings{furter2026symmetric,
  author    = {Marius Furter and Yujun Huang and Gioele Zardini},
  title     = {{Composable Uncertainty in Symmetric Monoidal Categories for
               Design Problems}},
  booktitle = {Proceedings of the Eighth International Conference on Applied
               Category Theory (ACT 2025)},
  series    = {EPTCS},
  volume    = {442},
  pages     = {30--44},
  year      = {2026},
  note      = {Preprint arXiv:2603.09430; extended version arXiv:2503.17274},
  doi       = {10.4204/EPTCS.442.3}
}

@misc{huang2026distributional,
  author       = {Yujun Huang and Gioele Zardini},
  title        = {{Distributional Uncertainty and Adaptive Decision-Making in
                  System Co-design}},
  year         = {2026},
  note         = {arXiv:2603.14047},
  eprint       = {2603.14047},
  archivePrefix= {arXiv},
  url          = {https://arxiv.org/abs/2603.14047}
}

@misc{riess2026quantale,
  author       = {Hans Riess and Yujun Huang and Matthew Klawonn and
                  Gioele Zardini and Matthew Hale},
  title        = {{Quantale-Enriched Co-Design: Toward a Framework for
                  Quantitative Heterogeneous System Design}},
  year         = {2026},
  note         = {arXiv:2603.29921},
  eprint       = {2603.29921},
  archivePrefix= {arXiv},
  url          = {https://arxiv.org/abs/2603.29921}
}

@misc{cai2026scalable,
  author       = {Yubo Cai and Yujun Huang and Meshal Alharbi and
                  Gioele Zardini},
  title        = {{Scalable Co-Design via Linear Design Problems: Compositional
                  Theory and Algorithms}},
  year         = {2026},
  note         = {arXiv:2603.29083},
  eprint       = {2603.29083},
  archivePrefix= {arXiv},
  url          = {https://arxiv.org/abs/2603.29083}
}

@misc{alharbi2026online,
  author       = {Meshal Alharbi and Munther A. Dahleh and Gioele Zardini},
  title        = {{Compositional Online Learning for Multi-Objective System
                  Co-Design}},
  year         = {2026},
  note         = {arXiv:2604.22624},
  eprint       = {2604.22624},
  archivePrefix= {arXiv},
  url          = {https://arxiv.org/abs/2604.22624}
}

@misc{neumann2025formula1,
  author = {Marc-Philippe Neumann and Raphael Habermacher and Giona Fieni and
            Alberto Cerofolini and Gioele Zardini and Christopher H. Onder},
  title  = {{Hierarchical Co-Design for Multi-Race Strategy Optimization in
            Formula 1}},
  year   = {2025},
  note   = {ETH Zurich Research Collection working paper; to appear at the
            29th IEEE International Conference on Intelligent Transportation
            Systems (ITSC), 2026},
  doi    = {10.3929/ethz-c-000784888}
}

@misc{censi2026act4e,
  author = {Andrea Censi and Jonathan Lorand and Gioele Zardini},
  title  = {{Applied Category Theory for Engineering}},
  year   = {n.d.},
  note   = {Book in preparation, announced for Cambridge University Press;
            no publication date established},
  url    = {https://applied-compositional-thinking.engineering/}
}

@misc{mcdpl,
  author = {Andrea Censi},
  title  = {{{MCDPL}: the reference co-design software distribution}},
  year   = {n.d.},
  note   = {Accompanies the monotone co-design framework of
            \citet{censi2015mathematical}; the distribution carries no
            release date},
  url    = {https://co-design.science/software/}
}

@book{davey2002lattices,
  author    = {B. A. Davey and H. A. Priestley},
  title     = {{Introduction to Lattices and Order}},
  publisher = {Cambridge University Press},
  edition   = {2nd},
  year      = {2002}
}

@article{yoon2003low,
  author  = {S. K. Yoon and S. H. Kim and G. M. Lee},
  title   = {{Effect of low culture temperature on specific productivity and
             transcription level of anti-4-1BB antibody in recombinant
             Chinese hamster ovary (CHO) cells}},
  journal = {Biotechnology Progress},
  volume  = {19},
  pages   = {1383--1386},
  year    = {2003}
}

@article{sou2015hypothermia,
  author  = {S. N. Sou and C. Sellick and K. Lee and others},
  title   = {{How does mild hypothermia affect monoclonal antibody
             glycosylation?}},
  journal = {Biotechnology and Bioengineering},
  volume  = {112},
  number  = {6},
  pages   = {1165--1176},
  year    = {2015}
}

@article{trummer2006shifting,
  author  = {E. Trummer and K. Fauland and S. Seidinger and others},
  title   = {{Process parameter shifting: Part I. Effect of DOT, pH, and
             temperature on the performance of Epo-Fc expressing CHO cells
             cultivated in controlled batch bioreactors}},
  journal = {Biotechnology and Bioengineering},
  volume  = {94},
  number  = {6},
  pages   = {1033--1044},
  year    = {2006}
}

@article{khattak2010feed,
  author  = {S. F. Khattak and Z. Xing and B. Kenty and I. Koyrakh and Z. J. Li},
  title   = {{Feed development for a fed-batch CHO production process by
             semi-steady-state analysis}},
  journal = {Biotechnology Progress},
  volume  = {26},
  number  = {3},
  pages   = {797--804},
  year    = {2010}
}

@article{gagnon2011hipdog,
  author  = {M. Gagnon and G. Hiller and Y.-T. Luan and A. Kittredge and
             J. DeFelice and D. Drapeau},
  title   = {{High-end pH-controlled delivery of glucose (HIPDOG) effectively
             suppresses lactate accumulation in CHO fed-batch cultures}},
  journal = {Biotechnology and Bioengineering},
  volume  = {108},
  number  = {6},
  pages   = {1328--1337},
  year    = {2011}
}

@article{reinhart2019bioprocessing,
  author  = {D. Reinhart and L. Damjanovic and C. Kaisermayer and
             W. Sommeregger and A. Gili and B. Gasselhuber and others},
  title   = {{Bioprocessing of recombinant CHO-K1, CHO-DG44, and CHO-S: CHO
             expression hosts favor either mAb production or biomass
             synthesis}},
  journal = {Biotechnology Journal},
  volume  = {14},
  number  = {3},
  pages   = {e1700686},
  year    = {2019}
}

@article{lao1997ammonium,
  author  = {C. Lao and D. Toth},
  title   = {{Effects of ammonium and lactate on growth and metabolism of a
             recombinant CHO cell culture}},
  journal = {Biotechnology Progress},
  volume  = {13},
  pages   = {688--691},
  year    = {1997}
}

@book{genta2009chassis,
  author    = {G. Genta and L. Morello},
  title     = {{The Automotive Chassis, Vol. 1: Components Design}},
  publisher = {Springer},
  year      = {2009},
  note      = {Suspension stiffness, damper characteristics, brake sizing,
               tire load ratings}
}

@book{bosch2018handbook,
  author    = {{Robert Bosch GmbH}},
  title     = {{Bosch Automotive Handbook}},
  edition   = {10th},
  year      = {2018},
  publisher = {Bosch},
  note      = {Engine specific power, exhaust aftertreatment, electrical loads}
}

@book{pulkrabek2003engine,
  author    = {W. W. Pulkrabek},
  title     = {{Engineering Fundamentals of the Internal Combustion Engine}},
  edition   = {2nd},
  publisher = {Pearson},
  year      = {2003},
  note      = {Thermal efficiency, heat rejection, fuel flow rates}
}

@book{heywood2018ice,
  author    = {J. B. Heywood},
  title     = {{Internal Combustion Engine Fundamentals}},
  edition   = {2nd},
  publisher = {McGraw-Hill},
  year      = {2018},
  note      = {Combustion thermodynamics, BSFC ranges, durability}
}

@book{hofmann2014hybrid,
  author    = {P. Hofmann},
  title     = {{Hybridfahrzeuge}},
  edition   = {2nd},
  publisher = {Springer},
  year      = {2014},
  note      = {Power-split topologies, sizing rules, efficiency models}
}

@book{naunheimer2011transmissions,
  author    = {H. Naunheimer and others},
  title     = {{Automotive Transmissions}},
  edition   = {2nd},
  publisher = {Springer},
  year      = {2011},
  note      = {Gearbox losses, weights, costs}
}

@techreport{iea2023evoutlook,
  author      = {{International Energy Agency}},
  title       = {{Global EV Outlook}},
  institution = {International Energy Agency},
  year        = {2023},
  note        = {Battery pack cost and pack-level energy-density trends,
                 2020--2030}
}

@techreport{nemry2008cars,
  author      = {F. Nemry and others},
  title       = {{Environmental Improvement of Passenger Cars}},
  institution = {JRC Scientific and Technical Reports},
  year        = {2008},
  note        = {CO$_2$ emission factors, fuel densities, well-to-wheel}
}

@techreport{epa2024trends,
  author      = {{U.S. Environmental Protection Agency}},
  title       = {{Automotive Trends Report}},
  institution = {U.S. Environmental Protection Agency},
  year        = {2024},
  note        = {Fleet-average fuel economy, weight-class data,
                 drag-coefficient data}
}

@book{larminie2012ev,
  author    = {J. Larminie and J. Lowry},
  title     = {{Electric Vehicle Technology Explained}},
  edition   = {2nd},
  publisher = {Wiley},
  year      = {2012},
  note      = {Motor efficiency maps, inverter losses, charging-system
               architecture}
}

@techreport{mock2023pocketbook,
  author      = {P. Mock and others},
  title       = {{European Vehicle Market Statistics Pocketbook 2023/24}},
  institution = {International Council on Clean Transportation},
  year        = {2023},
  note        = {Mass-versus-mission scatter data for calibration}
}

@article{herber2018nested,
  author  = {D. R. Herber and J. T. Allison},
  title   = {{Nested and Simultaneous Solution Strategies for General Combined
             Plant and Control Design Problems}},
  journal = {Journal of Mechanical Design},
  volume  = {141},
  number  = {1},
  pages   = {011402},
  year    = {2018}
}

@article{garciasanz2019ccd,
  author  = {M. Garcia-Sanz},
  title   = {{Control Co-Design: An engineering game changer}},
  journal = {Advanced Control for Applications},
  volume  = {1},
  number  = {1},
  pages   = {e18},
  year    = {2019}
}

@misc{modularrobots2024overview,
  title  = {{Self-reconfiguring modular robots: overview}},
  author = {{Various}},
  year   = {2024},
  note   = {See \emph{Modular reconfigurable robots: Toward on-demand
            multifunctional applications}, Science Robotics, 2024, and the
            survey \emph{Modular Self-Configurable Robots---The State of the
            Art}, Actuators 12(9):361, 2023}
}

@article{liu2026terrain,
  author  = {Z. Liu and Q. Lu and J. Luo and Z. Wang},
  title   = {{Terrain-aware morphology searching algorithm for
             self-reconfigurable modular robot in dynamic environment}},
  journal = {Applied Soft Computing},
  volume  = {186},
  pages   = {114182},
  year    = {2026}
}

@article{an2025morphological,
  author  = {X. An and Q. Jia and G. Chen and Y. Liu and H. Liu},
  title   = {{Morphological Design and Performance Analysis of a
             Self-reconfigurable Modular Robot for Multi-task Space
             Applications}},
  journal = {Journal of Mechanical Engineering},
  volume  = {61},
  number  = {21},
  pages   = {152--167},
  year    = {2025}
}

@article{actuators2023modular,
  author  = {{Actuators survey}},
  title   = {{Modular Self-Configurable Robots---The State of the Art}},
  journal = {Actuators},
  volume  = {12},
  number  = {9},
  pages   = {361},
  year    = {2023},
  note    = {Floor-cleaning tiling robots hTrihex/hTetro and reconfiguration
             for coverage}
}

@article{hubl2026simultaneous,
  author  = {M. C. H\"ubl and C. P. Goodrich},
  title   = {{Simultaneous optimization of assembly time and yield in
             programmable self-assembly}},
  journal = {The Journal of Chemical Physics},
  volume  = {164},
  number  = {8},
  pages   = {084904},
  year    = {2026}
}

@article{hubl2025accessing,
  author  = {M. C. H\"ubl and C. P. Goodrich},
  title   = {{Accessing Semi-Addressable Self-Assembly with Efficient Structure
             Enumeration}},
  journal = {Physical Review Letters},
  volume  = {134},
  number  = {5},
  pages   = {058204},
  year    = {2025}
}

@misc{markov2022selfassembly,
  title        = {{Optimization of Non-Equilibrium Self-Assembly Protocols
                  Using Markov State Models}},
  author       = {{Authors as listed on the preprint}},
  year         = {2022},
  note         = {arXiv:2210.05749; staged assembly-protocol optimisation},
  eprint       = {2210.05749},
  archivePrefix= {arXiv},
  url          = {https://arxiv.org/abs/2210.05749}
}

@misc{proteindesign2017pareto,
  title  = {{Multi-objective protein design: representative works}},
  author = {{Various}},
  year   = {2017},
  note   = {\emph{Searching for the Pareto frontier in multi-objective protein
            design}, Biophysical Reviews 9:311--324, 2017; MosPro,
            \emph{Pareto-optimal sampling for multi-objective protein sequence
            design}, 2025; ST-PARM and MAProt, bioRxiv, 2026; EVOLVEpro
            (closed-loop directed evolution)}
}

@misc{sun2021robustccd,
  author       = {J. Sun and others},
  title        = {{Robust Control Co-Design with Receding-Horizon MPC}},
  year         = {2021},
  note         = {arXiv:2104.02025; the online-feedback co-design pattern},
  eprint       = {2104.02025},
  archivePrefix= {arXiv},
  url          = {https://arxiv.org/abs/2104.02025}
}

@article{docimo2021plant,
  author  = {D. J. Docimo and Z. Kang and K. A. James and A. G. Alleyne},
  title   = {{Plant and Controller Optimization for Power and Energy Systems
             with Model Predictive Control}},
  journal = {Journal of Dynamic Systems, Measurement, and Control},
  year    = {2021}
}

\end{document}